\PassOptionsToPackage{top=1cm, bottom=1cm, left=2cm, right=2cm}{geometry}
\documentclass[final,3p,times]{elsarticle}

\makeatletter
\def\ps@pprintTitle{%
   \let\@oddhead\@empty
   \let\@evenhead\@empty
   \def\@oddfoot{\reset@font\hfil\thepage\hfil}
   \let\@evenfoot\@oddfoot
}
\makeatother

\usepackage{amsthm, amssymb, amsfonts, mathrsfs, amsmath}
\everymath{\displaystyle}

\usepackage{palatino}

\usepackage[T5]{fontenc}

\usepackage{amscd}
\usepackage{graphicx}

\usepackage{pb-diagram}

\usepackage[x11names]{xcolor}
\usepackage[
  colorlinks,
]{hyperref}

\usepackage{longtable}
\usepackage{array}          
\renewcommand{\arraystretch}{1.15}
\allowdisplaybreaks[4]

\newcolumntype{L}{>{\raggedright\arraybackslash}p{0.23\textwidth}}

\newcommand\rurl[1]{%
  \href{http://#1}{\nolinkurl{#1}}%
}

\AtBeginDocument{\hypersetup{citecolor=blue, urlcolor=blue, linkcolor=blue}}

\theoremstyle{plain}
\newtheorem{thm}{Theorem}[subsection]

\newtheorem{corl}[thm]{Corollary}

\theoremstyle{definition}

\renewcommand{\theequation}{\arabic{section}.\arabic{equation}}

\newtheorem*{acknow}{Acknowledgments}

\theoremstyle{plain}
\newtheorem{therm}{Theorem}[section]
\newtheorem{conj}[therm]{Conjecture}

\newtheorem{corls}[therm]{Corollary}

\theoremstyle{definition}

\newtheorem{rems}[therm]{Remark}
\newtheorem{cy}[therm]{Note}
\newtheorem{notas}[therm]{Notation}

\theoremstyle{definition}

\newtheorem{nx}[thm]{Remark}

\everymath{\displaystyle}

\def\leq{\leqslant}
\def\geq{\geqslant}
\def\DD{D\kern-.7em\raise0.4ex\hbox{\char '55}\kern.33em}

\renewcommand{\baselinestretch}{1.2}

\usepackage{sectsty}
\subsectionfont{\fontsize{11.5}{11.5}\selectfont}

\makeatletter
\def\blfootnote{\xdef\@thefnmark{}\@footnotetext}
\makeatother

\begin{document}
\fontsize{11.5pt}{11.5}\selectfont

\begin{frontmatter}

\title{{\bf On the hit problem for the polynomial algebra and the algebraic transfer}}

\author{\DD\d{\u A}NG V\~O PH\'UC}
\address{{\fontsize{10pt}{10}\selectfont Department of Information Technology, FPT University, Quy Nhon A.I Campus,\\ An Phu Thinh New Urban Area, Quy Nhon City, Binh Dinh, Vietnam\\[2.5mm]  \textbf{\textit{(Dedicated to my beloved wife and cherished son)}}}\\[2.5mm]
\textit{\fontsize{11pt}{11}\selectfont \textbf{ORCID: \url{https://orcid.org/0000-0002-6885-3996}}}

}
\cortext[]{\href{mailto:dangphuc150488@gmail.com, phucdv14@fe.edu.vn}{\texttt{Email address: dangphuc150488@gmail.com, phucdv14@fpt.edu.vn}}}

\begin{abstract}
Let $\mathcal A$ be the classical, singly-graded Steenrod algebra over the prime order field $\mathbb F_2$ and let $P^{\otimes h}: = \mathbb F_2[t_1, \ldots, t_h]$ denote the polynomial algebra on $h$ generators, each of degree $1.$ Write $GL_h$ for the usual general linear group of rank $h$ over $\mathbb F_2.$ Then, $P^{\otimes h}$ is an $\mathcal A[GL_h]$-module. As is well known, for all homological degrees $h \geq 6$, the cohomology groups ${\rm Ext}_{\mathcal A}^{h, h+\bullet}(\mathbb F_2, \mathbb F_2)$ of the algebra $\mathcal A$ are still shrouded in mystery. The algebraic transfer $Tr_h^{\mathcal A}: (\mathbb F_2\otimes_{GL_h}{\rm Ann}_{\overline{\mathcal A}}[P^{\otimes h}]^{*})_{\bullet}\longrightarrow {\rm Ext}_{\mathcal A}^{h, h+\bullet}(\mathbb F_2, \mathbb F_2)$ of rank $h,$ constructed by W. Singer [Math. Z. \textbf{202} (1989), 493-523], is a beneficial technique for describing the Ext groups. Singer's conjecture about this transfer states that \textit{it is always a one-to-one map}. Despite significant effort, neither a complete proof nor a counterexample has been found to date.  The unresolved nature of the conjecture makes it an interesting topic of research in Algebraic topology in general and in  homotopy theory in particular.

\medskip

The objective of this paper is to investigate Singer's conjecture, with a focus on all $h\geq 1$ in degrees $n\leq 10 = 6(2^{0}-1) + 10\cdot 2^{0}$ and for $h=6$ in the general degree $n:=n_s=6(2^{s}-1) + 10\cdot 2^{s},\, s\geq 0.$ Our methodology relies on the hit problem techniques for the polynomial algebra $P^{\otimes h}$, which allows us to investigate the Singer conjecture in the specified degrees. Our work is a continuation of the work presented by Mothebe et al. [J. Math. Res. \textbf{8} (2016), 92-100] with regard to the hit problem for $P^{\otimes 6}$ in degree $n_s$, expanding upon their results and providing novel contributions to this subject. More generally, for $h\geq 6,$ we show that the dimension of the cohit module $\mathbb F_2\otimes_{\mathcal A}P^{\otimes h}$ in degrees $2^{s+4}-h$ is equal to the order of the factor group of $GL_{h-1}$ by the Borel subgroup $B_{h-1}$ for every $s\geq h-5.$ Especially, for the Galois field $\mathbb F_{q}$ ($q$ denoting the power of a prime number), based on Hai's recent work [C. R. Math. Acad. Sci. Paris \textbf{360} (2022), 1009-1026], we claim that the dimension of the space of the indecomposable elements of $\mathbb F_q[t_1, \ldots t_h]$ in general degree $q^{h-1}-h$ is equal to the order of the factor group of $GL_{h-1}(\mathbb F_q)$ by a subgroup of the Borel group $B_{h-1}(\mathbb F_q).$ As applications, we establish the dimension result for the cohit module $\mathbb F_2\otimes_{\mathcal A}P^{\otimes 7}$ in degrees $n_{s+5},\, s > 0.$ Simultaneously, we demonstrate that the non-zero elements $h_2^{2}g_1 = h_4Ph_2\in {\rm Ext}_{\mathcal A}^{6, 6+n_1}(\mathbb F_2, \mathbb F_2)$ and $D_2\in {\rm Ext}_{\mathcal A}^{6, 6+n_2}(\mathbb F_2, \mathbb F_2)$ do not belong to the image of the sixth Singer algebraic transfer, $Tr_6^{\mathcal A}.$ This discovery holds significant implications for Singer's conjecture concerning algebraic transfers. We further deliberate on the correlation between these conjectures and antecedent studies, thus furnishing a comprehensive analysis of their implications.

\end{abstract}



\end{frontmatter}

\tableofcontents

\section{Introduction}\label{s1}
\setcounter{equation}{0}

Throughout this text, we adopt the convention of working over the prime field $\mathbb F_2$ and denote the Steenrod algebra over this field by $\mathcal A$, unless otherwise stated. Let $V_h$ be the $h$-dimensional $\mathbb F_2$-vector space. It is well-known that the mod 2 cohomology of the classifying space $BV_h$ is given by $P^{\otimes h} := H^{*}(BV_h) = H^{*}((\mathbb RP^{\infty})^{\times h}) =  \mathbb F_2[t_1, \ldots, t_h],$ where $(t_1, \ldots, t_h)$  is a basis of $H^{1}(BV_h) = {\rm Hom}(V_h, \mathbb F_2).$ The polynomial algebra $P^{\otimes h}$ is a connected $\mathbb Z$-graded algebra, that is, $P^{\otimes h} = \{P_n^{\otimes h}\}_{n\geq 0},$ in which $P_n^{\otimes h}:= (P^{\otimes h})_n$ denotes the vector subspace of homogeneous polynomials of degree $n$ with $P_0^{\otimes h}\equiv \mathbb F_2$ and $P_n^{\otimes h} = \{0\}$ for all $n < 0.$ We know that the algebra $\mathcal A$ is generated by the Steenrod squares $Sq^{k}$ ($k\geq 0$) and subject to the Adem relations (see e.g., \cite{Walker-Wood}). The Steenrod squares are the cohomology operations satisfying the naturality property. Moreover, they commute with the suspension maps, and therefore, they are stable. In particular, these squares $Sq^k$ were applied to the vector fields on spheres and the Hopf invariant one problem, which asks for which $k$ there exist maps of Hopf invariant $\pm 1$. The action of $\mathcal A$ over  $P^{\otimes h} $ is determined by instability axioms. By Cartan's formula, it suffices to determine $Sq^{i}(t_j)$ and the instability axioms give $Sq^{1}(t_j)= t_j^{2}$ while $Sq^{k}(t_j) = 0$ if $k > 1.$ The investigation of the Steenrod operations and related problems has been undertaken by numerous authors (for instance \cite{Brunetti1, Singer1, Silverman, Monks, Wood}). An illustration of the importance of the Steenrod operators is their stability property, which, when used in conjunction with the Freudenthal suspension theorem \cite{Freudenthal}, enables us to make the claim that  the homotopy groups $\pi_{k+1}(\mathbb S^{k})$ are non-trivial for $k\geq 2$ (see also \cite{Phuc13} for further details).  The Steenrod algebra naturally acts on the cohomology ring $H^{*}(X)$ of a CW-complex $X.$ In several cases, the resulting $\mathcal A$-module structure on $H^{*}(X)$ provides additional information about $X$ (for instance the CW-complexes $\mathbb{C}P^4/\mathbb{C}P^2$ and $\mathbb{S}^6\vee \mathbb{S}^8$ have cohomology rings that agree as graded commutative $\mathbb F_2$-algebras, but are different as modules over $\mathcal A.$ We also refer the readers to \cite{Phuc10-0} for an explicit proof.)  Hence the Steenrod algebra is one of the important tools in Algebraic topology. Especially, its cohomology ${\rm Ext}_{\mathcal A}^{*, *}(\mathbb F_2, \mathbb F_2)$ is an algebraic object that serves as the input to the Adams (bigraded) spectral sequence \cite{J.A} and therefore, computing this cohomology is of fundamental importance to the study of the stable homotopy groups of spheres. 
\medskip

The identification of a minimal generating set for the $\mathcal A$-module $P^{\otimes h}$ has been a significant and challenging open problem in Algebraic topology in general and in homotopy theory in particular for several decades. This problem, famously known as the "hit" problem, was first proposed by Frank Peterson \cite{Peterson, Peterson2} through computations for cases where $h<2$ and has since captured the attention of numerous researchers in the field, as evidenced by works such as Kameko \cite{Kameko}, Repka and Selick \cite{Repka-Selick}, Singer \cite{Singer1}, Wood \cite{Wood}, Mothebe et al. \cite{Mothebe0, Mothebe, MKR}, Walker and Wood \cite{Walker-Wood, Walker-Wood2}, Sum \cite{Sum00, Sum1, Sum2, Sum2-0, Sum3, Sum4, Sum5}, Ph\'uc and Sum \cite{P.S1, P.S2}, Ph\'uc \cite{Phuc4, Phuc6, Phuc10, Phuc10-0, Phuc11, Phuc13, Phuc16}, Hai \cite{Hai}), and others. Peterson himself, as well as several works such as \cite{Priddy, Singer, Wood}, have shown that the hit problem is closely connected to some classical problems in homotopy theory.
To gain a deeper understanding of this problem and its numerous applications, readers are cordially invited to refer to the excellent volumes written by Walker and Wood \cite{Walker-Wood, Walker-Wood2}. An interesting fact is that, if $\mathbb F_2$ is a trivial $\mathcal A$-module, then the hit problem is essentially the problem of finding a monomial basis for the graded vector space $$QP^{\otimes h}=\big\{QP_n^{\otimes h} := (QP^{\otimes h})_n =P_n^{\otimes h}/\overline{\mathcal A}P_n^{\otimes h} =  (\mathbb F_2\otimes_{\mathcal A}P^{\otimes h})_n = {\rm Tor}_{0, n}^{\mathcal A}(\mathbb F_2, P^{\otimes h})\big\}_{n\geq 0}.$$ Here $\overline{\mathcal A}P_n^{\otimes h}:=P_n^{\otimes h}\cap \overline{\mathcal A}P^{\otimes h}$ and $\overline{\mathcal A}$ denotes the set of positive degree elements in $\mathcal A.$ 
The investigation into the structure of $QP^{\otimes h}$ has seen significant progress in recent years. Notably, it has been able to explicitly describe $QP^{\otimes h}$ for $h\leq 4$ and all $n > 0$ through the works of Peterson \cite{Peterson} for $h = 1, 2$, Kameko \cite{Kameko} for $h = 3$, and Sum \cite{Sum1, Sum2} for $h = 4$. Even so, current techniques have yet to full address the problem. 

While the information that follows may not be integral to the chief content of this paper, it will be beneficial for readers who desire a more in-depth comprehension of the hit problem.
When considering the field $\mathbb F_p$, where $p$ is an odd prime, one must address the challenge of the "hit problem," which arises in the polynomial algebra $\mathbb F_p[t_1, \ldots, t_h] = H^{*}((\mathbb CP^{\infty})^{\times h}; \mathbb F_p)$ on generators of degree $2$. This algebra is viewed as a module over the mod $p$ Steenrod algebra $\mathcal A_p$. Here $\mathbb CP^{\infty}$ denotes the infinite complex projective space. The action of $\mathcal A_p$ on $\mathbb F_p[t_1, \ldots, t_h]$ can be succinctly expressed by $\mathscr P^{p^{j}}(t_i^{r}) = \binom{r}{p^{j}}t_i^{r+p^{j+1}-p^{j}},\, \beta(t_i) = 0$ ($\beta\in \mathcal A_p$ being the Bockstein operator) and the usual Cartan formula. In particular, if we write $r = \sum_{j\geq 0}\alpha_j(r)p^{j}$ for the $p$-adic expansion of $r,$ then $\mathscr P^{p^{j}}(t_i^{r}) \neq 0$ if and only if $\binom{r}{p^{j}}\equiv \alpha_j(r)\, ({\rm mod}\, p)\, \neq 0.$ Since each Steenrod reduced power $\mathscr P^{j}$ is decomposable unless $j$ is a power of $p,$ a homogeneous polynomial $f$ is hit if and only if it can be represented as $\sum_{j \geq 0}\mathscr P^{p^{j}}(f_j)$ for some homogeneous polynomials $f_j\in \mathbb F_p[t_1, \ldots, t_h].$ In other words, $f$ belongs to $\overline{\mathcal A_p}\mathbb F_p[t_1, \ldots, t_h].$ (This is analogous to the widely recognized case when $p = 2.$) To illustrate, let us consider the monomial $t^{p(p+1)-1}\in \mathbb F_p[t].$ Since $\binom{2p-1}{p}\equiv \binom{p-1}{0}\binom{1}{1} \equiv 1\, ({\rm mod}\, p),$ $t^{p(p+1)-1} = \mathscr P^{p}(t^{2p-1}),$ i.e., $t^{p(p+1)-1}$ is hit. Actually, the hit problem for the algebra $ \mathbb F_p[t_1, \ldots, t_h]$ is an intermediate problem of identifying a minimal set of generators for the ring $H^{*}(V; \mathbb F_p) = H^{*}(BV; \mathbb F_p) = \Lambda(V^{\sharp})\otimes_{\mathbb F_p}S(V^{\sharp})$ as a module over $\mathcal A_p.$ Here $\Lambda(V^{\sharp})$ is an exterior algebra on generators of degree $1$ while $S(V^{\sharp})$ is a symmetric algebra on generators of degree $2.$ In both situations, the generators may be identified as a basis for $V^{\sharp}$, the linear dual of an elementary abelian $p$-group $V$ of rank $h$, which can be regarded as an $h$-dimensional vector space over the field $\mathbb F_p.$ Viewed as an algebra over the Steenrod algebra, $S(V^{\sharp})$ can be identified with $H^{*}((\mathbb CP^{\infty})^{\times h}; \mathbb F_p)$. Consequently, the cohomology of $V$ over the field $\mathbb F_p$ can be expressed as $\Lambda(V^{\sharp})\otimes_{\mathbb F_p} \mathbb F_p[t_1, \ldots, t_h].$ Thus, the information about the hit problem for $H^{*}(V; \mathbb F_p)$  as an $\mathcal A_p$-module can usuallly be obtained from the similar information about the hit problem for the $\mathcal A_p$-module $\mathbb F_p[t_1, \ldots, t_h]$ without much difficulty. With a monomial $f = t_1^{a_1}t_2^{a_2}\ldots t_h^{a_h}\in \mathbb F_p[t_1, \ldots, t_h],$ we denote its degree by $\deg(f) = \sum_{1\leq i\leq h}a_i.$ This coincides with the usual grading of $P^{\otimes h}$ for $p  = 2.$ Notwithstanding, it is one half of the usual grading of $\mathbb F_p[t_1, \ldots, t_h]$ for $p$ odd. With respect to this grading, Peterson's conjecture \cite{Peterson2} is no longer true for $p$ odd in general. As a case in point, our work \cite{Phuc14} provides a detailed proof that  $\alpha((i+1)p^{r}-1+1) = \alpha_r((i+1)p^{r}-1+1) =  i+1 > 1,$ but the monomials $t^{(i+1)p^{r}-1}\in \mathbb F_p[t],$ for $1\leq i< p-1,\, r\geq 0,$ are not hit.

Returning to the topic of the indecomposables $QP^{\otimes h}_n,$ let $\mu: \mathbb{N} \longrightarrow \mathbb{N}$ be defined by $\mu(n) = \min\bigg\{k \in \mathbb{N}: \alpha(n+k) \leq k\bigg\}$, where $\alpha(n)$ denotes the number of ones in the binary expansion of $n$. In the work of Sum \cite{Sum4}, it has been demonstrated that $\mu(n) = h$ if and only if there exists uniquely a sequence of integers $d_1 > d_2 > \cdots > d_{h-1}\geq d_h > 0$ such that $n  = \sum_{1\leq j\leq h}(2^{d_j} - 1).$ On the other side, according to Wood \cite{Wood}, if $\mu(n) > h,$ then $\dim QP^{\otimes h}_{n} = 0.$ This validates also Peterson's conjecture \cite{Peterson2} in general. 
Singer \cite{Singer1} later proved a generalization of Wood's result, identifying a larger class of hit monomials. In \cite{Silverman}, Silverman makes progress toward proving a conjecture of Singer which would identify yet another class of hit monomials. In \cite{Monks}, Monks extended Wood's result to determine a new family of hit polynomials in $P^{\otimes h}.$ Notably, Kameko \cite{Kameko} showed that if $\mu(2n+h) = h,$ then $QP^{\otimes h}_{2n+h}\cong QP^{\otimes h}_{n}.$ This occurrence elucidates that the surjective map 
$(\widetilde {Sq^0_*})_{2n+h}: QP^{\otimes h}_{2n+h}  \longrightarrow QP^{\otimes h}_{n},\ \mbox{[}u\mbox{]}\longmapsto \left\{\begin{array}{ll}
\mbox{[}y\mbox{]}& \text{if $u = \prod_{1\leq j\leq h}t_jy^{2}$},\\
0& \text{otherwise},
\end{array}\right.$ 
defined by Kameko himself, transforms into a bijective map when $\mu(2n+h) = h.$ Thus it is only necessary to calculate $\dim QP^{\otimes h}_n$ for degrees $n$ in the "generic" form:
\begin{equation}\label{pt}
n= k(2^{s} - 1) + r\cdot 2^{s}
\end{equation}
 whenever $k,\, s,\, r$ are non-negative integers satisfying $\mu(r) < k \leq h.$ (For more comprehensive information regarding this matter, kindly refer to Remark \ref{nxpt} in Sect.\ref{s2}.) The dual problem to the hit problem for the algebra $P^{\otimes h}$ is to ascertain a subring consisting of elements of the Pontrjagin ring $H_*(BV_h) = [P^{\otimes h}]^{*}$ that are mapped to zero by all Steenrod squares of positive degrees. This subring is commonly denoted by  ${\rm Ann}_{\overline{\mathcal A}}[P^{\otimes h}]^{*}.$ Let $GL_h = GL(V_h)$ be the general linear group. This $GL_h$ acts on $V_h$ and then on $QP_n^{\otimes h}.$ For each positive integer $n,$ denote by $[QP_n^{\otimes h}]^{GL_h}$ the subspace of elements that are invariant under the action of $GL_h.$ It is known that there exists an isomorphism between $(\mathbb F_2\otimes_{GL_h}{\rm Ann}_{\overline{\mathcal A}}[P^{\otimes h}]^{*})_n$ and $[QP_n^{\otimes h}]^{GL_h}$, which establishes a close relationship between the hit problem and the $h$-th algebraic transfer \cite{Singer},
$$Tr_h^{\mathcal A}: (\mathbb F_2\otimes_{GL_h}{\rm Ann}_{\overline{\mathcal A}}[P^{\otimes h}]^{*})_n\longrightarrow {\rm Ext}_{\mathcal A}^{h, h+n}(\mathbb F_2, \mathbb F_2).$$ The homomorphism $Tr_h^{\mathcal A}$ was constructed by William Singer while studying the Ext groups, employing the modular invariant theory. One notable aspect is that the Singer transfer can be regarded as an algebraic formulation of the stable transfer $B(V_h)_+^{S}\longrightarrow \mathbb S^{0}.$ It is a well-established fact, as demonstrated by Liulevicius \cite{Liulevicius}, that there exist squaring operations $Sq^{i}$ for $i\geq 0$ that act on the $\mathbb F_2$-cohomology of the Steenrod algebra $\mathcal A$. These operations share many of the same properties as the Steenrod operations $Sq^{i}$ that act on the $\mathbb F_2$-cohomology of spaces. Nonetheless, $Sq^{0}$ is not the identity. On the other side, there exists an analogous squaring operation $Sq^{0}$, called the Kameko operation, which acts on the domain of the algebraic transfer and commutes with the classical $Sq^{0}$ on ${\rm Ext}_{\mathcal A}^{*,*}(\mathbb F_2, \mathbb F_2)$ thourgh Singer's transfer (see Sect.\ref{s2} for its precise meaning). Hence, the highly non-trivial character of the algebraic transfer establishes it as a tool of potential in the study of the inscrutable Ext groups. Moreover, the hit problem and the Singer transfer have been shown in the papers \cite{Minami0, Minami} to be significant tools for investigating the Kervaire invariant one problem. It is noteworthy that Singer made the following prediction.

\begin{conj}[see \cite{Singer}]\label{gtSinger}
The transfer $Tr_h^{\mathcal A}$ is a one-to-one map for any $h.$
\end{conj}
 
Despite not necessarily resulting in a one-to-one correspondence, Singer's transfer is a valuable tool for analyzing the structure of the Ext groups. It is established, based on the works of Singer \cite{Singer} and Boardman \cite{Boardman}, that the Singer conjecture is true for homological degrees up to $3$. In these degrees, the transfer is known to be an isomorphism. We are thrilled to announce that our recent works \cite{Phuc12, Phuc10-3, Phuc10-2} has finally brought closure to the complex and long-standing problem of verifying Singer's conjecture in the case of rank four. Our study, detailed in \cite{Phuc12, Phuc10-3, Phuc10-2}, specifically establishes the truth of Conjecture \ref{gtSinger} in the case where $h=4.$ For some information on the interesting case of rank five, we recommend consulting works such as \cite{Phuc4, Phuc6, Phuc10, Phuc10-0, Sum4}. It is essential to underscore that the isomorphism between the domain of the  homomorphism $Tr_h^{\mathcal A}$ and $(QP_n^{\otimes h})^{GL_h}$ (the subspace of $GL_h$-invariants of $QP_n^{\otimes h}$) implies that it is sufficient to explore Singer's transfer in internal degrees $n$ of the form \eqref{pt}. 
\medskip

Despite extensive research, no all-encompassing methodology exists for the investigation of the hit problem and Singer's algebraic transfer in every positive degree. Therefore, each computation holds considerable importance and serves as an independent contribution to these subjects. By this reason, our primary objective in this work is to extend the findings of Mothebe et al. \cite{MKR} regarding the hit problem of six variables, while simultaneously verifying Singer's conjecture for all ranks $h\geq 1$ in certain internal degrees. Our methodology is based on utilizing the techniques developed for the hit problem, which have proven to be quite effective in determining the Singer transfer. 
More precisely, using the calculations in \cite{MKR}, we embark on an investigation of Singer's conjecture for bidegrees $(h, h+n)$, where $h\geq 1$ and $1\leq n\leq 10 = 6(2^{0}-1) + 10\cdot 2^{0}$. Subsequently, we proceed to solve the hit problem for $P^{\otimes 6} = \mathbb F_2[t_1, \ldots, t_6]$ in degrees of the form \eqref{pt}, with $k = 6$ and $r = 1$ (i.e., degree $n:=n_s= 6(2^{s}-1) + 10\cdot 2^{s}$, $s\geq 0$). Furthermore, for $h\geq 6$ and degree $2^{s+4}-h$, we establish that for each $s\geq h-5$. the dimension of the cohit module $QP^{\otimes h}_{2^{s+4}-h}$ is equal to the order of the factor group of $GL_{h-1}$ by the Borel subgroup $B_{h-1}.$  Additionally, utilizing the algebra $\pmb{A}_q$ of Steenrod reduced powers over the Galois field $\mathbb F_{q}$ of $q = p^{m}$ elements, and based on Hai's recent work \cite{Hai}, we assert that for any $h\geq 2,$ the dimension of the cohit module $\mathbb F_q[t_1, \ldots t_h]/\overline{\pmb{A}}_q\mathbb F_q[t_1, \ldots t_h]$ in degree $q^{h-1}-h$ is equal to the order of the factor group of $GL_{h-1}(\mathbb F_q)$ by a subgroup of the Borel group $B_{h-1}(\mathbb F_q).$ As a result, we establish the dimension result for the space $QP^{\otimes 7}$ in degrees $n_{s+5}$, where $s > 0$, and explicitly determine the dimension of the domain of $Tr_6^{\mathcal A}$ in degrees $n_s$. Our findings reveal that $Tr_h^{\mathcal A}$ is an isomorphism in some degrees $\leq 10$ and that $Tr_6^{\mathcal A}$ does not detect the non-zero elements $h_2^{2}g_1 = Sq^{0}(h_1Ph_1) = h_4Ph_2$ and $D_2$ in the sixth cohomology groups ${\rm Ext}_{\mathcal A}^{6, 6+n_s}(\mathbb F_2, \mathbb F_2)$. This finding carries significant implications for Singer's conjecture on algebraic transfers. Specifically, we affirm the validity of Conjecture \ref{gtSinger} for the bidegrees $(h, h+n)$ where $h\geq 1$, $1\leq n\leq n_0$,  as well as for any bidegree $(6, 6+n_s).$  

\medskip

{\bf Organization of the rest of this work.}\ 
In Sect.\ref{s2}, we provide a brief overview of the necessary background material. The main findings are then presented in Sect.\ref{s3}, with the proofs being thoroughly explained in Sect.\ref{s4}. As an insightful consolidation, Sect.\ref{s5} will encapsulate the core essence of the paper by distilling its key discoveries and notable contributions. Finally, in the appendix (referred to as Sect.\ref{s6}), we provide an extensive list of admissible monomials of degree $n_1 =6(2^{1}-1) + 10\cdot 2^{1}$ in the $\mathcal A$-module $P^{\otimes 6}$  and some $\Sigma_6$-invariants of $QP^{\otimes 6}_{n_1}$ corresponding to certain weight vectors.

\begin{acknow}
I would like to extend my heartfelt gratitude and appreciation to the editor and the anonymous referees for their meticulous review of my work and insightful recommendations, which greatly contributed to the enhancement of this manuscript. 
Furthermore, their constructive feedback has served as a source of inspiration for me to elevate the outcomes of my research beyond the initial expectations. I am greatly appreciative of Bill Singer for his keen interest, constructive criticism, and enlightening e-mail correspondence. My sincere thanks are also due to Bob Bruner, John Rognes and Weinan Lin for helpful discussions on the Ext groups. \textbf{A condensed version of this article was published in the Journal of Algebra, 613 (2023), 1--31.}

\end{acknow}

\section{Several fundamental facts}\label{s2}


In order to provide the reader with necessary background information for later use, we present some related knowledge before presenting the main results. Readers may also refer to \cite{Kameko, Phuc4, Sum2} for further information. 
We will now provide an overview of some fundamental concepts related to the hit problem.

\medskip

{\bf Weight and exponent vectors of a monomial.}\ For a natural number $k,$ writing $\alpha_j(k)$ and $\alpha(k)$ for the $j$-th coefficients and the number of $1$'s in dyadic expansion of $k$, respectively. Thence, $\alpha(k) = \sum_{j\geq 0}\alpha_j(k),$ and $k$ can be written in the form $\sum_{j\geq 0}\alpha_j(k)2^j.$  For a monomial $t = t_1^{a_1}t_2^{a_2}\ldots t_h^{a_h}$ in $P^{\otimes h},$ we consider a sequence associated with $t$ by $\omega(t) :=(\omega_1(t), \omega_2(t), \ldots, \omega_i(t), \ldots)$ where $\omega_i(t)=\sum_{1\leq j\leq h}\alpha_{i-1}(a_j)\leq h,$ for all $i\geq 1.$ This sequence is called the {\it weight vector} of $t.$ One defines $\deg(\omega(t)) = \sum_{j\geq 1}2^{j-1}\omega_j(t).$ We use the notation $a(t) = (a_1, \ldots, a_h)$ to denote the exponent vector of $t.$ Both the sets of weight vectors and exponent vectors are assigned the left lexicographical order.

\medskip

{\bf The linear order on \mbox{\boldmath $P^{\otimes h}$.}}\ Consider the monomials $t = t_1^{a_1}t_2^{a_2}\ldots t_h^{a_h}$ and $t' = t_1^{a'_1}t_2^{a'_2}\ldots t_h^{a'_h}$ in the $\mathcal A$-module $P^{\otimes h}$. We define the relation "$<$" between these monomials as follows: $t < t'$ if and only if either $\omega(t) < \omega(t')$, or $\omega(t) = \omega(t')$ and $a(t) < a'(t').$

\medskip

{\bf  The equivalence relations on \mbox{\boldmath $P^{\otimes h}$.}}\ Let $\omega$ be a weight vector of degree $n$. We define two subspaces of $P_n^{\otimes h}$ associated with $\omega$ as follows: $P_n^{\otimes h}(\omega) = {\rm span}\{ t\in P_n^{\otimes h}|\, \deg(t) = \deg(\omega) = n,\  \omega(t)\leq \omega\},$ and $P_n^{\otimes h}(< \omega) = {\rm span}\{t\in P_n^{\otimes h}|\, \deg(t) = \deg(\omega) = n,\  \omega(t) < \omega\}.$ Let $u$ and $v$ be two homogeneous polynomials in $P_n^{\otimes h}.$ We define the equivalence relations  "$\equiv$" and  "$\equiv_{\omega}$" on $P_n^{\otimes h}$ by setting: $u\equiv v$ if and only if $(u+v)\in \overline{\mathcal {A}}P_n^{\otimes h}$, while $u \equiv_{\omega} v$ if and only if $u,\, v\in  P_n^{\otimes h}(\omega) $ and $$(u +v)\in \overline{\mathcal {A}}P_n^{\otimes h} \cap P_n^{\otimes h}(\omega) + P_n^{\otimes h}(< \omega).$$ (In particular, if $u\equiv 0,$ then $u$ is a hit monomial. If $u \equiv_{\omega} 0,$ then $u$ is called {\it $\omega$-hit}.) 

We will denote the factor space of $P_n^{\otimes h}(\omega)$ by the equivalence relation $\equiv_{\omega}$ as $QP_n^{\otimes h}(\omega)$. According to \cite{Sum4, Walker-Wood}, this $QP_n^{\otimes h}(\omega)$ admits a natural $\mathbb F_2[GL_h]$-module structure, and the reader is recommended to \cite{Sum4} for a detailed proof. It is noteworthy that if we define $\widetilde{(QP^{\otimes h}_{n})^{\omega}}:= \langle \{[t]\in QP^{\otimes h}_{n}:\ \omega(t) = \omega,\, \mbox{ and $t$ is admissible} \}\rangle,$ then this $\widetilde{(QP^{\otimes h}_{n})^{\omega}}$ is an $\mathbb F_2$-subspace of $QP^{\otimes h}_n.$ Furthermore, the mapping $QP^{\otimes h}_{n}(\omega)\longrightarrow \widetilde{(QP^{\otimes h}_{n})^{\omega}}$ determined by $[t]_{\omega}\longmapsto [t]$ is an isomorphism. This implies that $ \dim QP^{\otimes h}_{n} =  \sum_{\deg(\omega) = n}\dim\widetilde{(QP^{\otimes h}_{n})^{\omega}} = \sum_{\deg(\omega) = n}\dim QP^{\otimes h}_{n}(\omega).$

\medskip

\begin{cy}\label{cyP}

(i) The conjecture proposed by Kameko in the thesis \cite{Kameko} asserts that $\dim QP^{\otimes h}_{n}\leq \prod_{1\leq j\leq h}(2^{j}-1)$ for all values of $h$ and $n$. While this inequality has been proven for $h \leq 4$ and all $n$, counterexamples provided by Sum in \cite{Sum00, Sum2} demonstrate that it is wrong when $h > 4$. It is worth noting, however, that the local version of Kameko's conjecture, which concerns the inequality $\dim QP^{\otimes h}_{n}(\omega)\leq \prod_{1\leq j\leq h}(2^{j}-1)$, remains an open question.

(ii) As it is known, the algebra of divided powers $[P^{\otimes h}]^{*} = H_{*}(BV_h)= \Gamma(a_1, a_2, \ldots, a_h)$ is generated by $a_1, \ldots, a_h,$ each of degree $1.$ Here $a_i = a_i^{(1)}$ is dual to $t_i\in P^{\otimes h}_1,$ with duality taken with respect to the basis of $P^{\otimes h}$ consisting of all monomials in $t_1, \ldots, t_h.$ Kameko defined in \cite{Kameko} a homomorphism of $\mathbb F_2[GL_h]$-modules  $\widetilde{Sq}^{0}: [P^{\otimes h}]^{*}=H_{*}(BV_h)\longrightarrow [P^{\otimes h}]^{*} = H_{*}(BV_h),$ which is determined by $\widetilde{Sq}^{0}(a_1^{(i_1)}\ldots a_h^{(i_h)}) = a_1^{(2i_1+1)}\ldots a_h^{(2i_h+1)}.$ The dual of this $\widetilde{Sq}^{0}$ induced the homomorphism $(\widetilde {Sq^0_*})_{2n+h}: QP^{\otimes h}_{2n+h}  \longrightarrow QP^{\otimes h}_{n}$ (see Sect.\ref{s1}). Further, as $Sq^{2k+1}_*\widetilde{Sq}^{0} = 0$ and $Sq_{*}^{2k}\widetilde{Sq}^{0} = \widetilde{Sq}^{0}Sq_{*}^{k},$ $\widetilde{Sq}^{0}$ maps ${\rm Ann}_{\overline{\mathcal A}}[P^{\otimes h}]^{*}$ to itself. Here we write $Sq^u_*: H_{*}(BV_h)\longrightarrow H_{*-u}(BV_h)$ for the operation on homology which by duality of vector spaces is induced by the square $Sq^{u}: H_{*}(BV_h)\longrightarrow H_{*+u}(BV_h).$ The Kameko  $Sq^{0}$ is defined by $Sq^{0}: \mathbb F_2\otimes_{GL_h}{\rm Ann}_{\overline{\mathcal A}}[P^{\otimes h}]^{*}\longrightarrow \mathbb F_2\otimes_{GL_h}{\rm Ann}_{\overline{\mathcal A}}[P^{\otimes h}]^{*},$ which commutes with the classical $Sq^{0}$ on the $\mathbb F_2$-cohomology of $\mathcal A$ through the Singer algebraic transfer. Thus, for any integer $n\geq 1,$ the following diagram commutes:
$$ \begin{diagram}
\node{(\mathbb F_2\otimes_{GL_h}{\rm Ann}_{\overline{\mathcal A}}[P^{\otimes h}]^{*})_n} \arrow{e,t}{Tr_h^{\mathcal A}}\arrow{s,r}{Sq^0} 
\node{{\rm Ext}_{\mathcal A}^{h, h+n}(\mathbb F_2, \mathbb F_2)} \arrow{s,r}{Sq^0}\\ 
\node{(\mathbb F_2\otimes_{GL_h}{\rm Ann}_{\overline{\mathcal A}}[P^{\otimes h}]^{*})_{2n+h}} \arrow{e,t}{Tr_h^{\mathcal A}} \node{{\rm Ext}_{\mathcal A}^{h, 2h+2n}(\mathbb F_2, \mathbb F_2)}
\end{diagram}$$

Thus, Kameko's $Sq^0$ is known to be compatible via the Singer transfer with $Sq^0$ on ${\rm Ext}_{\mathcal A}^{*,*}(\mathbb F_2, \mathbb F_2)$. Moreover, the $GL_h$-coinvariants $(\mathbb F_2\otimes_{GL_h}{\rm Ann}_{\overline{\mathcal A}}[P^{\otimes h}]^{*})_n$ form a bigraded algebra and the Singer algebraic transfers $Tr_*^{\mathcal A}$ yield a morphism of bigraded algebras with values in ${\rm Ext}_{\mathcal A}^{*,*}(\mathbb F_2, \mathbb F_2).$ These compatibilities are suggestive of a far closer relationship between these structures. In addition, the operations $Sq^0$ and the algebra structure on ${\rm Ext}_{\mathcal A}^{*,*}(\mathbb F_2, \mathbb F_2)$ are key ingredients in understanding the image of the algebraic transfer. Unfortunately, detecting the image of the Singer transfer by mapping \textit{out} of ${\rm Ext}_{\mathcal A}^{*,*}(\mathbb F_2, \mathbb F_2)$ is not easy. For example, Lannes and Zarati  \cite{LZ} constructed an algebraic approximation to the Hurewicz map: for an unstable $\mathcal A$-module $M$ this is of the form ${\rm Ext}_{\mathcal A}^{h, h+n}(\Sigma^{-h}M, \mathbb F_2)\longrightarrow [(\mathbb F_2\otimes_{\mathcal A}\mathscr R_hM)_n]^*,$  where $\mathscr R_h$ is the $h$-th Singer functor (as defined by Lannes and Zarati). However, it is conjectured by H\uhorn ng \cite{Hung0} that this vanishes for $h>2$ in positive stem, an algebraic version of the long-standing and difficult \textit{generalized spherical class conjecture} in Algebraic topology, due to Curtis \cite{Curtis}. In \cite{Hung2}, H\uhorn ng and Powell proved the weaker result that this holds on the image of the transfer homomorphism.  This illustrates the difficulty of studying the Singer transfer.

An analogous diagram has also been established for the case of odd primes $p$ \cite{Minami}:
$$ \begin{diagram}
\node{(\mathbb F_p\otimes_{GL_h(\mathbb F_p)}{\rm Ann}_{\overline{\mathcal A_p}}H_*(V; \mathbb F_p))_n} \arrow{e,t}{Tr_h^{\mathcal A_p}}\arrow{s,r}{Sq^{0}} 
\node{{\rm Ext}_{\mathcal A_p}^{h, h+n}(\mathbb F_p, \mathbb F_p)} \arrow{s,r}{Sq^{0}}\\ 
\node{(\mathbb F_p\otimes_{GL_h(\mathbb F_p)}{\rm Ann}_{\overline{\mathcal A_p}}H_*(V; \mathbb F_p))_{p(n+h)-h}} \arrow{e,t}{Tr_h^{\mathcal A_p}} \node{{\rm Ext}_{\mathcal A_p}^{h, p(h+n)}(\mathbb F_p, \mathbb F_p)} 
\end{diagram}$$
Here, the left vertical arrow represents the Kameko $Sq^0,$ and the right vertical one represents the classical squaring operation. Our recent work \cite{Phuc14} proposes a conjecture that \textit{the transfer $Tr_h^{\mathcal A_p}$ is an injective map for all $1\leq h\leq 4$ and odd primes $p$.} We have also established the validity of this conjecture in certain generic degrees.
\end{cy}

 {\bf (Strictly) inadmissible monomial.}\ We say that a monomial $t\in P_n^{\otimes h}$ is {\it inadmissible} if there exist monomials $z_1, z_2,\ldots, z_k\in P_n^{\otimes h}$ such that $z_j < t$ for $1\leq j\leq k$ and  $t = \sum_{1\leq j\leq k}z_j  + \sum_{m > 0}Sq^{m}(z_m),$ for some $m\in \mathbb N$ and $z_m\in P_{n-m}^{\otimes h},\, m<n.$ Then, $t$ is said to be {\it admissible} if it is not inadmissible. A monomial $t\in P_n^{\otimes h}$ is said to be {\it strictly inadmissible} if and only if there exist monomials $z_1, z_2,\ldots, z_k$ in $P_n^{\otimes h}$ such that $z_j < t$ for $1\leq j \leq k$ and $t = \sum_{1\leq j\leq k}z_j + \sum_{0\leq  m \leq s-1}Sq^{2^{m}}(z_m),$ where $s = {\rm max}\{i\in\mathbb N: \omega_i(t) > 0\}$ and suitable polynomials $z_m\in P^{\otimes h}_{n-2^{m}}.$

\medskip
Note that every strictly inadmissible monomial is inadmissible but the converse is not generally true. For example, consider the monomial $t = t_1t_2^{2}t_3^{2}t_4^{2}t_5^{6}t_6\in P^{\otimes 6}_{14},$ we see that this monomial is not strictly inadmissible, despite its inadmissibility. This can be demonstrated through the application of the Cartan formula, which yields $$t = Sq^{1}(t_1^{4}t_2t_3t_4t_5^{5}t_6) + Sq^{3}(t_1^{2}t_2t_3t_4t_5^{5}t_6) + Sq^{6}(t_1t_2t_3t_4t_5^{3}t_6) + \sum_{X<t}X.$$ 
\begin{therm}[{\bf Criteria on inadmissible monomials}]\label{dlKS}
The following claims are each true:
\begin{itemize}
\item[(i)]  Let $t$ and $z$ be monomials in $P^{\otimes h}.$ For an integer $r >0,$ assume that there exists an index $i>r$ such that $\omega_i(t) = 0.$ If $z$ is inadmissible, then $tz^{2^r}$ is, too {\rm (see Kameko \cite{Kameko})};

\medskip

\item[(ii)] Let $z, w$ be monomials in $P^{\otimes h}$ and let $r$ be a positive integer. Suppose that there is an index $j > r$ such that $\omega_j(z) = 0$ and $\omega_r(z)\neq 0.$ If $z$ is strictly inadmissible, then so is, $zw^{2^{r}}$ {\rm (see Sum \cite{Sum2})}.
\end{itemize}
\end{therm}

\medskip

We shall heavily rely on the arithmetic function $\mu: \mathbb{N}\longrightarrow \mathbb{N},$ as well as Kameko's map $(\widetilde{Sq^0_*})_{2n+h}: QP^{\otimes h}_{2n+h}\longrightarrow QP^{\otimes h}_{n}$, both of which are elucidated in Sect. \ref{s1}. The technical theorem below related to the $\mu$-function holds crucial significance.

\begin{therm}\label{dlWS}
The following statements are each true:

\begin{itemize}

\item[(i)] {\rm (cf. Sum \cite{Sum4})}.   $\mu(n) = r\leq h$ if and only if there exists uniquely a sequence of positive integers $d_1 > d_2 > \cdots > d_{r-1} \geq d_r$ such that $n = \sum_{1\leq i\leq r}(2^{d_i} - 1).$ 

\medskip

\item[(ii)] {\rm (cf. Wood \cite{Wood})}.  For each positive integer $n,$ the space $QP^{\otimes h}_n$ is trivial if and only if $\mu(n) > h.$ 

\medskip

\item[(iii)]  {\rm (cf. Kameko \cite{Kameko})}.  The homomorphism $(\widetilde {Sq^0_*})_{2n+h}$ is an isomorphism of $\mathbb F_2$-vector spaces if and only if $\mu(2n+h) = h.$ 

\end{itemize}
\end{therm}

\begin{rems}\label{nxpt}
In Sect.\ref{s1}, it was noted that the hit problem needs to be solved only in degrees of the form \eqref{pt}. Furthermore, in \cite[Introduction]{Sum2}, Sum made the remark that for every positive integer $n$, the condition $3\leq \mu(n)\leq h$ holds if and only if there exist uniquely positive integers $s$ and $r$ satisfying $1\leq \mu(n)-2\leq \mu(r) = \alpha(r+\mu(r))\leq \mu(n)-1$ and $n = \mu(n)(2^{s}-1) + r\cdot2^{s}.$ This can be demonstrated straightforwardly by utilizing Theorem \ref{dlWS}(i). Suppose that $\mu(n) = k\geq 3.$ Then, the "only if" part has been shown in \cite{Phuc4}. The "if" part is established as follows: if $n = k(2^s-1) + r\cdot 2^s$ and $1\leq k-2\leq \mu(r) = \alpha(r+\mu(r)) \leq k-1.$ Then, either $\mu(r) = k-2$ or $\mu(r) = k-1.$  We set  $\mu(r) = \ell$ and see that by Theorem \ref{dlWS}(i),  there exist uniquely a sequence of integers $c_1 > c_2>\ldots>c_{\ell-1} \geq c_{\ell} > 0$ such that $r = 2^{c_1} + 2^{c_2} + \cdots + 2^{c_{\ell-1}} + 2^{c_{\ell}} - \ell.$ Obviously, $\alpha(r + \ell) = \ell,$ and so, $n =k(2^s-1) + r\cdot2^{s} = 2^{c_1 + s} + 2^{c_2 + s} + \cdots + 2^{c_{\ell} + s} + 2^{s}(s-\ell) - s.$ Now, if $\ell = k-2,$ then $n = 2^{c_1 + s} + 2^{c_2 + s} + \cdots + 2^{c_{\ell} + s} + 2^{s}(s-\ell) - s = 2^{c_1 + s} + 2^{c_2 + s} + \cdots + 2^{c_{k-2} + s} + 2^{s} + 2^{s} - k.$ Let $u_i = c_i + s$ with $1\leq i\leq k-2$ and let $u_{k-1} = u_k = s.$ Since $u_1 > u_2 > \cdots > u_{k-2} > u_{k-1} = u_k,$ by Theorem \ref{dlWS}(i), $\mu(n) = k.$ Finally, if $\ell = k-1$ then $ n = 2^{c_1 + s} + 2^{c_2 + s} + \cdots + 2^{c_{k-1} + s} + 2^{s} - k.$ We put $v_i = c_i + s$ where $1\leq i\leq k-1$ and $v_k = s.$ Since $v_1 > v_2 > \cdots > v_{k-1} >  v_k,$ according to Theorem \ref{dlWS}(i), one derives $\mu(n) = k.$
\end{rems}

{\bf Spike monomial.}\ A monomial $t_1^{a_1}t_2^{a_2}\ldots t_h^{a_h}$ in $P^{\otimes h}$ is called a {\it spike} if every exponent $a_j$ is of the form $2^{\beta_j} - 1.$ In particular, if the exponents $\beta_j$ can be arranged to satisfy $\beta_1 > \beta_2 > \ldots > \beta_{r-1}\geq \beta_r \geq 1,$ where only the last two smallest exponents can be equal, and $\beta_j = 0$ for $ r < j  \leq h,$ then the monomial $t_1^{a_1}t_2^{a_2}\ldots t_h^{a_h}$ is called a {\it minimal spike}. 

\medskip

\begin{therm}[see Ph\'uc and Sum \cite{P.S1}]\label{dlPS}
All the spikes in $P^{\otimes h}$ are admissible and their weight vectors are weakly decreasing. Furthermore, if a weight vector $\omega = (\omega_1, \omega_2, \ldots)$ is weakly decreasing and $\omega_1\leq h,$ then there is a spike $z\in P^{\otimes h}$ such that $\omega(z) = \omega.$
\end{therm}

The subsequent information demonstrates the correlation between minimal spike and hit monomials.

\begin{therm}[{\bf  Singer's criterion on hit monomials} {\rm \cite{Singer1}}]\label{dlSin}
Suppose that $t\in P^{\otimes h}$ and $\mu(\deg(t))\leq h.$ Consequently, if $z$ is a minimal spike in $P^{\otimes h}$ such that $\omega(t) < \omega(z),$ then $t\equiv 0$ (or equivalently, $t$ is hit).
\end{therm}

It is of importance to observe that the converse of Theorem \ref{dlSin} is generally not valid. As a case in point, let us consider $z = t_1^{31}t_2^{3}t_3^{3}t_4^{0}t_5^{0}\in   P^{\otimes 5}_{37}$ and $t = t_1(t_2t_3t_4t_5)^{9}\in  P^{\otimes 5}_{37}.$ One has $\mu(37) = 3 < 5,$ and $t  = fg^{2^3},$ where $f = t_1t_2t_3t_4t_5$ and $g =t_2t_3t_4t_5.$ Then $\deg(f) = 5 < (2^{3}-1)\mu(\deg(g)),$ and so, due to Silverman \cite[Theorem 1.2]{Silverman2}, we must have $t\equiv 0.$ It can be observed that despite $z$ being the minimal spike of degree $37$ in the $\mathcal A$-module $P^{\otimes 5},$ the weight $\omega(t) = (5,0,0,4,0)$ exceeds the weight of $z,$ which is $\omega(z) = (3,3,1,1,1).$ The reader may also refer to \cite{Phuc16} for further information regarding the cohit module $QP^{\otimes 5}_{37}.$

\begin{notas}
We will adopt the following notations for convenience and consistency:

\begin{itemize}
\item [$\bullet$] 
Let us denote by $(P^{\otimes h})^{0}:= {\rm span}\bigg\{\prod_{1\leq j\leq h}t_j^{\alpha_j} \in P^{\otimes h}\bigg|\, \prod_{1\leq j\leq h}\alpha_j = 0\bigg\}$ and $(P^{\otimes h})^{> 0}:= {\rm span}\bigg \{\prod_{1\leq j\leq h}t_j^{\alpha_j} \in P^{\otimes h}\bigg|\, \prod_{1\leq j\leq h}\alpha_j > 0\bigg\}.$ It can be readily observed that these spaces are $\mathcal A$-submodules of $P^{\otimes h}.$ Moreover, for each positive degree $n,$ we have $QP_n^{\otimes h} \cong (QP_n^{\otimes h})^0\,\bigoplus\, (QP_n^{\otimes h})^{>0},$ where $(QP_n^{\otimes h})^0:= (Q(P^{\otimes h})^{0})_n = (\mathbb F_2\otimes_{\mathcal A} (P^{\otimes h})^{0})_n$ and $(QP_n^{\otimes h})^{>0}:= (Q(P^{\otimes h})^{>0})_n = (\mathbb F_2\otimes_{\mathcal A} (P^{\otimes h})^{>0})_n$  are the $\mathbb F_2$-subspaces of $QP_n^{\otimes h}.$ 

\medskip

\item [$\bullet$] Given a monomial $t\in P^{\otimes h}_n,$ we write $[t]$ as the equivalence class of $t$ in $QP^{\otimes h}_n.$ If $\omega$ is a weight vector of degree $n$ and $t\in P_n^{\otimes h}(\omega),$ we denote by $[t]_\omega$ the equivalence class of $t$ in $QP^{\otimes h}_n(\omega).$ Noting that if $\omega$ is a weight vector of a minimal spike, then $[t]_{\omega} = [t].$ For a subset $C\subset P^{\otimes h}_n,$ we will often write $|C|$ to denote the cardinality of $C$ and use notation $[C] = \{[t]\, :\, t\in C\}.$ If $C\subset  P_n^{\otimes h}(\omega),$ then we denote $[C]_{\omega} = \{[t]_{\omega}\, :\, t\in C\}.$

\medskip

\item [$\bullet$] Write $\mathscr {C}^{\otimes h}_{n},\, (\mathscr {C}^{\otimes h}_{n})^{0}$ and $(\mathscr {C}^{\otimes h}_{n})^{>0}$ as the sets of all the admissible monomials of degree $n$ in the $\mathcal A$-modules $P^{\otimes h},$\ $(P^{\otimes h})^{0}$ and $(P^{\otimes h})^{>0},$ respectively. If $\omega$ is a weight vector of degree $n.$ then we put $\mathscr {C}^{\otimes h}_{n}(\omega) := \mathscr {C}^{\otimes h}_{n}\cap P_n(\omega),$\ $(\mathscr {C}^{\otimes h}_{n})^{0}(\omega) := (\mathscr {C}^{\otimes h}_{n})^{0}\cap P_n^{\otimes h}(\omega),$ and $(\mathscr {C}^{\otimes h}_{n})^{>0}(\omega) := (\mathscr {C}^{\otimes h}_{n})^{>0}\cap P_n^{\otimes h}(\omega).$
\end{itemize}
\end{notas}

\section{Statement of main results}\label{s3}


We are now able to present the principal findings of this paper. 
The demonstration of these results will be exhaustively expounded in subsequent section. As previously alluded to, our examination commences with a critical analysis of the hit problem for the polynomial algebra $P^{\otimes 6}$ in degree $n_s:= 6(2^{s}-1) + 10\cdot 2^{s}$, where $s$ is an arbitrary non-negative integer.

\medskip
{\bf Case $\pmb{s = 0.}$}\ Mothebe et al. demonstrated in \cite{MKR} the following outcome.

\begin{therm}[see \cite{MKR}]\label{dlMKR}
For each integer $h\geq 2,$ $\dim QP^{\otimes h}_{n_0} = \sum_{2\leq j\leq n_0}C_j\binom{h}{j},$ where $\binom{h}{j} = 0$ if $h < j$ and $C_2 = 2,$ $C_3 = 8,$ $C_4 = 26,$ $C_5 = 50,$ $C_6 = 65,$ $C_7 = 55,$ $C_8 = 28,$ $C_9 = 8,$ $C_{n_0} = 1.$ This means that there exist exactly $945$ admissible monomials of degree $n_0$ in the $\mathcal A$-module $P^{\otimes 6}.$ 
\end{therm}

The following corollary is readily apparent.

\begin{corls}\label{hq10-1}
\begin{itemize}
\item[(i)]
One has an isomorphism of $\mathbb F_2$-vector spaces: $ QP^{\otimes 6}_{n_0}\cong \bigoplus_{1\leq j\leq 5} QP^{\otimes 6}_{n_0}(\overline{\omega}^{(j)}),$ where $\overline{\omega}^{(1)}:= (2,2,1),$\ $\overline{\omega}^{(2)}:=(2,4),$\ $\overline{\omega}^{(3)}:=(4,1,1),$\ $\overline{\omega}^{(4)}:=(4,3),$ and $\overline{\omega}^{(5)}:=(6,2).$ 

\item[(ii)] $(QP^{\otimes 6}_{n_0})^{0}\cong \bigoplus_{1\leq j\leq 4}(QP^{\otimes 6}_{n_0})^{0}(\overline{\omega}^{(j)})$ and $(QP^{\otimes 6}_{n_0})^{>0}\cong \bigoplus_{2\leq j\leq 5}(QP^{\otimes 6}_{n_0})^{>0}(\overline{\omega}^{(j)}),$ and 

\centerline{\begin{tabular}{c||cccccc}
$j$  &$1$ & $2$ & $3$ & $4$& $5$\cr
\hline
\hline
\ $\dim (QP^{\otimes 6}_{n_0})^{0}(\overline{\omega}^{(j)})$ & $400$ & $30$ & $270$ &$180$ &$0$  \cr
\hline
\hline
\ $\dim (QP^{\otimes 6}_{n_0})^{>0}(\overline{\omega}^{(j)})$ & $0$ & $4$ & $10$ &$36$ & $15$  \cr
\end{tabular}}
\end{itemize}
\end{corls}

It is worth noting that the epimorphism of $\mathbb{F}_2$-vector spaces, Kameko's squaring operation $(\widetilde {Sq^0_*})_{n_0}: QP^{\otimes 6}_{n_0} \longrightarrow (QP^{\otimes 6}_{2})^{0}$, implies that $QP^{\otimes 6}_{n_0}$ is isomorphic to ${\rm Ker}((\widetilde {Sq^0_*})_{n_0})\bigoplus \psi((QP_{2}^{\otimes 6})^{0})$. Here, $\psi: (QP^{\otimes 6}_{2})^{0}\longrightarrow QP^{\otimes 6}_{n_0}$ is induced by the up Kameko map $\psi: (P^{\otimes 6}_2)^{0}\longrightarrow P^{\otimes 6}_{n_0},\, t\longmapsto t_1t_2\ldots t_6t^{2}$. Hence, by virtue of Corollary \ref{hq10-1}, one has the isomorphisms: ${\rm Ker}((\widetilde {Sq^0_*})_{n_0})\cong \bigoplus_{1\leq j\leq 4}QP^{\otimes 6}_{n_0}(\overline{\omega}^{(j)}),$ and $\psi((QP_{2}^{\otimes 6})^{0}) \cong QP^{\otimes 6}_{n_0}(\overline{\omega}^{(5)}).$

\begin{rems}\label{nxp10}
Let  us consider the set $\mathcal L_{h, k} = \{J = (j_1, \ldots, j_k):\, 1\leq j_1 < j_2 < \cdots < j_k\leq h \},\ 1\leq k < h.$ Obviously, $|\mathcal L_{h, k}| = \binom{h}{k}.$  For each $J\in \mathcal L_{h, k},$ we define the the homomorphism $\varphi_J: P^{\otimes k}\longrightarrow P^{\otimes h}$ of algebras by setting $\varphi_J(t_{u}) = t_{j_{u}},\ 1\leq u\leq h.$ It is straightforward to see that this homomorphism is also a homomorphism of $\mathcal A$-modules. For each $1\leq k<h,$ we have the isomorphism of $\mathbb F_2$-vector spaces $Q(\varphi_J((P^{\otimes k})^{>0}))_{n}(\omega) = (\mathbb F_2\otimes_{\mathcal A} \varphi_J((P^{\otimes k})^{>0}))_{n}\cong (QP_{n}^{\otimes k})^{>0}(\omega),$ where $\omega$ is a weight vector of degree $n.$ As a consequence of this, and based on the work of \cite{Walker-Wood}, we get
$$ 
 (QP_{n}^{\otimes h})^{0}(\omega)\cong \bigoplus_{1\leq k\leq h-1}\bigoplus_{J\in\mathcal L_{h, k}}Q(\varphi_J((P^{\otimes k})^{>0}))_{n}(\omega) \cong\bigoplus_{1\leq k\leq h-1}\bigoplus_{1\leq d\leq \binom{h}{k}}(QP_{n}^{\otimes k})^{>0}(\omega),
$$
which implies $ \dim (QP_{n}^{\otimes h})^{0}(\omega) = \sum_{1\leq k\leq h-1}\binom{h}{k}\dim (QP_{n}^{\otimes k})^{> 0}(\omega).$ By utilizing Theorem \ref{dlWS}(ii) in combination, we obtain
$$
\dim (QP_{n}^{\otimes h})^{0}(\omega) = \sum_{\mu(n)\leq k\leq h-1}\binom{h}{k}\dim (QP_{n}^{\otimes k})^{> 0}(\omega).
$$
\end{rems}

Through a straightforward calculation utilizing Theorem \ref{dlMKR}, we can claim the following.

\begin{corls}\label{hq10-1-0}
Let $\overline{\omega}^{(j)}$ be the weight vectors as in Corollary \ref{hq10-1} with $1\leq j\leq 5.$ Then, for each $h\geq 7,$ the dimension of $(QP^{\otimes h}_{n_0})^{>0}(\overline{\omega}^{(j)})$ is determined by the following table:

\centerline{\begin{tabular}{c||cccccccccccccccc}
$j$  &&$1$ && $2$ && $3$ && $4$&& $5$\cr
\hline
\hline
\ $\dim (QP^{\otimes 7}_{n_0})^{>0}(\overline{\omega}^{(j)})$ && $0$ && $0$ && $0$ &&$20$ &&$35$   \cr
\hline
\hline
\ $\dim (QP^{\otimes 8}_{n_0})^{>0}(\overline{\omega}^{(j)})$ && $0$ && $0$ && $0$ &&$0$ && $20$   \cr
\hline
\hline
\ $\dim (QP^{\otimes h}_{n_0})^{>0}(\overline{\omega}^{(j)}),\, h \geq 9$ && $0$ && $0$ && $0$ &&$0$ && $0$  \cr
\end{tabular}}
\end{corls}

Through a basic computation, in conjunction with Remark \ref{nxp10}, Corollaries \ref{hq10-1}, \ref{hq10-1-0}, as well as the preceding outcomes established by \cite{Peterson}, \cite{Kameko}, and \cite{Sum2}, we are able to deduce the subsequent corollary.

\begin{corls}\label{hq10-1-1}
Let $\overline{\omega}^{(j)}$ be the weight vectors as in Corollary \ref{hq10-1} with $1\leq j\leq 5.$ Then, for each $h\geq 7,$ the dimension of $(QP^{\otimes h}_{n_0})^{0}(\overline{\omega}^{(j)})$ is given as follows:
$$ \dim (QP^{\otimes h}_{n_0})^{0}(\overline{\omega}^{(j)}) =  \left\{\begin{array}{ll}
 2\bigg[\binom{h}{2} + 4\binom{h}{3} + 6\binom{h}{4}\bigg] + 5\binom{h}{5}&\mbox{if $j = 1$},\\[1mm]
 5\binom{h}{5}+ 4\binom{h}{6} &\mbox{if $j = 2$},\\[1mm]
 10\bigg[\binom{h}{4}+2\binom{h}{5} + \binom{h}{6}\bigg]&\mbox{if $j = 3$},\\[1mm]
 812&\mbox{if $j = 4,\, h =7$},\\[1mm]
 4\bigg[\binom{h}{4}+ 5\binom{h}{5}+9\binom{h}{6} + 5\binom{h}{7}\bigg] &\mbox{if $j = 4,\, h \geq 8$},\\[1mm]
 105 &\mbox{if $j = 5,\, h =7$},\\[1mm]
 700&\mbox{if $j = 5,\, h = 8$},\\[1mm]
 5\bigg[3\binom{h}{6}+ 7\binom{h}{7}+ 4\binom{h}{8}\bigg]&\mbox{if $j = 5,\, h \geq 9$}.
\end{array}\right.$$

\end{corls}

As a direct implication of the findings presented in \cite{MKR}, we derive

\begin{corls}\label{hq10-2}
For each integer $h\geq 6,$ consider the following weight vectors of degree $n = h+4$:
$$\overline{\omega}^{(1,\, h)}:= (h-4, 4),\ \ \overline{\omega}^{(2,\, h)}:= (h-2, 1,1),\ \ \overline{\omega}^{(3,\, h)}:= (h-2, 3),\ \ \overline{\omega}^{(4,\, h)}:= (h, 2).$$
Then, for each rank $h\geq 6,$ the dimension of $(QP^{\otimes h}_{h+4})^{>0}(\overline{\omega}^{(j,\, h)})$ is determined by the following table:

\centerline{\begin{tabular}{c||cccccccccccc}
$j$  &&$1$ && $2$ && $3$ && $4$\cr
\hline
\hline
\ $\dim (QP^{\otimes h}_{h+4})^{>0}(\overline{\omega}^{(j,\, h)})$ && $\binom{h-1}{4}-1$ && $\binom{h-1}{2}$ && $h\binom{h-2}{2}$ &&$\binom{h}{2}.$  \cr
\end{tabular}}
\end{corls}


Owing to Corollary \ref{hq10-1}, one has $\overline{\omega}^{(j,\, h)} = \overline{\omega}^{(j+1)}$ for $h = 6$ and $1\leq j\leq 4.$ So we can infer that the dimension of $(QP^{\otimes 6}_{n_0})^{>0}(\overline{\omega}^{(j)}),\, 2\leq j\leq 5$ in Corollary \ref{hq10-1} can be derived from Corollary \ref{hq10-2}. Furthermore, in light of Corollaries \ref{hq10-1-0}, \ref{hq10-1-1}, \ref{hq10-2}, as well as the previous results established by Peterson \cite{Peterson}, Kameko \cite{Kameko}, Sum \cite{Sum2}, and Mothebe et al. \cite{MKR}, we are also able to confirm that a local version of Kameko's conjecture (as articulated in Note \ref{cyP}) holds true for certain weight vectors of degrees $h+4$, where $h\geq 1.$
\medskip

As is well-known, Mothebe et al. \cite{MKR} computed the dimension of $QP_n^{\otimes h}$ for $h\geq 1$ and degrees $n$ satisfying $1\leq n\leq 9$. The following theorem provides more details.

\begin{therm}[see \cite{MKR}]\label{dlMM}
Given any $h\geq 1,$ the dimension of $QP_n^{\otimes h}$ is determined as follows:
$$ \begin{array}{ll}
\medskip
\dim QP_1^{\otimes h} &=h ,\ \ \dim QP_2^{\otimes h} = \binom{h}{2},\ \ \dim QP_3^{\otimes h} = \sum_{1\leq j\leq 3}\binom{h}{j},\\
\medskip
\dim QP_4^{\otimes h} &= 2\binom{h}{2} + 2\binom{h}{3}+\binom{h}{4},\ \ \dim QP_5^{\otimes h} = 3\binom{h}{3} + 3\binom{h}{4}+\binom{h}{5},\\
\medskip
\dim QP_6^{\otimes h} &=\binom{h}{2} + 3\binom{h}{3}+6\binom{h}{4} + 4\binom{h}{5} + \binom{h}{6},\\
\medskip
\dim QP_7^{\otimes h} &= \binom{h}{1} + \binom{h}{2}+4\binom{h}{3} + 9\binom{h}{4} + 10\binom{h}{5}+ 5\binom{h}{6} + \binom{h}{7},\\
\medskip
\dim QP_8^{\otimes h} &= 3\binom{h}{2}+6\binom{h}{3} + 13\binom{h}{4} + 19\binom{h}{5}+ 15\binom{h}{6} + 6\binom{h}{7} + \binom{h}{8},\\
\medskip
\dim QP_9^{\otimes h} &= 7\binom{h}{3} + 18\binom{h}{4} + 31\binom{h}{5}+ 34\binom{h}{6} + 21\binom{h}{7} + 7\binom{h}{8}+\binom{h}{9},
\end{array}$$
where the binomial coefficients $\binom{h}{k}$ are to be interpreted modulo 2 with the usual convention $\binom{h}{k} = 0$ if $k \geq h+1.$
\end{therm} 
The theorem has also been demonstrated by Peterson \cite{Peterson} for $h\leq 2$, by Kameko's thesis \cite{Kameko} for $h = 3$ and by Sum \cite{Sum2} for $h  =4.$

Using Theorems \ref{dlMKR}, \ref{dlMM}, and Corollaries \ref{hq10-1-0}, \ref{hq10-1-1}, we aim to analyze the behavior of the Singer transfer in bidegree $(h, h+n)$ for $1\leq n\leq n_0$ and any $h\geq 1.$ As a result of our investigation, we establish the following first main result.

\begin{therm}\label{dlc0}
For any integer $n$ satisfying $1\leq n\leq n_0,$ the algebraic transfer $$Tr_h^{\mathcal A}: (\mathbb F_2\otimes_{GL_h}{\rm Ann}_{\overline{\mathcal A}}[P^{\otimes h}]^{*})_{n}\longrightarrow {\rm Ext}_{\mathcal A}^{h, h+n}(\mathbb F_2, \mathbb F_2)$$ is a trivial isomorphism for all $h\geq 1$, except for the cases of rank 5 in degree 9 and rank 6 in degree $n_0$. In these exceptional cases, $Tr_h^{\mathcal A}$ is a monomorphism. Consequently, Singer's Conjecture \ref{gtSinger} holds true in bidegrees $(h, h+n)$ for $h\geq 1$ and $1\leq n\leq n_0.$ 
\end{therm}

The theorem has been proven by Singer \cite{Singer} for $1\leq h\leq 2$, by Boardman \cite{Boardman} for $h = 3$, by Sum \cite{Sum2-1} and the author \cite{Phuc12} for $h = 4,$ by Sum \cite{Sum0, Sum2-0, Sum3} for $h =5$ and $n = 4,\, 5,\, n_0,$ by Sum and T\'in \cite{Sum-Tin0, Sum-Tin} for $h = 5$ and $n =1,\, 2,\, 3,\, 7,\, 9.$ The present writer has established the theorem for the case $h=5$ and degree $n=6, 8,$ as well as for the cases $6\leq h\leq 8$ and any degree $n$, as shown in \cite{Phuc0, Phuc6, Phuc11}. It should be brought to the attention of the readers that, $Tr_h^{\mathcal A}$ is not an epimorphism for rank 5 in degree 9 \cite{Singer}, and also for rank 6 in degree $n_0$ \cite{CHa, CHa0}. These imply that $Ph_1\not\in {\rm Im}(Tr_5^{\mathcal A})$ and $h_1Ph_1\not\in {\rm Im}(Tr_6^{\mathcal A}),$ where $\{Ph_1\}\subset {\rm Ext}_{\mathcal A}^{5, 5+9}(\mathbb F_2, \mathbb F_2)$  and $\{h_1Ph_1\}\subset {\rm Ext}_{\mathcal A}^{6, 6+n_0}(\mathbb F_2, \mathbb F_2)$ are sets that generate ${\rm Ext}_{\mathcal A}^{5, 5+9}(\mathbb F_2, \mathbb F_2)$ and ${\rm Ext}_{\mathcal A}^{6, 6+n_0}(\mathbb F_2, \mathbb F_2),$ respectively.

\medskip

{\bf Case $\pmb{s = 1}.$}\ We notice that $n_1 = 6(2^{1}-1) + 10\cdot 2^{1} = 26$ and make the following observation.

\begin{rems}\label{nxp0}

Let us consider the Kameko map $(\widetilde {Sq^0_*})_{n_1}: QP^{\otimes 6}_{n_1}  \longrightarrow QP^{\otimes 6}_{n_0}$ which is an epimorphism of the $\mathbb F_2$-vector spaces and is determined by $(\widetilde {Sq^0_*})_{n_1}([u]) = [t]$ if $u = t_1t_2\ldots t_6t^{2}$ and $(\widetilde {Sq^0_*})_{n_1}([u]) = 0$ otherwise. Then, since the homomorphism $q: {\rm Ker}((\widetilde {Sq^0_*})_{n_1})\longrightarrow  QP^{\otimes 6}_{n_1}$ is an embedding, we have a short exact sequence of $\mathbb F_2$-vector spaces.: $0\longrightarrow {\rm Ker}((\widetilde {Sq^0_*})_{n_1})\longrightarrow  QP^{\otimes 6}_{n_1}\longrightarrow QP^{\otimes 6}_{n_0}\longrightarrow 0.$  Let us consider the up Kameko map $\psi: P_{n_0}^{\otimes 6}\longrightarrow P_{n_1}^{\otimes 6},$ which is determined by $\psi(t) = t_1t_2\ldots t_6t^{2}$ for any $t\in P_{n_0}^{\otimes 6}.$ This $\psi$ induces a homomorphism $\psi: QP_{n_0}^{\otimes 6}\longrightarrow QP_{n_1}^{\otimes 6}.$ These data imply that the above exact sequence is split and so, $QP^{\otimes 6}_{n_1}  \cong {\rm Ker}((\widetilde {Sq^0_*})_{n_1})\bigoplus QP^{\otimes 6}_{n_0}.$ Furthermore, as well known, $(QP^{\otimes 6}_{n_1})^{0}$ and ${\rm Ker}((\widetilde {Sq^0_*})_{n_1})\cap (QP^{\otimes 6}_{n_1})^{>0}$ are the $\mathbb F_2$-vector subspaces of ${\rm Ker}((\widetilde {Sq^0_*})_{n_1}$ and $QP^{\otimes 6}_{n_1}\cong (QP^{\otimes 6}_{n_1})^{0}\bigoplus (QP^{\otimes 6}_{n_1})^{>0},$ one gets  
$$ QP^{\otimes 6}_{n_1}  \cong  (QP^{\otimes 6}_{n_1})^{0}\bigoplus \big({\rm Ker}((\widetilde {Sq^0_*})_{n_1})\cap (QP^{\otimes 6}_{n_1})^{>0}\big)\bigoplus QP^{\otimes 6}_{n_0}.
$$
Through the combination of the previously mentioned remark along with the utilization of Corollary \ref{hq10-1}, we arrive at the following conclusion.

\begin{corls}\label{hqMKR}
We have an isomorphism of $\mathbb F_2$-vector spaces: $$ QP_{n_0}^{\otimes 6}\cong \langle \{[t_1t_2\ldots t_6t^{2}]:\, t\in \mathscr C^{\otimes 6}_{n_0}\} \rangle \cong \bigoplus_{1\leq j\leq 5}(QP^{\otimes 6}_{n_1})^{>0}(\overline{\omega}^{(j)}),$$
where $\overline{\omega}^{(1)}:= (6,2,2,1),$\ $\overline{\omega}^{(2)}:=(6,2,4),$\ $\overline{\omega}^{(3)}:= (6,4,1,1),$\ $\overline{\omega}^{(4)}:=(6,4,3)$ and $\overline{\omega}^{(5)}:=(6,6,2),$ and the dimension of $(QP^{\otimes 6}_{n_1})^{>0}(\overline{\omega}^{(j)}) $ is determined by $$\dim (QP^{\otimes 6}_{n_1})^{>0}(\overline{\omega}^{(j)}) = \dim (QP^{\otimes 6}_{n_0})^{0}(\overline{\omega}^{(j)})+\dim (QP^{\otimes 6}_{n_0})^{>0}(\overline{\omega}^{(j)}),\ \mbox{for $1\leq j\leq 5.$}$$ Here the dimensions of $(QP^{\otimes 6}_{n_0})^{0}(\overline{\omega}^{(j)})$ and $(QP^{\otimes 6}_{n_0})^{>0}(\overline{\omega}^{(j)})$ are given as in Corollary \ref{hq10-1}.
\end{corls}
\end{rems}

We must now determine the dimensions of $(QP^{\otimes 6}_{n_1})^{0}$ and ${\rm Ker}((\widetilde {Sq^0*})_{n_1})\cap (QP^{\otimes 6}_{n_1})^{>0}.$ To accomplish this, we invoke a well-known outcome concerning the dimension of $QP^{\otimes 5}$ at degree $n_1$.

\begin{therm}[see Walker and Wood \cite{Walker-Wood2}]\label{dlWW}
In any minimal generating set for the $\mathcal A$-module $P^{\otimes h},$ there are $2^{\binom{h}{2}}$ elements in degree $2^{h}-(h+1).$ Consequently, $QP^{\otimes 5}_{n_1}$ is an $\mathbb F_2$-vector of dimension $1024.$
\end{therm}

Walker and Wood proved this by considering the special case of the Steinberg representation ${\rm St}_h,$ using the hook formula to count the number of semistandard Young tableaux. More precisely, they claim that by the hook formula, the dimension of the cohit module $QP^{\otimes h}_{2^{h}-h-1}$ is upper bounded by $2^{\binom{h}{2}}.$ The equality then follows from the first occurrence of the Steinberg representation in this degree. Thus, $QP^{\otimes h}_{2^{h}-h-1}\cong {\rm St}_h$ for the first occurrence degree $2^{h}-h-1.$ It would also be interesting to see that the dimension of this cohit module is equal to the order of the Borel subgroup $B_h$ of $GL_h.$

\medskip

We also employ the following homomorphisms: Let $h$ be a fixed integer with $5\leq h\leq 6$, and for each $l\in\mathbb Z$ such that $1\leq l\leq h$, we define a homomorphism $\mathsf{q}_{l}: P^{\otimes (h-1)}\longrightarrow P^{\otimes h}$ of algebras by setting $\mathsf{q}_{l}(t_j) = t_j$ for $1\leq j \leq l-1$ and $\mathsf{q}_{l}(t_j) = t_{j+1}$ for $l\leq j \leq h-1.$ Obviously, this $\mathsf{q}_{l}$ is also a homomorphism of $\mathcal A$-modules. The following comment is a crucial factor in computing $(QP^{\otimes 6}_{n_1})^{0}$ and ${\rm Ker}((\widetilde {Sq^0*})_{n_1})\cap (QP^{\otimes 6}_{n_1})^{>0}.$

\begin{rems}\label{nxp1}
\begin{itemize}
\item[(i)] It is patently obvious that the weight vector of the minimal spike $t_1^{15}t_2^{7}t_3^{3}t_4$ of degree $n_1$ in the $\mathcal A$-module $P^{\otimes 6}$ is $(4,3,2,1),$ In \cite{MKR}, Mothebe et.al proved that the cohit $QP^{\otimes 6}$ has dimension $1205$ in degree $11.$ So, $\omega(\underline{t})\in\big\{(3,2,1),\, (3,4),\, (5,1,1),\, (5,3)\big\}.$ As an immediate consequence of these results and Theorems \ref{dlKS}(i), \ref{dlSin}, we state that if $t$ is an admissible monomial in $P^{\otimes 6}_{n_1}$ such that $[t]$ belongs to the kernel of Kameko's map $(\widetilde {Sq^0_*})_{n_1},$ then $\omega(t)\in \big\{(4, 3,2,1),\ (4, 3,4),\ (4, 5,1,1),\ (4, 5,3)\big\}$ and $t$ can represented as $t_it_jt_kt_{\ell}\underline{t}^2,$ where $\underline{t}$ is an admissible monomial of degree $11$ in $P^{\otimes 6}$ and $1\leq i<j<k<\ell\leq 6.$ 

\item[(ii)] 
Since $QP^{\otimes 6}_{n_1}(\omega) \cong (QP^{\otimes 6}_{n_1})^{0}(\omega)\bigoplus (QP^{\otimes 6}_{n_1})^{> 0}(\omega),$ where $\omega$ is a weight vector of degree $n_1,$ one obtains an isomorphism: $QP^{\otimes 6}_{n_1} \cong (QP^{\otimes 6}_{n_1})^{0}\bigoplus \big(\bigoplus_{\deg(\omega)=n_1}(QP^{\otimes 6}_{n_1})^{>0}(\omega)\big).$ On the other hand, by Remark \ref{nxp0} and Corollary \ref{hqMKR}, we infer that $$QP^{\otimes 6}_{n_1} \cong  (QP^{\otimes 6}_{n_1})^{0}\bigoplus \big({\rm Ker}((\widetilde {Sq^0_*})_{n_1})\cap (QP^{\otimes 6}_{n_1})^{>0}\big)\bigoplus \big(\bigoplus_{1\leq j\leq 5}(QP^{\otimes 6}_{n_1})^{>0}(\omega^{(j)})\big),$$ where $(QP^{\otimes 6}_{n_1})^{0}\subset {\rm Ker}((\widetilde {Sq^0_*})_{n_1}).$  Hence, by invoking the aforementioned argument (i), an isomorphism will be established as follows: ${\rm Ker}((\widetilde {Sq^0_*})_{n_1})\cap (QP^{\otimes 6}_{n_1})^{>0}\cong  U_1\bigoplus U_2,$ where $$U_1:=(QP^{\otimes 6}_{n_1})^{>0}(4,5,1,1)\bigoplus (QP^{\otimes 6}_{n_1})^{>0}(4,5,3),\  \mbox{ and }\ U_2:=  (QP^{\otimes 6}_{n_1})^{>0}(4,3,2,1)\bigoplus (QP^{\otimes 6}_{n_1})^{>0}(4,3,4).$$
\end{itemize}
\end{rems}

Drawing on Remark \ref{nxp1}, we establish the second main result of this paper.

\begin{therm}\label{dlc1}
With the above notation, the following assertions are true:

\begin{itemize}
\item[(i)]
$ \dim (QP^{\otimes  6}_{n_1})^{0}(\omega) = \left\{\begin{array}{ll}
5184&\mbox{if $\omega = (4,3,2,1)$},\\[1mm]
0&\mbox{if $\omega\neq (4,3,2,1)$}.
\end{array}\right.$\\[1mm]
Consequently,  $(QP^{\otimes  6}_{n_1})^{0}$ is isomorphic to $(QP^{\otimes  6}_{n_1})^{0}(4,3,2,1)$ and $(\mathscr C^{\otimes 6}_{n_1})^{0} = (\mathscr C^{\otimes 6}_{n_1})^{0}(4,3,2,1)$ has all $5184$ admissible monomials.

\item[(ii)] $\dim U_1 = 546$ and $\dim U_2 = 3090.$ These imply that there exist exactly $9765$ admissible monomials of degree $n_1$ in the $\mathcal A$-module $P^{\otimes 6}.$ Consequently, the cohit $QP^{\otimes 6}_{n_1}$ is $9765$-dimensional.
\end{itemize}
\end{therm}

We will now recall a previously established result on the Kameko squaring operation.

\begin{therm}[see Kameko \cite{Kameko}]\label{dlMK}
The homomorphism $(\widetilde {Sq^0_*})_{2n+h}: QP^{\otimes h}_{2n+h}  \longrightarrow QP^{\otimes h}_{n}$ is an isomorphism of the $\mathbb F_2$-vector spaces if and only if $\mu(2n+h) = h.$ Then, one has an inverse homomorphism $ \psi: QP^{\otimes h}_{n}\longrightarrow QP^{\otimes h}_{2n+h}$ of $(\widetilde {Sq^0_*})_{2n+h},$ which is induced by the mapping $\psi: P^{\otimes h}\longrightarrow P^{\otimes h},$ $t\longmapsto \prod_{1\leq j\leq h}t_jt^{2}.$
\end{therm}

Write $\mathbb F_q$ for the Galois field of size $q$ ($q$ being a power of the prime characteristic $p$ of this field), let $B_h(\mathbb F_q)$ be the Borel subgroup of the general linear group $GL_h(\mathbb F_q)$ over $\mathbb F_q.$ When $q = 2,$ we put $GL_h:= GL_h(\mathbb F_2)$ and $B_h:= B_h(\mathbb F_2).$ Note that  $\mathbb F_{q = p^{m}}\cong \mathbb F_p^{\oplus m}$ as groups (in fact as $\mathbb F_p$-modules). For the sake of completeness, let us remind the readers that the algebra of Steenrod $q$-th reduced powers $\pmb {A}_q$ can be defined as an algebra over $\mathbb F_q$ by generators $\mathscr P^{j},\, j\geq 0,$ subject to the relation $\mathscr P^{0} = 1$ and the Adem relations, $\mathscr P^{a}\mathscr P^{b} = \sum_{0\leq j\leq [a/q]}(-1)^{i+j}\binom{(q-1)(b-j)-1}{a-qj}\mathscr P^{a+b-j}\mathscr P^{j},\ a<qb.$ For $q = p,$ as a subalgebra of the mod $p$ Steenrod algebra $\mathcal A_p,$ the element $\mathscr P^{k}$ is given the degree $2k(p-1),$ but for simplicity, one regrade $\pmb{A}_q$ by giving $\mathscr P^{j}$ the "reduced" degree $k.$ So, when $q = p =2,$ $\mathscr P^{k}$ will mean Steenrod squares $Sq^{k},$ and not $Sq^{2k}$ (see also \cite{Smith}). Consider an $h$-dimensional vector space $\pmb{V}_h$ over $\mathbb F_q,$ the symmetric power algebra $S(\pmb{V}^{*}_h)$ on the dual $\pmb{V}_h^{*} = {\rm Hom}_{\mathbb F_q}(\pmb{V}_h, \mathbb F_q)$ of $\pmb{V}_h$ is identified with the polynomial algebra $\mathbb F_q[t_1, \ldots, t_h],$ where $\deg(t_i) = 1$ for every $i$ and $\{t_1, \ldots, t_h\}$ is a basis of $\pmb{V}^{*}_h.$ Applying Theorem \ref{dlMK} in conjunction with the work by Hai \cite{Hai}, we derive the following corollaries.

\begin{corls}
In degree $q^{h-1}-h,$ we have $$ \dim_{\mathbb F_q} (\mathbb F_q[t_1, \ldots t_h]/\overline{\pmb{A}}_q\mathbb F_q[t_1, \ldots t_h])_{q^{h-1}-h} = {\rm ord}(GL_{h-1}(\mathbb F_q)/B^{*}_{h-1}(\mathbb F_q)) = \prod_{1\leq j\leq h-1}(q^{j}-1),$$
on which $B_{h-1}^{*}(\mathbb F_q)\subset B_{h-1}(\mathbb F_q)\cap {\rm Ker}(\det)$ and each element of $B_{h-1}^{*}(\mathbb F_q)$ has 1's in the main diagonal. Here $\det$ denotes the $\mathbb F_q$-linear map $GL_{h-1}(\mathbb F_q)\longrightarrow \mathbb F_q^{*}.$
\end{corls}

\begin{corls}\label{hqs22}
Let $h \geq 6$ be a given fixed integer. Setting $n_{h, s}:=2^{s+4}-h,$ then, for each $s\geq h-5,$ we have
$$ \dim QP^{\otimes h}_{n_{h, s}}  =  {\rm ord}(GL_{h-1}/B_{h-1})  =  \prod_{1\leq j\leq h-1}(2^{j}-1).$$
Moreover, $QP^{\otimes h}_{n_{h, s}}\cong {\rm Ker}((\widetilde {Sq^0_*})_{n_{h, h-5}})\bigoplus QP^{\otimes h}_{2^{h-2}-h}$ for any $s\geq h-5.$
\end{corls}

Indeed, we have that the order of the Borel subgroup $B_h(\mathbb F_q)$ is $q^{\binom{h}{2}}\prod_{1\leq j\leq h}(q-1),$ since elements in the main diagonal are taken from $\mathbb F_q^{*}$ and elements above to the main diagonal can be any element of $\mathbb F_q.$ The order of $GL_h(\mathbb F_q)$ is determined as follows: the first row $u_1$ of the matrix can be anything but the $0$-vector, so there are $q^{h}-1$ possibilities for the first row. For any one of these possibilities, the second row $u_2$ can be anything but a multiple of the first row, giving $q^{h}-q$ possibilities. For any choice $u_1,\, u_2$ of the first two rows, the third row can be anything but a linear combination of $u_1$ and $u_2.$ The number of linear combinations $\sum_{1\leq i\leq 2}\gamma_iu_i$ is just the number of choices for the pair $(\gamma_1,\gamma_2),$ and there are $q^{2}$ of these. It follows that for every $u_1$ and $u_2,$ there are $q^{h}-q^{2}$ possibilities for the third row. For any allowed choice  $u_1, u_2, u_3,$ the fourth row can be anything except a linear combination $\sum_{1\leq i\leq 3}\gamma_iu_i$ of the first three rows. Thus for every allowed $u_1, u_2, u_3$ there are $q^{3}$ forbidden fourth rows, and therefore $q^{h}-q^{3}$ allowed fourth rows. In the same way, the number of non-singular matrices is $(q^{h}-1)(q^{h}-q)\ldots (q^{h}-q^{h-1}),$  and so, $${\rm ord}(GL_h(\mathbb F_q)) = \prod_{0\leq j\leq h-1}(q^{h}-q^{j}) = q^{\binom{h}{2}}\prod_{1\leq j\leq h}(q^{j}-1).$$ 
Given the $\mathbb F_q$-linear $\det: GL_h(\mathbb F_q)\longrightarrow \mathbb F_q^{*},$ consider the subsets $B^{*}_h(\mathbb F_q)$ of the groups $B_h(\mathbb F_q)\cap {\rm Ker}(\det),$ where each element of $B^{*}_h(\mathbb F_q)$ has 1's in the main diagonal. Then, $B^{*}_h(\mathbb F_q)$ is also a self-conjugate subgroup of $GL_h(\mathbb F_q).$ It is straightforward to see that the order of $B^{*}_h(\mathbb F_q)$ is $q^{\binom{h}{2}}.$ In particular, when $q = 2,$ we have ${\rm Ker}(\det) = GL_h$ and $B^{*}_h = B_h.$ Thus ${\rm ord}(B_h(\mathbb F_q)) = {\rm ord}(B^{*}_h(\mathbb F_q))\prod_{1\leq j\leq h}(q-1)$ and ${\rm ord}(GL_h(\mathbb F_q)) = {\rm ord}(B^{*}_h(\mathbb F_q))\prod_{1\leq j\leq h}(q^{j}-1).$ In \cite[Theorem 4]{Hai}, by considering a variant of a family of finite quotient rings of $\mathbb F_q[t_1, \ldots, t_h]$, Hai proved that the space of the indecomposable elements of $\mathbb F_q[t_1, \ldots, t_h]$ has dimension $(q-1)(q^{2}-1)\ldots (q^{h-1}-1)$ in degree $q^{h-1}-h.$ From these data, we get 
$$  \dim_{\mathbb F_q}(\mathbb F_q[t_1, \ldots t_h]/\overline{\pmb{A}}_q\mathbb F_q[t_1, \ldots t_h])_{q^{h-1}-h}  = \prod_{1\leq j\leq h-1}(q^{j}-1) = {\rm ord}(GL_{h-1}(\mathbb F_q)/B^{*}_{h-1}(\mathbb F_q)).$$ 
(The reader should also keep in mind that the product $\prod_{1\leq j\leq h-1}(q^{j}-1)$ is also a well known formula for the degree of a cuspidal character of $GL_h(\mathbb F_q).$ The cuspidal characters are of great importance for characters of $GL_h(\mathbb F_q)$ since each character of this linear group is build up from cuspidal characters.)
 Now, with the field $\mathbb F_2$ and degree $n_{h, s} = 2^{s+4}-h,$ since
$$n_{h, s} = (2^{s+3}-1) + (2^{s+2}-1)  + (2^{s+1}-1) +  \cdots + (2^{s-(h-5)}-1)+(2^{s-(h-5)}-1),$$ $\mu(n_{h, s}) = h$ for any $s \geq h-4,$ and so, by Theorem \ref{dlMK}, the iterated Kameko squaring operation $(\widetilde {Sq^0_*})^{s-h+5}_{n_{h, s}}: QP^{\otimes h}_{n_{h, s}}  \longrightarrow QP^{\otimes h}_{n_{h, h-5}}$ is an isomorphism for every $s \geq h-5.$ Combining this with the facts that $\dim QP^{\otimes h}_{n_{h, h-5}}  = \prod_{1\leq j\leq h-1}(2^{j}-1)$ and $B_{h-1}^{*} = B_{h-1},$ we must have
$$ \dim QP^{\otimes h}_{n_{h, s}} = \prod_{1\leq j\leq h-1}(2^{j}-1) = {\rm ord}(GL_{h-1}/B^{*}_{h-1})= {\rm ord}(GL_{h-1}/B_{h-1}),\ \mbox{for all $s\geq h-5.$}$$ Moreover, as the Kameko homomorphism $(\widetilde {Sq^0_*})_{n_{h, s}}: QP^{\otimes h}_{n_{h, s}}\longrightarrow QP^{\otimes h}_{n_{h, s-1}}$ is an epimorphism and $QP^{\otimes h}_{n_{h, s}}\cong QP^{\otimes h}_{n_{h, h-5}},$ one gets $QP^{\otimes h}_{n_{h, s}}\cong {\rm Ker}((\widetilde {Sq^0_*})_{n_{h, h-5}})\bigoplus QP^{\otimes h}_{2^{h-2}-h}$ for arbitrary $s\geq h-5.$ 

\medskip

Let us take notice that, in the case of $q=2$ and $h=6$, the dimensionality of $QP^{\otimes 6}_{n_1}$ is equal to $(2^{1}-1)\ldots (2^{6-1}-1) = 9765$, a result that can be gleaned from Theorem \ref{dlc1}. Therefore, our research stands independently of Hai's, and our approach is completely distinct. Furthermore our work offers a precise and unambiguous description of a monomial basis for the cohit module $QP^{\otimes 6}_{2^{6-1}-6= n_1}$, which serves as a representation of $GL_6$. Theoretically, our technique can be extended to any values of $h$ and $n.$ Nonetheless, the process of calculation becomes increasingly complex as the dimensions of $QP^{\otimes h}_n$ grow larger with increasing $h$ and $n$.

\begin{rems}
Consider general degree $n_h = 2^{h-2}-h,\, h\geq 4,$ we have $\dim QP^{\otimes 4}_{n_4} = \dim \mathbb F_2 = 1$, $\dim QP^{\otimes 5}_{n_5} = 7$ (see Theorem \ref{dlMM}) and $\dim QP^{\otimes 6}_{n_6} = 945$ (see Theorem \ref{dlMKR}). Given any $h\geq 7,$ by Corollary \ref{hqs22}, $\dim QP^{\otimes h}_{n_h} = 3.7\ldots (2^{h-1}-1)-\dim {\rm Ker}((\widetilde {Sq^0_*})_{n_{h, h-5}}).$ Hence, in order to determine the dimension of $QP^{\otimes h}_{n_h}$ for all $h>6,$ it suffices to calculate the dimension of the kernel of Kameko's map $(\widetilde {Sq^0_*})_{n_{h, h-5}}$. However, this aspect will be investigated in a separate study. Utilizing a result from Hai \cite[Corollary 3]{Hai}, we have $QP^{\otimes (h-2)}_{n_h}\cong {\rm St}_{h-2}\otimes_{\mathbb F_2}{\rm det}^{1},$ where ${\rm det}^{1}$ denotes the first power of the determinant representation of $GL_{h-2}$ and ${\rm St}_{h-2}$ is the Steinberg module (a.k.a the Steinberg representation). Remarkably, for $h = 8,\, 9,$ since the cohomology groups ${\rm Ext}_{\mathcal A}^{h-2, n_h + h-2}(\mathbb F_2, \mathbb F_2)$ are trivial \cite{Bruner}, the Singer conjecture is wrong if $\dim[QP^{\otimes (h-2)}_{n_h}]^{GL_{h-2}} > 0.$ For $h  =4,$ one has an isomorphism $(\mathbb F_2 \otimes_{GL_2} {\rm Ann}_{\overline{\mathcal A}}[P^{\otimes 2}]^{*})_{0} \cong \mathbb F_2\cong {\rm Ext}_{\mathcal A}^{2, 2}(\mathbb F_2, \mathbb F_2),$ which implies that the Singer conjecture holds for bidegree $(2, 2).$ For $h = 5,\, 6,\, 7,$ Singer's conjecture for bidegree $(h-2, n_h+h-2)$ has been verified by Boardman \cite{Boardman} for $h = 5,$ by the present author \cite{Phuc12} for $h = 6$ and by Sum \cite{Sum3} for $h = 7.$ By these, it would also be of significant interest to determine explicit generators of $QP^{\otimes (h-2)}_{n_h}.$ The dimension of this cohit module was determined by Peterson \cite{Peterson} for $h = 4,$ by Kameko \cite{Kameko} for $h  = 5,$ and by Sum \cite{Sum2, Sum3} for $h = 6,\, 7.$ (See also Theorems \ref{dlMKR} and \ref{dlMM} for the cases where $4\leq h\leq 6.$)  
\end{rems}

Building upon Corollary \ref{hqs22} and the calculations in \cite{Tangora, Bruner, Bruner2, Lin2}, we can see that with degree $n_{h, s}$ as in Corollary \ref{hqs22}, 
 $$ \begin{array}{ll}
 {\rm Ext}_{\mathcal A}^{7, 7+n_{7, s}}(\mathbb F_2, \mathbb F_2)=\left\{\begin{array}{ll}
0 &\mbox{if $s = 1,\, 4$},\\[1mm]
\langle Q_2(0)\rangle &\mbox{if $s = 2$},\\[1mm]
\langle \{Q_2(1), h_6D_2 \} \rangle &\mbox{if $s = 3$},\\[1mm]
\end{array}\right.\\
\\
 {\rm Ext}_{\mathcal A}^{8, 8+n_{8, s}}(\mathbb F_2, \mathbb F_2)=\left\{\begin{array}{ll}
0 &\mbox{if $s = 1,\, 2$},\\[1mm]
\langle h_{6}Q_2(0)\rangle &\mbox{if $s = 3$},\\[1mm]
\langle x_{n_{8, 4}, 8}\rangle &\mbox{if $s = 4$},
\end{array}\right.
\end{array}$$
where $x_{n_{8, 4}, 8}$ is an indecomposable element. We believe that the following prediction would be of significant interest to investigate regarding Conjecture \ref{gtSinger} in high homological degrees.
\newpage
\begin{conj}\label{gtbsm}
The family $\big\{Q_2(k):\, k\geq 0\big\}$ is a finite $Sq^{0}$-family. Furthermore, we have that:
\begin{itemize}
\item[(i)] the transfer $Tr_7^{\mathcal A}$ does not detect the non-zero elements $Q_2(0),\, Q_2(1)$ and $h_6D_2;$  

\item[(ii)] the transfer $Tr_8^{\mathcal A}$ does not detect the non-zero elements $h_{6}Q_2(0)$ and $x_{n_{8, 4}, 8}.$  
\end{itemize}
\end{conj}


Note that an $Sq^{0}$ -family is called \textit{finite} if it has only finitely many nonzero elements, \textit{infinite} if all of its elements are nonzero \cite{Hung}. Due to Corollary \ref{hqs22}, it is observed that the conjecture for items (i) and (ii) is valid under the following circumstances: 
$$ \begin{array}{ll}
 (\mathbb F_2 \otimes_{GL_7} {\rm Ann}_{\overline{\mathcal A}}[P^{\otimes 7}]^{*})_{n_{7, 1}}\cong  (\mathbb F_2 \otimes_{GL_7} {\rm Ann}_{\overline{\mathcal A}}[P^{\otimes 7}]^{*})_{n_{7, 2}} \cong [{\rm Ker}((\widetilde {Sq^0_*})_{n_{7, 2}})]^{GL_7} = 0,\\
 (\mathbb F_2 \otimes_{GL_8} {\rm Ann}_{\overline{\mathcal A}}[P^{\otimes 8}]^{*})_{n_{8, 2}}\cong (\mathbb F_2 \otimes_{GL_8} {\rm Ann}_{\overline{\mathcal A}}[P^{\otimes 8}]^{*})_{n_{8, 3}} \cong [{\rm Ker}((\widetilde {Sq^0_*})_{n_{8, 3}})]^{GL_8} = 0.
\end{array}$$  
Our approach for determining the domain of $Tr_7^{\mathcal A}$ with respect to degrees $n_{7, 1}$ and $n_{7, 2}$, as well as the domain of $Tr_8^{\mathcal A}$ with respect to degree $n_{8, 3}$, will involve the use of Theorems \ref{dlMM} and \ref{dlWW}, alongside Corollary \ref{hqs22}. Nevertheless, the calculation at hand seems to be rather daunting.

Adopting an alternative perspective, H\uhorn ng \cite{Hung} proposed an interesting notion concerning a \textit{critical element} that exists within ${\rm Ext}_{\mathcal A}^{h, h+n}(\mathbb F_2, \mathbb F_2)$. Specifically, a non-zero element $\zeta$ in ${\rm Ext}_{\mathcal A}^{h, h+n}(\mathbb F_2, \mathbb F_2)$ is deemed \textit{critical} if it satisfies two conditions: firstly, $\mu(2n+h) = h$, and secondly, the image of $\zeta$ under the classical squaring operation $Sq^{0}$ is zero.
It is well-established that $Sq^0$ is a monomorphism in positive stems of ${\rm Ext}_{\mathcal A}^{h, *}(\mathbb F_2, \mathbb F_2)$ for $h < 5,$  thereby implying the absence of any critical element for $h < 5.$ Remarkably, H\uhorn ng's work \cite[Theorem 5.9]{Hung} states that Singer's Conjecture \ref{gtSinger} is not valid, if the algebraic transfer detects the critical elements. Now, given the non-zero elements $Q_2(1)$ and $h_6D_2$, we are able to deduce that $\mu(2{\rm Stem}(Q_2(1))+7) = \mu(2{\rm Stem}(h_6D_2)+7)  = 7.$ Furthermore, it is worth noting that $Sq^{0}(Q_2(1)) = 0 = Sq^{0}(h_6D_2)$, an observation which can be attributed to the fact that ${\rm Ext}_{\mathcal A}^{7, 7+n_{7,4}}(\mathbb F_2, \mathbb F_2) = 0$, as previously discussed. Thus $Q_2(1)$ and $h_6D_2$ must be critical elements. By this reason, in the event that Conjecture \ref{gtbsm}(i) is proven to be false,  it would entail the refutation of Singer's conjecture in general.

\medskip

We now turn our attention to the hit problem for $P^{\otimes 6}$ in degree $n_s$ with $s > 1.$

\medskip

{\bf Cases $\pmb{s > 1}.$}\ By Theorem \ref{dlc1} and Corollary \ref{hqs22}, $\dim QP^{\otimes 6}_{n_s} = \dim QP^{\otimes 6}_{n_1}  = 9765$ for any $s > 0.$ Moreover, since the iterated homomorphism $((\widetilde {Sq^0_*})_{n_{s}})^{s-1}: QP^{\otimes 6}_{n_s}  \longrightarrow QP^{\otimes 6}_{n_1}$ is an isomorphism, for every positive integer $s,$ we have immediately the below corollary.

\begin{corls}\label{dlc2}
For each integer $s\geq 2,$ the set $\big\{[t]:\ t\in \psi^{s-1}(\mathscr C^{\otimes 6}_{n_1})\big\}$ is a monomial basis of the $\mathbb F_2$-vector space $QP^{\otimes 6}_{n_s},$ on which $\psi: P^{\otimes 6}\longrightarrow P^{\otimes 6},\, t\longmapsto t_1t_2\ldots t_6t^{2}$ and  $\psi^{s-1}(\mathscr C^{\otimes 6}_{n_1}) = \bigg\{\prod_{1\leq j\leq 6}t^{2^{s-1}-1}_ju^{2^{s-1}}:\ u\in \mathscr C^{\otimes 6}_{n_1}\bigg\}.$
\end{corls}

The next contribution of this work is to apply the aforementioned results into the investigation of the cohit module $QP^{\otimes 7}$ in general degree $n_{s+5}$ and the behavior of the sixth algebraic transfer in internal degrees $n_s$ for any $s > 0.$ To achieve this goal, we will begin by recalling an interesting result on an inductive formula for the dimension of $QP^{\otimes h}_n.$

\begin{therm}[see Sum \cite{Sum2}]\label{dlS}
Consider the degree $n$ of the form \eqref{pt} with $k = h-1,$ and $s, r$ positive integers such that $1\leq h-3\leq \mu(r)\leq h-2,$ and $\mu(r) = \alpha(r + \mu(r)).$ Then for each $s\geq h-1,$ we have $\dim QP^{\otimes h}_n  = (2^{h}-1)\dim QP^{\otimes (h-1)}_r.$
\end{therm}

\begin{rems}\label{nxP}

With the general degrees $n_s:= (h-1)(2^{s}-1) + r\cdot 2^{s},$ assume there is a non-negative integer $\zeta$ such that $\zeta < s$ and $1\leq h-3\leq \mu(n_{\zeta}) = \alpha(n_{\zeta}+ \mu(n_{\zeta}))\leq h-2.$ Let us consider generic degrees of the form $k(2^{s-\zeta + h-1}-1) + r\cdot 2^{s-\zeta + h-1},$ where $k = h-1$, $r = n_{\zeta}$ and $s\geq \zeta\geq 0.$ Consequently, due to $\mu(r) = \alpha(r+\mu(r)),$ we have the following inductive formula, which is deduced from Theorem \ref{dlS} and the proof of this theorem on pages 445-446 of \cite{Sum2}:
$$ \dim QP^{\otimes h}_{(h-1)(2^{s-\zeta + h-1}-1) + n_{\zeta}2^{s-\zeta + h-1}} = (2^{h}-1)\dim QP^{\otimes (h-1)}_{n_s},\ \mbox{ for every $s\geq \zeta.$}$$
Now, with $h = 7,$ $r = 10$, $\zeta  = 1,$ and degree $n_s,$ we have $\mu(n_{1}) = 4 = \alpha(n_{1}+ \mu(n_{1})).$ Hence the following is immediate from Corollary \ref{hqs22} and Remark \ref{nxP}.
\end{rems}

\begin{corls}\label{dlc3}
For every positive integer $s,$ the cohit module $QP^{\otimes 7}$ has dimension $1240155$ in degree $n_{s+5} = 6(2^{s+5}-1) + 10\cdot 2^{s+5}.$
\end{corls}

As a consequence of Theorem \ref{dlc1} and the computations done in \cite{Tangora, Bruner, Bruner2}, we are able to establish the third main result of this paper.

\begin{therm}\label{dlc4}
For each integer $s > 0,$ the coinvariant $(\mathbb F_2 \otimes_{GL_6} {\rm Ann}_{\overline{\mathcal A}}[P^{\otimes 6}]^{*})_{n_s}$ is trivial. Consequently, the algebraic transfer $Tr_6^{\mathcal A}: (\mathbb F_2 \otimes_{GL_6} {\rm Ann}_{\overline{\mathcal A}}[P^{\otimes 6}]^{*})_{n_s}\longrightarrow {\rm Ext}_{\mathcal A}^{6, 6+n_s}(\mathbb F_2, \mathbb F_2)$ is a monomorphism, but it is not an epimorphism for $0< s < 3.$ This means that the transfer $Tr_6^{\mathcal A}$ does not detect the non-zero elements $h_4Ph_2 = h_2^{2}g_1\in {\rm Ext}_{\mathcal A}^{6, 6+n_1}(\mathbb F_2, \mathbb F_2)$ and $D_2\in {\rm Ext}_{\mathcal A}^{6, 6+n_2}(\mathbb F_2, \mathbb F_2).$ When $s = 3,$ the transfer $Tr_6^{\mathcal A}$ is a trivial isomorphism. 
\end{therm}


By invoking Theorems \ref{dlc0} and \ref{dlc4}, it is possible to deduce that in the bidegrees $(6, 6+n_s),$ Singer's transfer is a monomorphism, but not an epimorphism for $0\leq s\leq 2.$ In the case where $s\geq 3,$ the transfer is a trivial isomorphism if its codomain is zero, and a monomorphism otherwise. These lead to an immediate consequence

\begin{corls}\label{hqbs}
Conjecture \ref{gtSinger} is valid in the bidegrees of $(6, 6+n_s)$ for any non-negative integer $s$.
\end{corls}

{\bf Final remarks.} Drawing upon the findings in \cite{Tangora, Bruner, Bruner2, Lin, Chen, Chen2, Lin2} on the structure of ${\rm Ext}_{\mathcal A}^{s, *}(\mathbb F_2, \mathbb F_2)$ for $s\leq 6,$ we end this section by presenting the following conjecture, which predicts the structure of the sixth cohomology group ${\rm Ext}_{\mathcal A}^{6, 6+n_s}(\mathbb F_2, \mathbb F_2)$ for all $s\geq 0$. 

\begin{conj}\label{gtP}
With degree $n_s = 6(2^{s}-1) + 10\cdot 2^{s},$ we have 
\begin{equation*}
 {\rm Ext}_{\mathcal A}^{6, 6+n_s}(\mathbb F_2, \mathbb F_2)=\left\{\begin{array}{ll}
\langle h_1Ph_1 \rangle &\mbox{if $s = 0$},\\[1mm]
\langle h_2^{2}g_1\rangle &\mbox{if $s = 1$},\\[1mm]
\langle D_2  \rangle &\mbox{if $s = 2$},\\[1mm]
\langle \{h_s^{2}h_{s+1}^{2}h_{s+2}^{2},\, h_{s+1}^{2}g_s= h_{s+1}h_{s+3}g_{s-1},\, h_sh_{s+2}f_{s-1},\, h_sh_{s+1}^{2}c_s\}  \rangle = 0 &\mbox{if $s \geq 3$},\\[1mm]
\end{array}\right.
\end{equation*}
where $P$ denotes the Adams periodicity operator and
\begin{equation*}
 \begin{array}{ll}
\medskip
 &h_1Ph_1=[\lambda_1\lambda_2 \lambda_0^3 \lambda_7 + \lambda_1^2 \lambda_2 \lambda_4 \lambda_1^2 + \lambda_1^3 \lambda_2 \lambda_4 \lambda_1 + \lambda_1^2 \lambda_2 \lambda_1^2 \lambda_4]\neq 0,\\
& h_2^{2}g_1= Sq^{0}(h_1Ph_1) =  h_4Ph_2= [\lambda_{15}\lambda_4\lambda_0^3 \lambda_7 + \lambda_{15}\lambda_3\lambda_5 \lambda_1^3 + \lambda_{15}\lambda_3\lambda_2\lambda_4 \lambda_1^2 +  \lambda_{15}\lambda_3\lambda_1\lambda_2 \lambda_4 \lambda_1 \\
\medskip
&\hspace{5cm} + \lambda_{15}\lambda_3\lambda_2\lambda_1^2\lambda_4 + \lambda_{15}\lambda_2\lambda_2 \lambda_0^2 \lambda_7 + \lambda_{15}\lambda_1\lambda_1 \lambda_2\lambda_0 \lambda_7]\neq 0,\\
&D_2 = [\lambda_0^{4}\lambda_{11}\lambda_{47}]\neq 0.
\end{array}
\end{equation*}
\end{conj}
Note that $h_{s+1}g_s = h_{s+3}g_{s-1}$ for any $s\geq 2.$
 Given the calculations presented in \cite{Tangora, Bruner, Lin2}, it has been unequivocally established that the conjecture holds for $s\leq 3.$ Additionally, if the conjecture is confirmed to be accurate in general, then Singer's algebraic transfer would be a trivial isomorphism in the bidegrees $(6, 6+n_s)$ for $s\geq 3$. The readers will observe that Singer's Conjecture \ref{gtSinger} in these bidegrees would be disproven if the dimension of the invariant $[QP^{\otimes 6}_{n_1}]^{GL_6}$ is equal to 1. However, as demonstrated in Theorem \ref{dlc4}, this eventuality did not transpire. Inspired by the calculations in \cite{Bruner2}, we are confident that Conjecture \ref{gtP} also holds true for all $s > 3$. 

From another perspective, by virtue of the calculations set forth in the works \cite{Lin, Chen2, Bruner2}, and by making use of the fundamental property that $Sq^{0}$ is an algebraic homomorphism, it follows that $Sq^{0}(h_2^{2}g_1) = Sq^{0}(h_4Ph_2)= h_5h_3g_1 = 0.$ On the other hand, in \cite{Bruner}, Bruner claimed that ${\rm Ext}_{\mathcal A}^{6, 6+n_3}(\mathbb F_2, \mathbb F_2)$ is trivial. So, $Sq^{0}$ must send the indecomposable element $D_2$ to zero. Thus, since $\mu(2n_1+6) = 6= \mu(2n_2+6),$ both $h_2^{2}g_1$ and $D_2$ are critical elements. (Additionally, considering the fact that $2^{4} = {\rm Stem}(Ph_2) +5 < 4({\rm Stem}(Ph_2))^2$ and that $Ph_2\in {\rm Ext}_{\mathcal A}^{5, 16}(\mathbb F_2, \mathbb F_2)$ is a critical element, as noted in \cite{Hung}, it follows that $h_2^{2}g_1=h_4Ph_2 $ is a critical element. This explanation further strengthens the aforementioned assertion.) An interesting observation from H\uhorn ng's paper \cite[Lemma 5.3]{Hung} is that the condition $2^{m}\geq \max\big\{4d^{2}, d+h\big\}$ is insufficient to identify all critical elements of the form $h_mx$ in ${\rm Ext}_{\mathcal A}^{h, *}(\mathbb F_2, \mathbb F_2)$ when $x$ is also a critical element and $d = {\rm Stem}(x).$ This can be inferred from the foregoing facts that $D_2\in {\rm Ext}_{\mathcal A}^{6, 6+n_2}(\mathbb F_2, \mathbb F_2)$ and $h_6D_2\in {\rm Ext}_{\mathcal A}^{7, 7+n_{7, 3}}(\mathbb F_2, \mathbb F_2)$ are critical elements, but $2^{6} < \max\big\{4({\rm Stem}(D_2))^{2},\, {\rm Stem}(D_2)+6\big\}.$ In summary, even thourgh the elements $h_4Ph_2$ and $D_2$ are critical, they cannot be detected by the algebraic transfer. (It should be brought to the attention of the readers that the research conducted by Qu\`ynh \cite{Quynh} demonstrates that the indecomposable element $Ph_2$ does not belong to the image of the fifth transfer.) This reinforces the conclusion that Conjecture \ref{gtSinger} continues to hold for the bidegrees $(6, 6+n_1)$ and $(6, 6+n_2),$ as established in Theorem \ref{dlc4} and Corollary \ref{hqbs}.

\section{Proofs of the main results}\label{s4}

This section is devoted to proving Theorems \ref{dlc0}, \ref{dlc1} and \ref{dlc4}. To begin with, we need the following homomorphisms and a helpful remark below. One should note that $V_h\cong \langle \{t_1, \ldots, t_h\}\rangle \subset  P^{\otimes h}.$ For $1\leq d\leq h,$  we define the $\mathbb F_2$-linear map $\sigma_d: V_h\longrightarrow V_h$ by setting  
$$ \left\{\begin{array}{ll}
\sigma_d(t_d) = t_{d+1},\\[1mm]
\sigma_d(t_{d+1}) = t_d,\\[1mm]
\sigma_d(t_i) = t_i,\ \mbox{for $i\neq d, d +1,\; 1\leq d\leq h-1,$}\\[1mm]
\sigma_h(t_1) = t_1 + t_2,\ \ \sigma_h(t_i) = t_i,\ \mbox{for $2\leq i \leq h.$}
\end{array}\right.$$
\medskip

Denote by $\Sigma_h\subset GL_h$ the symmetric group of degree $h.$ Then, $\Sigma_h$ is generated by the ones associated with $\sigma_1,\, \ldots, \sigma_{h-1}.$ For each permutation in $\Sigma_h$, consider corresponding permutation matrix; these form a group of matrices isomorphic to $\Sigma_h.$ Indeed, consider the following map $\Delta: \Sigma_h\longrightarrow \mathcal P_{h\times h},$ where the latter is the set of permutation matrices of order $h.$ This map is defined as follows: given $\sigma\in \Sigma_h,$ the $i$-th column of $\Delta(\sigma)$ is the column vector with a $1$ in the $\rho(i)$-th position, and $0$ elsewhere. It is easy to see that $\Delta(\rho)$ is indeed a permutation matrix, since a $1$ occurs in any position if and only if that position is described by $(\rho(i), i),$ for any $1\leq i\leq h.$ The map $\Delta$ is clearly multiplicative. (It is to be noted that because these are matrices, it is enough to show that each corresponding entry is equal. So let us take the entry $(i, j)$ of each matrix.) Then, $\Delta(\rho\circ\rho')_{ij}  =1$ if and only if $i = \rho\circ \rho'(j).$ Note also that by ordinary matrix multiplication, one has $(\Delta(\rho)\Delta(\rho'))_{ij} = \sum_{1\leq k\leq h}\Delta(\rho)_{ik}\Delta(\rho')_{kj}.$ Now, we know that $\Delta(\rho)_{ik} = 1$ only when $i = \rho(k).$ Similarly, $\Delta(\rho')_{kj} = 1$ only when $k = \rho'(j).$ Hence, their product is one precisely when both of these happen: $i = \rho(k),$ and $k = \rho'(j).$ If both these do not happen simultaneously, then whenever one of $\Delta(\rho)_{ik},\, \Delta(\rho')_{kj}$ is one of the other will be zero, so the whole sum will be zero. However, this is the same as saying that the sum is one exactly when $i = \rho\circ\rho'(j).$ This description matches with the description for $\Delta(\rho\circ \rho')_{ij}$ given earlier. Hence, entry by entry these matrices are the same. Therefore the matrices are the same, and hence $\Delta$ is a homomorphism between the two spaces, an isomorphism as it has trivial kernel and the sets are of the same cardinality. Thus, $GL_h\cong GL(V_h),$ and $GL_h$ is generated by the matrices associated with $\sigma_1, \ldots, \sigma_h.$ 

\medskip

Let $T = t_1^{a_1}t_2^{a_2}\ldots t_h^{a_h}$ be a monomial in $P_n^{\otimes h}$. Then, the weight vector $\omega(T)$ is invariant under the permutation of the generators $t_j,\ j = 1, 2, \ldots, h;$ hence $QP_n^{\otimes h}(\omega(T))$ also has a $\Sigma_h$-module structure. We see that the linear map $\sigma_d$ induces a homomorphism of $\mathcal A$-algebras which is also denoted by $\sigma_d: P^{\otimes h}\longrightarrow P^{\otimes h}.$ So, a class $[T]_{\omega(T)}\in QP_n^{\otimes h}(\omega)$ is an $GL_h$-invariant if and only if $\sigma_d(T) \equiv_{\omega(T)} T$ for $1\leq d\leq h.$ If  $\sigma_d(T) \equiv_{\omega(T)} T$ for $1\leq d\leq h-1,$ then $[T]_{\omega(T)}$ is an $\Sigma_h$-invariant. (We must stress that the explicit calculation of the $GL_h$-invariants of $QP_n^{\otimes h}(\omega)$ in every positive degree $n$ is a non-trivial undertaking. Nonetheless, this computation becomes significantly more tractable when a monomial basis of $QP_n^{\otimes h}(\omega)$ is precisely determined.)

\subsection{Proof of Theorem \ref{dlc0}}\label{s2.0}

Undoubtedly if $h > n,$ then $QP^{\otimes h}_{n} \cong (QP^{\otimes h}_{n})^{0}.$ So, the coinvariant $(\mathbb F_2\otimes_{GL_h}{\rm Ann}_{\overline{\mathcal A}}[P^{\otimes h}]^{*})_{n}$ vanishes for any $1\leq n\leq n_0$ and $h\geq n+1.$ Let us consider the following weight vectors:
$$ {\omega}^{*}_{(1)}:=(3,1,1), \ \  {\omega}^{*}_{(2)}:=(3,3), \ \  {\omega}^{*}_{(3)}:=(5,2), \ \  {\omega}^{*}_{(4)}:=(7,1),\ \ {\omega}^{*}_{(5)}:=(9,0).$$
Consequently $\deg({\omega}^{*}_{(1)}) = \deg({\omega}^{*}_{(2)}) = \deg({\omega}^{*}_{(3)}) = \deg({\omega}^{*}_{(4)}) = \deg({\omega}^{*}_{(5)}) =  9.$ It can be seen that $(QP^{\otimes 9}_{9})^{>0}\cong (\widetilde {Sq^0_*})_{9}(QP^{\otimes 9}_{9}) \cong  \mathbb F_2$, and so, $(QP^{\otimes 9}_{9})^{>0}= \mathbb F_2\big[\prod_{1\leq i\leq 9}t_i\big]_{\omega^{*}_{(5)}}.$ In combination with the earlier studies by Peterson \cite{Peterson}, Kameko \cite{Kameko}, Sum \cite{Sum2}, Sum and T'in \cite{Sum-Tin}, and Mothebe et al. \cite{MKR}, the following isomorphisms are obtained:
$$ (QP^{\otimes h}_{9})^{>0}\cong \left\{\begin{array}{ll}
(QP^{\otimes h}_{9})^{>0}({\omega}^{*}_{(1)})\bigoplus (QP^{\otimes h}_{9})^{>0}({\omega}^{*}_{(2)})&\mbox{if $3\leq h\leq 4$},\\[2mm]
(QP^{\otimes h}_{9})^{>0}({\omega}^{*}_{(1)})\bigoplus (QP^{\otimes h}_{9})^{>0}({\omega}^{*}_{(2)})\bigoplus (QP^{\otimes h}_{9})^{>0}({\omega}^{*}_{(3)})&\mbox{if $h = 5$},\\[2mm]
(QP^{\otimes h}_{9})^{>0}({\omega}^{*}_{(2)})\bigoplus (QP^{\otimes h}_{9})^{>0}({\omega}^{*}_{(3)})&\mbox{if $h = 6$},\\[2mm]
(QP^{\otimes h}_{9})^{>0}({\omega}^{*}_{(3)})\bigoplus (QP^{\otimes h}_{9})^{>0}({\omega}^{*}_{(4)})&\mbox{if $h = 7$},\\[2mm]
(QP^{\otimes h}_{9})^{>0}({\omega}^{*}_{(4)})&\mbox{if $h = 8$},\\[2mm]
(QP^{\otimes h}_{9})^{>0}({\omega}^{*}_{(5)})&\mbox{if $h = 9$},\\[2mm]
O&\mbox{if $h \geq n_0.$}
\end{array}\right.$$
Hence the dimensions of the indecomposables $(QP^{\otimes h}_{9})^{>0}({\omega}^{*}_{(j)})$ are determined as follows:

\centerline{\begin{tabular}{c||cccccccccccccccc}
$j$  &&$1$ && $2$ && $3$ && $4$ && $5$ \cr
\hline
\hline
\ $\dim (QP^{\otimes 3}_{9})^{>0}(\omega^{*}_{(j)})$ && $6$ && $1$ && $0$ &&$0$ &&$0$    \cr
\hline
\hline
\ $\dim (QP^{\otimes 4}_{9})^{>0}({\omega}^{*}_{(j)})$ && $12$ && $6$ && $0$ &&$0$  &&$0$  \cr
\hline
\hline
\ $\dim (QP^{\otimes 5}_{9})^{>0}({\omega}^{*}_{(j)})$ && $6$ && $15$ && $10$ &&$0$  &&$0$ \cr
\hline
\hline
\ $\dim (QP^{\otimes 6}_{9})^{>0}({\omega}^{*}_{(j)})$ && $0$ && $10$ && $24$ &&$0$ &&$0$   \cr
\hline
\hline
\ $\dim (QP^{\otimes 7}_{9})^{>0}({\omega}^{*}_{(j)})$ && $0$ && $0$ && $14$ &&$7$ &&$0$  \cr
\hline
\hline
\ $\dim (QP^{\otimes 8}_{9})^{>0}({\omega}^{*}_{(j)})$ && $0$ && $0$ && $0$ &&$7$ &&$0$  \cr
\hline
\hline
\ $\dim (QP^{\otimes 9}_{9})^{>0}({\omega}^{*}_{(j)})$ && $0$ && $0$ && $0$ &&$0$ &&$1$  \cr
\end{tabular}}

Through a straightforward calculation utilizing the aforementioned data and Remark \ref{nxp10}, we obtain
$$ \dim (QP^{\otimes h}_{9})^{0}(\omega^{*}_{(j)})= \left\{\begin{array}{ll}
6\binom{h}{3} &\mbox{if $j = 1$ and $h=4$},\\[2mm]
6\binom{h}{3} + 12\binom{h}{4}&\mbox{if $j = 1$ and $h=5$},\\[2mm]
6\binom{h}{3} + 12\binom{h}{4} + 6\binom{h}{5}&\mbox{if $j = 1$ and $h\geq 6$},\\[2mm]
\binom{h}{3} &\mbox{if $j = 2$ and $h=4$},\\[2mm]
\binom{h}{3} + 6\binom{h}{4}&\mbox{if $j = 2$ and $h=5$},\\[2mm]
\binom{h}{3} + 6\binom{h}{4} + 15\binom{h}{5}&\mbox{if $j = 2$ and $h=6$},\\[2mm]
\binom{h}{3} + 6\binom{h}{4} + 15\binom{h}{5} + 10\binom{h}{6}&\mbox{if $j = 2$ and $h\geq 7$},\\[2mm]
10\binom{h}{5}  &\mbox{if $j = 3$ and $h=6$},\\[2mm]
10\binom{h}{5} + 24\binom{h}{6}&\mbox{if $j = 3$ and $h=7$},\\[2mm]
10\binom{h}{5} + 24\binom{h}{6} + 14\binom{h}{7} &\mbox{if $j = 3$ and $h\geq 8$},\\[2mm]
7\binom{h}{7} &\mbox{if $j = 4$ and $h=8$},\\[2mm]
7\bigg(\binom{h}{7} + \binom{h}{8}\bigg)   &\mbox{if $j = 4$ and $h\geq 9.$}
\end{array}\right.$$
Then, for each $h\geq n_0$ we have an isomorphism $QP^{\otimes h}_{9} \cong \bigoplus_{1\leq j\leq 5}(QP^{\otimes h}_{9})^{0}(\omega^{*}_{(j)}).$ 

\medskip

$\bullet$ For $h = 9$ and $n = 9,$ since $(\widetilde {Sq^0_*})_{9}: QP^{\otimes 9}_{9}  \longrightarrow \mathbb F_2$ is an epimorphism, $$QP^{\otimes 9}_{9} \cong \mathbb F_2\bigoplus {\rm Ker}((\widetilde {Sq^0_*})_{9}) \cong \mathbb F_2\big[\prod_{1\leq i\leq 9}t_i\big]_{\omega^{*}_{(5)}}\bigoplus \big(\bigoplus_{1\leq j\leq 4}(QP^{\otimes 9}_{9})^{0}(\omega^{*}_{(j)})\big).$$ This shows that $(QP^{\otimes 9}_{9})^{0}\cong \bigoplus_{1\leq j\leq 4}(QP^{\otimes 9}_{9})^{0}(\omega^{*}_{(j)}).$  Hence, $\dim (QP^{\otimes 9}_{9})^{0} = \sum_{1\leq j\leq 4}\dim (QP^{\otimes 9}_{9})^{0}(\omega^{*}_{(j)}) = 10437$ and $\dim QP^{\otimes 9}_{9} = 10438.$ Using this result and the homomorphisms $\sigma_d: P^{\otimes 9}\longrightarrow P^{\otimes 9},\, 1\leq d\leq 9,$ we claim that $[QP^{\otimes 9}_{9}]^{GL_9}$  is zero, and so is, $(\mathbb F_2\otimes_{GL_9}{\rm Ann}_{\overline{\mathcal A}}[P^{\otimes 9}]^{*})_{9}.$ 

\medskip

$\bullet$ For $9\leq h\leq n_0$ and $n = n_0,$ by a simple computation using Theorem \ref{dlMKR} and Corollaries \ref{hq10-1-0}, \ref{hq10-1-1}, one has the following isomorphisms:
$$ \begin{array}{ll}
\medskip
QP^{\otimes h}_{n_0}&\cong \bigoplus_{1\leq j\leq 6}(QP^{\otimes h}_{n_0})^{0}(\overline{\omega}^{j})\bigoplus (QP^{\otimes h}_{n_0})^{>0}(\overline{\omega}^{6}),\ \mbox{for $h = 9$},\\ 
QP^{\otimes h}_{n_0}&\cong \bigoplus_{1\leq j\leq 6}(QP^{\otimes h}_{n_0})^{0}(\overline{\omega}^{j})\bigoplus (QP^{\otimes h}_{n_0})^{>0}(\overline{\omega}^{7}),\ \mbox{for $h = n_0,$} 
\end{array}$$
where $\overline{\omega}^{6} := (8,1)$ and $\overline{\omega}^{7} := (n_0,0).$ It is to be noted that $\bigoplus_{1\leq j\leq 6}(QP^{\otimes n_0}_{n_0})^{0}(\overline{\omega}^{j}) \cong {\rm Ker}((\widetilde {Sq^0_*})_{n_0}),$ where ${\rm Ker}((\widetilde {Sq^0_*})_{n_0})$  is the kernel of the Kameko homomorphism $(\widetilde {Sq^0_*})_{n_0}: QP^{\otimes n_0}_{n_0}\longrightarrow \mathbb F_2.$ The dimensions of the cohit spaces $(QP^{\otimes h}_{n_0})^{0}(\overline{\omega}^{j}),\, 1\leq j\leq 5$ are explicitly determined as in Corollary \ref{hq10-1-1}. Based on Theorem \ref{dlMKR} and direct calculations, we find that 
$$ \dim (QP^{\otimes h}_{n_0})^{0}(\overline{\omega}^{6}) =  \left\{\begin{array}{ll}
72&\mbox{if $h = 9$},\\[1mm]
8\bigg(\binom{n_0}{8} + \binom{n_0}{9}\bigg) = 125 &\mbox{if $h = n_0$},
\end{array}\right.$$
$$ \dim (QP^{\otimes h}_{n_0})^{>0}(\overline{\omega}^{j}) =  \left\{\begin{array}{ll}
8&\mbox{if $h = 9$ and $j = 6$},\\[1mm]
\dim \mathbb F_2 = 1 &\mbox{if $h = n_0$ and $j = 7$}.
\end{array}\right.$$
Using these data and the homomorphisms $\sigma_d,\, 1\leq d\leq n_0,$ we state that $[QP^{n_0}_{n_0}]^{GL_{n_0}} = 0$ and that for each $1\leq j\leq 6,$ the invariant $[(QP^{\otimes h}_{n_0})^{>0}(\overline{\omega}^{j})]^{GL_9}$ is zero, and so is, $(\mathbb F_2\otimes_{GL_9}{\rm Ann}_{\overline{\mathcal A}}[P^{\otimes 9}]^{*})_{n_0}.$ We will describe explicitly $[(QP^{\otimes 9}_{n_0})^{>0}(\overline{\omega}^{j})]^{GL_9}$ for $j = 6.$ The others can be obtained by similar calculations. As shown above, $(QP^{\otimes 9}_{n_0})^{>0}(\overline{\omega}^{6})$ is an $\mathbb F_2$-vector space of dimension $8$ with a monomial basis represented by the following admissible monomials:

\begin{center}
\begin{tabular}{llr}
${\rm a}_{73}=t_1t_2t_3t_4t_5t_6t_7t_8t_9^{2}$, & ${\rm a}_{74}=t_1t_2t_3t_4t_5t_6t_7t_8^{2}t_9$, & \multicolumn{1}{l}{${\rm a}_{75}=t_1t_2t_3t_4t_5t_6t_7^{2}t_8t_9$,} \\
${\rm a}_{76}=t_1t_2t_3t_4t_5t_6^{2}t_7t_8t_9$, & ${\rm a}_{77}=t_1t_2t_3t_4t_5^{2}t_6t_7t_8t_9$, & \multicolumn{1}{l}{${\rm a}_{78}=t_1t_2t_3t_4^{2}t_5t_6t_7t_8t_9$,} \\
${\rm a}_{79}=t_1t_2t_3^{2}t_4t_5t_6t_7t_8t_9$, & ${\rm a}_{80}=t_1t_2^{2}t_3t_4t_5t_6t_7t_8t_9.$ &  
\end{tabular}%
\end{center}

Suppose $[f]_{\overline{\omega}^{6}}\in [(QP^{\otimes 9}_{n_0})^{>0}(\overline{\omega}^{6})]^{\Sigma_9}.$ Then, we have $f\equiv_{\overline{\omega}^{6}} \sum_{73\leq i\leq 80} \gamma_{i}{\rm a}_{i}$ where $\gamma_i\in \mathbb F_2$ for every $i.$ Let us consider the homomorphisms $\sigma_d: P^{\otimes 9}\longrightarrow P^{\otimes 9},\, 1\leq d\leq 8.$ An easy calculation shows:
$$ \begin{array}{ll}
 \sigma_1(f) &\equiv_{\overline{\omega}^{6}}\sum_{73\leq j\leq 79} \gamma_{j}{\rm a}_{j} + \gamma_{80}t_1^{2}t_2t_3t_4t_5t_6^{2}t_7t_8t_9\\
&\equiv_{\overline{\omega}^{6}} \sum_{73\leq j\leq 79}(\gamma_{j} + \gamma_{80}){\rm a}_{j} + \gamma_{80}{\rm a}_{80},\\
&(\mbox{since $t_1^{2}t_2t_3t_4t_5t_6^{2}t_7t_8t_9 = Sq^{1}(t_1t_2t_3t_4t_5t_6^{2}t_7t_8t_9) + \sum_{73\leq j\leq 80}{\rm a}_{j}$}),\\
\medskip
\sigma_2(f) &\equiv_{\overline{\omega}^{6}}\sum_{73\leq j\leq 78} \gamma_{j}{\rm a}_{j} + \gamma_{79}{\rm a}_{80} + \gamma_{80}{\rm a}_{79},\ \ \sigma_3(f) \equiv_{\overline{\omega}^{6}}\sum_{j\neq 78,\, 79} \gamma_{j}{\rm a}_{j} + \gamma_{78}{\rm a}_{79} + \gamma_{79}{\rm a}_{78} ,\\
\medskip
\sigma_4(f) &\equiv_{\overline{\omega}^{6}}\sum_{j\neq 77,\, 78} \gamma_{j}{\rm a}_{j} + \gamma_{77}{\rm a}_{78} + \gamma_{78}{\rm a}_{77},\ \ \sigma_5(f) \equiv_{\overline{\omega}^{6}}\sum_{j\neq 76,\, 77} \gamma_{j}{\rm a}_{j} + \gamma_{76}{\rm a}_{77} + \gamma_{77}{\rm a}_{76},\\
\medskip
\sigma_6(f) &\equiv_{\overline{\omega}^{6}}\sum_{j\neq 75,\, 76} \gamma_{j}{\rm a}_{j} + \gamma_{75}{\rm a}_{76} + \gamma_{76}{\rm a}_{75},\ \ \sigma_7(f) \equiv_{\overline{\omega}^{6}}\sum_{j\neq 74,\, 75} \gamma_{j}{\rm a}_{j} + \gamma_{74}{\rm a}_{75} + \gamma_{75}{\rm a}_{74},\\
\sigma_8(f) &\equiv_{\overline{\omega}^{6}}\sum_{j\neq 73,\, 74} \gamma_{j}{\rm a}_{j} + \gamma_{73}{\rm a}_{74} + \gamma_{74}{\rm a}_{73}.
\end{array}$$
By these equalities and the relations $\sigma_d(f) + f\equiv_{\overline{\omega}^{6}} 0,\, 1\leq d\leq 8$ we get $\gamma_{i} = 0,\, 73\leq i\leq 80.$ Thus, $[QP^{\otimes 9}_{n_0}(\overline{\omega}^{6})]^{GL_9} = [(QP^{\otimes 9}_{n_0})^{0}(\overline{\omega}^{6})]^{GL_9}.$ Note that $(QP^{\otimes 9}_{n_0})^{0}(\overline{\omega}^{6})$ is an $\mathbb F_2$-vector space of dimension $72$ with a monomial basis represented by the following admissible monomials:

\begin{center}
\begin{tabular}{llll}
${\rm a}_{1}=t_2t_3t_4t_5t_6t_7t_8t_9^{3}$, & ${\rm a}_{2}=t_2t_3t_4t_5t_6t_7t_8^{3}t_9$, & ${\rm a}_{3}=t_2t_3t_4t_5t_6t_7^{3}t_8t_9$, & ${\rm a}_{4}=t_2t_3t_4t_5t_6^{3}t_7t_8t_9$, \\
${\rm a}_{5}=t_2t_3t_4t_5^{3}t_6t_7t_8t_9$, & ${\rm a}_{6}=t_2t_3t_4^{3}t_5t_6t_7t_8t_9$, & ${\rm a}_{7}=t_2t_3^{3}t_4t_5t_6t_7t_8t_9$, & ${\rm a}_{8}=t_2^{3}t_3t_4t_5t_6t_7t_8t_9$, \\
${\rm a}_{9}=t_1t_3t_4t_5t_6t_7t_8t_9^{3}$, & ${\rm a}_{10}=t_1t_3t_4t_5t_6t_7t_8^{3}t_9$, & ${\rm a}_{11}=t_1t_3t_4t_5t_6t_7^{3}t_8t_9$, & ${\rm a}_{12}=t_1t_3t_4t_5t_6^{3}t_7t_8t_9$, \\
${\rm a}_{13}=t_1t_3t_4t_5^{3}t_6t_7t_8t_9$, & ${\rm a}_{14}=t_1t_3t_4^{3}t_5t_6t_7t_8t_9$, & ${\rm a}_{15}=t_1t_3^{3}t_4t_5t_6t_7t_8t_9$, & ${\rm a}_{16}=t_1^{3}t_3t_4t_5t_6t_7t_8t_9$, \\
${\rm a}_{17}=t_1t_2t_4t_5t_6t_7t_8t_9^{3}$, & ${\rm a}_{18}=t_1t_2t_4t_5t_6t_7t_8^{3}t_9$, & ${\rm a}_{19}=t_1t_2t_4t_5t_6t_7^{3}t_8t_9$, & ${\rm a}_{20}=t_1t_2t_4t_5t_6^{3}t_7t_8t_9$, \\
${\rm a}_{21}=t_1t_2t_4t_5^{3}t_6t_7t_8t_9$, & ${\rm a}_{22}=t_1t_2t_4^{3}t_5t_6t_7t_8t_9$, & ${\rm a}_{23}=t_1t_2^{3}t_4t_5t_6t_7t_8t_9$, & ${\rm a}_{24}=t_1^{3}t_2t_4t_5t_6t_7t_8t_9$, \\
${\rm a}_{25}=t_1t_2t_3t_5t_6t_7t_8t_9^{3}$, & ${\rm a}_{26}=t_1t_2t_3t_5t_6t_7t_8^{3}t_9$, & ${\rm a}_{27}=t_1t_2t_3t_5t_6t_7^{3}t_8t_9$, & ${\rm a}_{28}=t_1t_2t_3t_5t_6^{3}t_7t_8t_9$, \\
${\rm a}_{29}=t_1t_2t_3t_5^{3}t_6t_7t_8t_9$, & ${\rm a}_{30}=t_1t_2t_3^{3}t_5t_6t_7t_8t_9$, & ${\rm a}_{31}=t_1t_2^{3}t_3t_5t_6t_7t_8t_9$, & ${\rm a}_{32}=t_1^{3}t_2t_3t_5t_6t_7t_8t_9$, \\
${\rm a}_{33}=t_1t_2t_3t_4t_6t_7t_8t_9^{3}$, & ${\rm a}_{34}=t_1t_2t_3t_4t_6t_7t_8^{3}t_9$, & ${\rm a}_{35}=t_1t_2t_3t_4t_6t_7^{3}t_8t_9$, & ${\rm a}_{36}=t_1t_2t_3t_4t_6^{3}t_7t_8t_9$, \\
${\rm a}_{37}=t_1t_2t_3t_4^{3}t_6t_7t_8t_9$, & ${\rm a}_{38}=t_1t_2t_3^{3}t_4t_6t_7t_8t_9$, & ${\rm a}_{39}=t_1t_2^{3}t_3t_4t_6t_7t_8t_9$, & ${\rm a}_{40}=t_1^{3}t_2t_3t_4t_6t_7t_8t_9$, \\
${\rm a}_{41}=t_1t_2t_3t_4t_5t_7t_8t_9^{3}$, & ${\rm a}_{42}=t_1t_2t_3t_4t_5t_7t_8^{3}t_9$, & ${\rm a}_{43}=t_1t_2t_3t_4t_5t_7^{3}t_8t_9$, & ${\rm a}_{44}=t_1t_2t_3t_4t_5^{3}t_7t_8t_9$, \\
${\rm a}_{45}=t_1t_2t_3t_4^{3}t_5t_7t_8t_9$, & ${\rm a}_{46}=t_1t_2t_3^{3}t_4t_5t_7t_8t_9$, & ${\rm a}_{47}=t_1t_2^{3}t_3t_4t_5t_7t_8t_9$, & ${\rm a}_{48}=t_1^{3}t_2t_3t_4t_5t_7t_8t_9$, \\
${\rm a}_{49}=t_1t_2t_3t_4t_5t_6t_7t_9^{3}$, & ${\rm a}_{50}=t_1t_2t_3t_4t_5t_6t_7^{3}t_9$, & ${\rm a}_{51}=t_1t_2t_3t_4t_5t_6^{3}t_7t_9$, & ${\rm a}_{52}=t_1t_2t_3t_4t_5^{3}t_6t_7t_9$, \\
${\rm a}_{53}=t_1t_2t_3t_4^{3}t_5t_6t_7t_9$, & ${\rm a}_{54}=t_1t_2t_3^{3}t_4t_5t_6t_7t_9$, & ${\rm a}_{55}=t_1t_2^{3}t_3t_4t_5t_6t_7t_9$, & ${\rm a}_{56}=t_1^{3}t_2t_3t_4t_5t_6t_7t_9$, \\
${\rm a}_{57}=t_1t_2t_3t_4t_5t_6t_8t_9^{3}$, & ${\rm a}_{58}=t_1t_2t_3t_4t_5t_6t_8^{3}t_9$, & ${\rm a}_{59}=t_1t_2t_3t_4t_5t_6^{3}t_8t_9$, & ${\rm a}_{60}=t_1t_2t_3t_4t_5^{3}t_6t_8t_9$, \\
${\rm a}_{61}=t_1t_2t_3t_4^{3}t_5t_6t_8t_9$, & ${\rm a}_{62}=t_1t_2t_3^{3}t_4t_5t_6t_8t_9$, & ${\rm a}_{63}=t_1t_2^{3}t_3t_4t_5t_6t_8t_9$, & ${\rm a}_{64}=t_1^{3}t_2t_3t_4t_5t_6t_8t_9$, \\
${\rm a}_{65}=t_1t_2t_3t_4t_5t_6t_7t_8^{3}$, & ${\rm a}_{66}=t_1t_2t_3t_4t_5t_6t_7^{3}t_8$, & ${\rm a}_{67}=t_1t_2t_3t_4t_5t_6^{3}t_7t_8$, & ${\rm a}_{68}=t_1t_2t_3t_4t_5^{3}t_6t_7t_8$, \\
${\rm a}_{69}=t_1t_2t_3t_4^{3}t_5t_6t_7t_8$, & ${\rm a}_{70}=t_1t_2t_3^{3}t_4t_5t_6t_7t_8$, & ${\rm a}_{71}=t_1t_2^{3}t_3t_4t_5t_6t_7t_8$, & ${\rm a}_{72}=t_1^{3}t_2t_3t_4t_5t_6t_7t_8.$
\end{tabular}%
\end{center}

\medskip

Suppose $[g]_{\overline{\omega}^{6}}\in [QP^{\otimes 9}_{n_0}(\overline{\omega}^{6}]^{\Sigma_9}.$ Then, one has $g\equiv_{\overline{\omega}^{6}} \sum_{1\leq j\leq 72} \beta_{i}{\rm a}_{i},$ in which $\beta_i\in \mathbb F_2,\, 1\leq i\leq 72.$ Using the homomorphisms $\sigma_d: P^{\otimes 9}\longrightarrow P^{\otimes 9},$ for $1\leq d\leq 8,$ we obtain the following equalities:
$$ \begin{array}{ll}
\sigma_1(g)&\equiv_{\overline{\omega}^{6}} \sum_{1\leq i\leq 8}\beta_{i}{\rm a}_{i+8} + \sum_{9\leq i\leq 16}\beta_{i}{\rm a}_{i-8} + \sum_{17\leq i\leq 22}\beta_{i}{\rm a}_{i} + \beta_{23}{\rm a}_{24} + \beta_{24}{\rm a}_{23}\\
&\quad + \sum_{25\leq i\leq 30}\beta_{i}{\rm a}_{i} +  \beta_{31}{\rm a}_{32} + \beta_{32}{\rm a}_{31} + \sum_{33\leq i\leq 38}\beta_{i}{\rm a}_{i} +  \beta_{39}{\rm a}_{40} + \beta_{40}{\rm a}_{39} \\
&\quad+ \sum_{41\leq i\leq 46}\beta_{i}{\rm a}_{i} +  \beta_{47}{\rm a}_{48} + \beta_{48}{\rm a}_{47} + \sum_{49\leq i\leq 54}\beta_{i}{\rm a}_{i} +  \beta_{55}{\rm a}_{56} + \beta_{56}{\rm a}_{55}\\
\medskip
&\quad + \sum_{57\leq i\leq 62}\beta_{i}{\rm a}_{i} +  \beta_{63}{\rm a}_{64} + \beta_{64}{\rm a}_{63}+\sum_{65\leq i\leq 70}\beta_{i}{\rm a}_{i} +  \beta_{71}{\rm a}_{72} + \beta_{72}{\rm a}_{71},\\
\sigma_2(g)&\equiv_{\overline{\omega}^{6}} \sum_{1\leq i\leq 6}\beta_{i}{\rm a}_{i} +\beta_{7}{\rm a}_{8} + \beta_{8}{\rm a}_{7} + 
 \sum_{9\leq i\leq 16}\beta_{i}{\rm a}_{i+8} + \sum_{17\leq i\leq 24}\beta_{i}{\rm a}_{i-8}\\
&\quad + \sum_{25\leq i\leq 29}\beta_{i}{\rm a}_{i} +  \beta_{30}{\rm a}_{31} + \beta_{31}{\rm a}_{30} + \sum_{32\leq i\leq 37}\beta_{i}{\rm a}_{i} +  \beta_{38}{\rm a}_{39} + \beta_{39}{\rm a}_{38}\\
&\quad + \sum_{40\leq i\leq 45}\beta_{i}{\rm a}_{i} +  \beta_{46}{\rm a}_{47} + \beta_{47}{\rm a}_{46} + \sum_{48\leq i\leq 53}\beta_{i}{\rm a}_{i} +  \beta_{54}{\rm a}_{55} + \beta_{55}{\rm a}_{54}\\
\medskip
&\quad+ \sum_{56\leq i\leq 61}\beta_{i}{\rm a}_{i} +  \beta_{62}{\rm a}_{63} + \beta_{63}{\rm a}_{62}+\sum_{64\leq i\leq 69}\beta_{i}{\rm a}_{i} +  \beta_{70}{\rm a}_{71} + \beta_{71}{\rm a}_{70} + \beta_{72}{\rm a}_{72},\\
\sigma_3(g)&\equiv_{\overline{\omega}^{6}} \sum_{1\leq i\leq 5}\beta_{i}{\rm a}_{i} +\beta_{6}{\rm a}_{7} + \beta_{7}{\rm a}_{6} + 
 \sum_{8\leq i\leq 13}\beta_{i}{\rm a}_{i} +\beta_{14}{\rm a}_{15} + \beta_{15}{\rm a}_{14} + \beta_{16}{\rm a}_{16} \\
&\quad+ \sum_{17\leq i\leq 24}\beta_{i}{\rm a}_{i+8} + \sum_{25\leq i\leq 32}\beta_{i}{\rm a}_{i-8}  + \sum_{33\leq i\leq 36}\beta_{i}{\rm a}_{i} +  \beta_{37}{\rm a}_{38} + \beta_{38}{\rm a}_{37}\\
&\quad + \sum_{39\leq i\leq 44}\beta_{i}{\rm a}_{i} +  \beta_{45}{\rm a}_{46} + \beta_{46}{\rm a}_{45} + \sum_{47\leq i\leq 52}\beta_{i}{\rm a}_{i} +  \beta_{53}{\rm a}_{54} + \beta_{54}{\rm a}_{53}\\
\medskip
&\quad+ \sum_{55\leq i\leq 60}\beta_{i}{\rm a}_{i} +  \beta_{61}{\rm a}_{62} + \beta_{62}{\rm a}_{61}+\sum_{63\leq i\leq 68}\beta_{i}{\rm a}_{i} +  \beta_{69}{\rm a}_{70} + \beta_{70}{\rm a}_{69} + \sum_{71\leq i\leq 72}\beta_{i}{\rm a}_{i},\\
\end{array}$$
\newpage
$$ \begin{array}{ll}
\sigma_4(g)&\equiv_{\overline{\omega}^{6}} \sum_{1\leq i\leq 4}\beta_{i}{\rm a}_{i} +\beta_{5}{\rm a}_{6} + \beta_{6}{\rm a}_{5} + 
 \sum_{7\leq i\leq 12}\beta_{i}{\rm a}_{i} +\beta_{13}{\rm a}_{14} + \beta_{14}{\rm a}_{13}\\
&\quad + \sum_{15\leq i\leq 20}\beta_{i}{\rm a}_{i} +\beta_{21}{\rm a}_{22} + \beta_{22}{\rm a}_{21}  + \sum_{23\leq i\leq 24}\beta_{i}{\rm a}_{i}  + \sum_{25\leq i\leq 32}\beta_{i}{\rm a}_{i+8} \\
&\quad+ \sum_{33\leq i\leq 40}\beta_{i}{\rm a}_{i-8}  + \sum_{41\leq i\leq 43}\beta_{i}{\rm a}_{i} +  \beta_{44}{\rm a}_{45} + \beta_{45}{\rm a}_{44}+ \sum_{46\leq i\leq 51}\beta_{i}{\rm a}_{i} +  \beta_{52}{\rm a}_{53} + \beta_{53}{\rm a}_{52}\\
\medskip
&\quad+\sum_{54\leq i\leq 59}\beta_{i}{\rm a}_{i} +  \beta_{60}{\rm a}_{61} + \beta_{61}{\rm a}_{60} + \sum_{62\leq i\leq 67}\beta_{i}{\rm a}_{i}+  \beta_{68}{\rm a}_{69} + \beta_{69}{\rm a}_{68} +  \sum_{70\leq i\leq 72}\beta_{i}{\rm a}_{i},\\
\sigma_5(g)&\equiv_{\overline{\omega}^{6}} \sum_{1\leq i\leq 3}\beta_{i}{\rm a}_{i} +\beta_{4}{\rm a}_{5} + \beta_{5}{\rm a}_{4} + 
 \sum_{6\leq i\leq 11}\beta_{i}{\rm a}_{i} +\beta_{12}{\rm a}_{13} + \beta_{13}{\rm a}_{12}\\
&\quad + \sum_{14\leq i\leq 19}\beta_{i}{\rm a}_{i} +\beta_{20}{\rm a}_{21} + \beta_{21}{\rm a}_{20}  + \sum_{22\leq i\leq 27}\beta_{i}{\rm a}_{i} +\beta_{28}{\rm a}_{29} + \beta_{29}{\rm a}_{28}  + \sum_{30\leq i\leq 32}\beta_{i}{\rm a}_{i} \\
&\quad+ \sum_{33\leq i\leq 40}\beta_{i}{\rm a}_{i+8}  + \sum_{41\leq i\leq 48}\beta_{i}{\rm a}_{i-8} +  \beta_{49}{\rm a}_{50} + \beta_{50}{\rm a}_{49}+  \beta_{51}{\rm a}_{52} + \beta_{52}{\rm a}_{51}\\
\medskip
&\quad + \sum_{53\leq i\leq 58}\beta_{i}{\rm a}_{i} +  \beta_{59}{\rm a}_{60} + \beta_{60}{\rm a}_{59}+\sum_{61\leq i\leq 66}\beta_{i}{\rm a}_{i} +  \beta_{67}{\rm a}_{68} + \beta_{68}{\rm a}_{67} + \sum_{69\leq i\leq 72}\beta_{i}{\rm a}_{i},\\
\sigma_6(g)&\equiv_{\overline{\omega}^{6}} \sum_{1\leq i\leq 2}\beta_{i}{\rm a}_{i} +\beta_{3}{\rm a}_{4} + \beta_{4}{\rm a}_{3} + 
 \sum_{5\leq i\leq 10}\beta_{i}{\rm a}_{i} +\beta_{11}{\rm a}_{12} + \beta_{12}{\rm a}_{11}\\
&\quad + \sum_{13\leq i\leq 18}\beta_{i}{\rm a}_{i} +\beta_{19}{\rm a}_{20} + \beta_{20}{\rm a}_{19}  + \sum_{21\leq i\leq 26}\beta_{i}{\rm a}_{i} +\beta_{27}{\rm a}_{28} + \beta_{28}{\rm a}_{27}  + \sum_{29\leq i\leq 34}\beta_{i}{\rm a}_{i} \\
&\quad +\beta_{35}{\rm a}_{36} + \beta_{36}{\rm a}_{35} + \sum_{37\leq i\leq 40}\beta_{i}{\rm a}_{i} + \sum_{41\leq i\leq 48}\beta_{i}{\rm a}_{i+16} +  \beta_{49}{\rm a}_{49} + \beta_{50}{\rm a}_{51}+  \beta_{51}{\rm a}_{50} \\
\medskip
&\quad  + \sum_{52\leq i\leq 56}\beta_{i}{\rm a}_{i}+ \sum_{57\leq i\leq 64}\beta_{i}{\rm a}_{i-16} +  \beta_{65}{\rm a}_{65} + \beta_{66}{\rm a}_{67} + \beta_{67}{\rm a}_{66} +\sum_{68\leq i\leq 72}\beta_{i}{\rm a}_{i},\\
\sigma_7(g)&\equiv_{\overline{\omega}^{6}} \beta_{1}{\rm a}_{1} + \beta_{2}{\rm a}_{3}+\beta_{3}{\rm a}_{2} +  \sum_{4\leq i\leq 9}\beta_{i}{\rm a}_{i} +\beta_{10}{\rm a}_{11} + \beta_{11}{\rm a}_{10} + \sum_{12\leq i\leq 17}\beta_{i}{\rm a}_{i} +\beta_{18}{\rm a}_{19} \\
&\quad + \beta_{19}{\rm a}_{18} + \sum_{20\leq i\leq 25}\beta_{i}{\rm a}_{i} +\beta_{26}{\rm a}_{27} + \beta_{27}{\rm a}_{26}  + \sum_{28\leq i\leq 33}\beta_{i}{\rm a}_{i} +\beta_{34}{\rm a}_{35} + \beta_{35}{\rm a}_{34} \\
&\quad   + \sum_{36\leq i\leq 41}\beta_{i}{\rm a}_{i} + \beta_{42}{\rm a}_{43}+\beta_{43}{\rm a}_{42} + \sum_{44\leq i\leq 48}\beta_{i}{\rm a}_{i}+ \sum_{49\leq i\leq 56}\beta_{i}{\rm a}_{i+8} + \sum_{57\leq i\leq 64}\beta_{i}{\rm a}_{i-8}  \\
\medskip
&\quad + \beta_{65}{\rm a}_{66} + \beta_{66}{\rm a}_{65} + \sum_{67\leq i\leq 72}\beta_{i}{\rm a}_{i},\\
\sigma_8(g)&\equiv_{\overline{\omega}^{6}} \beta_{1}{\rm a}_{2} + \beta_{2}{\rm a}_{1}+  \sum_{3\leq i\leq 8}\beta_{i}{\rm a}_{i} +\beta_{9}{\rm a}_{10} + \beta_{10}{\rm a}_{9} + \sum_{11\leq i\leq 16}\beta_{i}{\rm a}_{i} +\beta_{17}{\rm a}_{18} + \beta_{18}{\rm a}_{17} + \sum_{19\leq i\leq 24}\beta_{i}{\rm a}_{i} \\
&\quad +\beta_{25}{\rm a}_{26}  + \beta_{26}{\rm a}_{25}  + \sum_{27\leq i\leq 32}\beta_{i}{\rm a}_{i} +\beta_{33}{\rm a}_{34} + \beta_{34}{\rm a}_{33} + \sum_{35\leq i\leq 40}\beta_{i}{\rm a}_{i} + \beta_{41}{\rm a}_{42}+\beta_{42}{\rm a}_{41} \\
&\quad   + \sum_{43\leq i\leq 48}\beta_{i}{\rm a}_{i}+ \sum_{49\leq i\leq 56}\beta_{i}{\rm a}_{i+16} +\beta_{57}{\rm a}_{58}   + \beta_{58}{\rm a}_{57} + \sum_{59\leq i\leq 64}\beta_{i}{\rm a}_{i} + \sum_{65\leq i\leq 72}\beta_{i}{\rm a}_{i-16}.
\end{array}$$
Then, from the relations $\sigma_d(g)\equiv_{\overline{\omega}^{6}}  g,\, 1\leq d\leq 8,$ one gets $\beta_i = \beta_1$ for all $i,\, 2\leq i\leq 72.$ This means that $[QP^{\otimes 9}_{n_0}(\overline{\omega}^{6}]^{\Sigma_9} = \mathbb F_2\big[\sum_{1\leq i\leq 72}{\rm a}_{i}\big]_{\overline{\omega}^{6}}.$ Then, given any $[h]_{\overline{\omega}^{6}}\in [QP^{\otimes 9}_{n_0}(\overline{\omega}^{6}]^{GL_9},$ we have $h\equiv_{\overline{\omega}^{6}} \zeta\sum_{1\leq j\leq 72}{\rm a}_{i}$ with $\zeta\in \mathbb F_2.$ Since $\sigma_9(h) + h\equiv_{\overline{\omega}^{6}}  0,$ $\zeta = 0$, and therefore, $[QP^{\otimes 9}_{n_0}(\overline{\omega}^{6}]^{GL_9} = 0.$ Now, the theorem can be derived from the above results, in conjunction with the following facts: firstly, the transfer $\{Tr_h^{\mathcal A}\}_{h\geq 0}:  \{\mathbb F_2\otimes_{GL_h}{\rm Ann}_{\overline{\mathcal A}}[P^{\otimes h}]^{*}\}_{h\geq 0}\longrightarrow \{{\rm Ext}_{\mathcal A}^{h, *}(\mathbb F_2, \mathbb F_2)\}_{h\geq 0}$ is an algebra homomorphism (see Singer \cite{Singer}), and secondly, according to \cite{Tangora, Bruner, Lin, Chen, Lin2}, the cohomology groups ${\rm Ext}_{\mathcal A}^{h, h+n}(\mathbb F_2, \mathbb F_2),$ where $h\geq 1$ and $1\leq n\leq n_0$, can be identified as follows:
\newpage
$$ {\rm Ext}_{\mathcal A}^{h, h+n}(\mathbb F_2, \mathbb F_2)=\left\{\begin{array}{ll}
\mathbb F_2 h_1 &\mbox{if $h = 1$ and $n = 1,$}\\[1mm]
0 &\mbox{if $h \geq 2$ and $n = 1,$}\\[1mm]
\mathbb F_2 h^{2}_1 &\mbox{if $h = 2$ and $n = 2,$}\\[1mm]
0 &\mbox{if $h \geq 1,\, h\neq 2$ and $n = 2,$}\\[1mm]
\mathbb F_2 h_2 &\mbox{if $h = 1$ and $n = 3,$}\\[1mm]
\mathbb F_2 h_0h_2 &\mbox{if $h = 2$ and $n = 3,$}\\[1mm]
\mathbb F_2 h_1^{3} &\mbox{if $h = 3$ and $n = 3,$}\\[1mm]
0 &\mbox{if $h \geq 4$ and $n = 3,$}\\[1mm]
0 &\mbox{if $h \geq 1$ and $4\leq n\leq 5,$}\\[1mm]
\mathbb F_2 h^{2}_2 &\mbox{if $h = 2$ and $n = 6,$}\\[1mm]
0 &\mbox{if $h \geq 1,\, h\neq 2$ and $n = 6,$}\\[1mm]
\mathbb F_2 h_3 &\mbox{if $h = 1$ and $n = 7,$}\\[1mm]
\mathbb F_2 h_0h_3 &\mbox{if $h = 2$ and $n = 7,$}\\[1mm]
\mathbb F_2 h_0^{2}h_3 &\mbox{if $h = 3$ and $n = 7,$}\\[1mm]
\mathbb F_2 h_0^{3}h_3 &\mbox{if $h = 4$ and $n = 7,$}\\[1mm]
0 &\mbox{if $h \geq 5$ and $n = 7,$}\\[1mm]
\mathbb F_2 h_1h_3 &\mbox{if $h = 2$ and $n = 8,$}\\[1mm]
\mathbb F_2 c_0 &\mbox{if $h = 3$ and $n = 8,$}\\[1mm]
0 &\mbox{if $h\geq 1,\, h \neq 2,\, 3$ and $n = 8,$}\\[1mm]
\mathbb F_2 h^{3}_2 &\mbox{if $h = 3$ and $n = 9,$}\\[1mm]
\mathbb F_2 h_1c_0 &\mbox{if $h = 4$ and $n = 9,$}\\[1mm]
\mathbb F_2 Ph_1 &\mbox{if $h = 5$ and $n = 9,$}\\[1mm]
0 &\mbox{if $h\geq 1,\, h \neq 3,\, 4,\, 5$ and $n = 9,$}\\[1mm]
\mathbb F_2 h_1Ph_1 &\mbox{if $h = 6$ and $n = n_0,$}\\[1mm]
0&\mbox{if $h\geq 1,\, h \neq 6$ and $n = n_0.$}
\end{array}\right.$$

\subsection{Proof of Theorem \ref{dlc1}}\label{s2.1}
Recall that by Remark \ref{nxp10}, one gets $\dim (QP_{n_1}^{\otimes h})^{0}(\omega) = \sum_{\mu(n_1) = 4\leq k\leq h-1}\binom{h}{k}\dim (QP_{n_1}^{\otimes k})^{> 0}(\omega),$ where $\omega$ is a weight vector of degree $n_1.$ Sine $\mu(n_1) = 4,$ by Theorem \ref{dlWS}(iii), the Kameko homomorphism $(\widetilde {Sq^0_*})_{n_1}: QP^{\otimes 4}_{n_1}\longrightarrow QP^{\otimes 4}_{11}$ is an isomorphism. So, we have the inverse homomorphism $((\widetilde {Sq^0_*})_{n_1})^{-1}: QP^{\otimes 4}_{11}\longrightarrow QP^{\otimes 4}_{n_1},$ determined by $((\widetilde {Sq^0_*})_{n_1})^{-1}([y]) = [t_1t_2t_3t_4y^{2}]$ for all $[y]\in QP^{\otimes 4}_{11}.$ So, a basis of $QP^{\otimes 4}_{n_1} = (QP^{\otimes 4}_{n_1})^{>0}$ is the set of all the equivalence classes represented by the admissible monomials of the form $t_1t_2t_3t_4y^{2}$ where $y\in \mathscr C^{\otimes 4}_{11}.$ By Sum \cite{Sum2}, $|\mathscr C^{\otimes 4}_{11}| = 64$, which means that $\dim (QP^{\otimes 4}_{n_1})^{>0}(\omega) = \dim (QP^{\otimes 4}_{n_1})^{>0} = 64$ if $\omega = (4,3,2,1)$ and $(QP^{\otimes 4}_{n_1})^{>0}(\omega)  = 0$ otherwise. So, by the above formula, $\dim (QP_{n_1}^{\otimes 5})^{0}(\omega)  = \binom{5}{4}\dim (QP_{n_1}^{\otimes 4})^{> 0}(\omega)  = 320$ if $\omega = (4,3,2,1)$ and $(QP^{\otimes 5}_{n_1})^{>0}(\omega)  = 0$ otherwise.  On the other side, since $QP_{n_1}^{\otimes 5}\cong (QP_{n_1}^{\otimes 5})^{0}\bigoplus (QP_{n_1}^{\otimes 5})^{>0},$ by Theorem \ref{dlWW}, we derive $(QP^{\otimes 5}_{n_1})^{>0}(\omega)  = 0$ if $\omega\neq (4,3,2,1)$ and $\dim (QP_{n_1}^{\otimes 5})^{>0}(4,3,21) = \dim (QP_{n_1}^{\otimes 5})^{>0}(4,3,2,1) = 1024 - 320 = 704.$ Therefore,
$$\dim (QP_{n_1}^{\otimes 6})^{0}(4,3,2,1) = \sum_{\mu(n_1) = 4\leq k\leq 5}\binom{6}{k}\dim (QP_{n_1}^{\otimes k})^{> 0}(4,3,2,1)= 5184,
$$
and $(QP^{\otimes 6}_{n_1})^{0}(\omega)  = 0$ if $\omega\neq (4,3,2,1).$ Since $(QP^{\otimes 6}_{n_1})^{0}\cong \bigoplus_{\deg(\omega) = n_1}(QP^{\otimes 6}_{n_1})^{0}(\omega),$ by the above calculations, we obtain $(QP^{\otimes 6}_{n_1})^{0}\cong (QP^{\otimes 6}_{n_1})^{0}(4,3,2,1).$ This completes the proof of Part (i). 

We will now proceed to prove Part (ii) of the theorem. Our first step is to compute the space $U_1$ by explicitly determining the monomial bases of $(QP_{n_1}^{\otimes 6})^{>0}(4,5,1,1)$ and $(QP_{n_1}^{\otimes 6})^{>0}(4,5,3).$ It is important to note that according to Remark \ref{nxp1}(i), if $t\in (P_{n_1}^{\otimes 6})^{>0}(4,5,1,1)$ or $t\in (P_{n_1}^{\otimes 6})^{>0}(4,5,3)$ such that $[t]\in {\rm Ker}((\widetilde {Sq^0_*})_{n_1}),$ then $t$ can be represented as $t_it_jt_kt_l\underline{t}^2,$ where $1\leq i<j<k<l\leq 6$ and $\underline{t}\in \mathscr C^{\otimes 6}_{11}.$ Therefore, in order to compute a monomial basis of $U_1,$ we need to determine all admissible monomials of degree $11$ in $P^{\otimes 6}.$ In \cite{MKR}, Mothebe et al. showed that $QP^{\otimes 6}_{11}$ has dimension $1205.$ Utilizing this result, we can explicitly determine all monomials of the form $t_it_jt_kt_l\underline{t}^2\in P^{\otimes 6}_{n_1}.$ In particular, utilizing the dimension result for $QP^{\otimes h}_{11}$ in \cite{MKR} and conducting a straightforward computation with Remark \ref{nxp10}, Corollary \ref{hq10-2}, and our previous work in \cite{Phuc11}, we obtain the following corollary.
\begin{corl}\label{hqT}
Let us consider the weight vectors of degree $11$:
$$ \begin{array}{ll}
& \widehat{\omega}_{(1)}:= (3,2,1),\ \widehat{\omega}_{(2)}:= (3,4),\ \widehat{\omega}_{(3)}:= (5,1,1),\ \widehat{\omega}_{(4)}:= (5,3),\\
& \widehat{\omega}_{(5)}:= (7,2),\ \widehat{\omega}_{(6)}:= (9,1),\ \widehat{\omega}_{(7)}:= (11,0).
\end{array}$$
Then, for each $h\geq 6,$ we have the isomorphisms:
$$  QP^{\otimes h}_{11}\cong \left\{\begin{array}{ll}
 \bigoplus_{1\leq j\leq 4}QP^{\otimes h}_{11}(\widehat{\omega}_{(j)}) &\mbox{if $h = 6$},\\[1mm]
 \bigoplus_{1\leq j\leq 5}QP^{\otimes h}_{11}(\widehat{\omega}_{(j)}) &\mbox{if $7\leq h\leq 8$},\\[1mm]
 \bigoplus_{1\leq j\leq 6}QP^{\otimes h}_{11}(\widehat{\omega}_{(j)}) &\mbox{if $9\leq h\leq 10$},\\[1mm]
 \bigoplus_{1\leq j\leq 7}QP^{\otimes h}_{11}(\widehat{\omega}_{(j)}) &\mbox{if $h \geq 11,$}
\end{array}\right.$$
where $QP^{\otimes h}_{11}(\widehat{\omega}_{(j)})\cong (QP^{\otimes h}_{11})^{0}(\widehat{\omega}_{(j)})\bigoplus (QP^{\otimes h}_{11})^{>0}(\widehat{\omega}_{(j)}).$ The dimensions of $(QP^{\otimes h}_{11})^{0}(\widehat{\omega}_{(j)})$ and $(QP^{\otimes h}_{11})^{>0}(\widehat{\omega}_{(j)})$ are determined as follows:
$$ 
 \dim (QP^{\otimes h}_{11})^{0}(\widehat{\omega}_{(j)})=\left\{\begin{array}{ll}
8\binom{h}{3} + 32\binom{h}{4} + 40\binom{h}{5} &\mbox{if $h = 6,\, j=1$ {\rm(see \cite{Phuc11})}},\\[1mm]
10\binom{h}{5}  &\mbox{if $h = 6,\, j=2$ {\rm(see \cite{Phuc11})}},\\[1mm]
15\binom{h}{5}  &\mbox{if $h = 6,\, j=3$ {\rm(see \cite{Phuc11})}},\\[1mm]
10\binom{h}{5}  &\mbox{if $h = 6,\, j=4$ {\rm(see \cite{Phuc11})}},\\[1mm]
0  &\mbox{if $h = 6,\, 5\leq j\leq 7$},\\[1mm]
8\binom{h}{3} + 32\binom{h}{4} + 40\binom{h}{5} + 16\binom{h}{6} &\mbox{if $h\geq 7,\, j=1$ {\rm (see \cite{Phuc11})}},\\[1mm]
10\binom{h}{5} +24\binom{h}{6}   &\mbox{if $h = 7,\, j=2$ {\rm(see \cite{Phuc11})}},\\[1mm]
15\binom{h}{5}  + 30\binom{h}{6} &\mbox{if $h = 7,\, j=3$ {\rm(see \cite{Phuc11})}},\\[1mm]
10\binom{h}{5}+45\binom{h}{6}  &\mbox{if $h = 7,\, j=4$ {\rm(see \cite{Phuc11})}},\\[1mm]
0&\mbox{if $h = 7,\, 5\leq j\leq 7$},\\[1mm]
\end{array}\right.$$

\newpage
$$ \dim (QP^{\otimes h}_{11})^{> 0}(\widehat{\omega}_{(j)})=\left\{\begin{array}{ll}
10\binom{h}{5} +24\binom{h}{6} +14\binom{h}{7}   &\mbox{if $h \geq  8,\, j=2$},\\[1mm]
15\binom{h}{5}  + 30\binom{h}{6}  + 15\binom{h}{7} &\mbox{if $h \geq 8,\, j=3$},\\[1mm]
10\binom{h}{5}+45\binom{h}{6} + 70\binom{h}{7}  &\mbox{if $h = 8,\, j=4$},\\[1mm]
21\binom{h}{7}  &\mbox{if $h = 8,\, j=5$},\\[1mm]
0  &\mbox{if $h = 8,\, 6\leq j\leq 7$},\\[1mm]
10\binom{h}{5}+45\binom{h}{6} + 70\binom{h}{7} +35\binom{h}{8}  &\mbox{if $h \geq 9,\, j=4$},\\[1mm]
21\binom{h}{7}+48\binom{h}{8}   &\mbox{if $h \geq 9,\, j=5$},\\[1mm]
9\binom{h}{9}  &\mbox{if $h = 10,\, j=6$},\\[1mm]
9\bigg(\binom{h}{9} + \binom{h}{10}\bigg)  &\mbox{if $h \geq 11,\, j=6$},\\[1mm]
0  &\mbox{if $h  = 11,\, j=7$},\\[1mm]
\binom{h}{11}  &\mbox{if $h \geq 12,\, j=7$},\\[1mm]
16 &\mbox{if $h = 6,\, j = 1$ {\rm(see \cite{Phuc11})}},\\[1mm]
24 &\mbox{if $h = 6,\, j = 2$ {\rm(see \cite{Phuc11})}},\\[1mm]
30 &\mbox{if $h = 6,\, j = 3$ {\rm(see \cite{Phuc11})}},\\[1mm]
45 &\mbox{if $h = 6,\, j = 4$ {\rm(see \cite{Phuc11})}},\\[1mm]
0 &\mbox{if $h = 6,\, 5\leq j\leq 7$},\\[1mm]
0 &\mbox{if $h = 7,\, j = 1,\, 6,\, 7$},\\[1mm]
 \dim (QP^{\otimes h}_{h+4})^{> 0}(\overline{\omega}^{(j-1,\, h)})  &\mbox{if $h = 7,\, 2\leq j\leq 5$ {\rm(see Corollary \ref{hq10-2})}},\\[1mm]
0 &\mbox{if $h = 8,\, j\neq 4,\, 5$},\\[1mm]
35 &\mbox{if $h = 8,\, j = 4$},\\[1mm]
48 &\mbox{if $h = 8,\, j = 5$},\\[1mm]
0 &\mbox{if $h = 9,\, j\neq 5,\, 6$},\\[1mm]
27 &\mbox{if $h = 9,\, j = 5$},\\[1mm]
\dim (QP^{\otimes h}_{1})^{0} = 9 &\mbox{if $h = 9,\, j = 6$},\\[1mm]
9 &\mbox{if $h = 10,\, j = 6$},\\[1mm]
0 &\mbox{if $h \geq 10,\, j \neq 6$}.
\end{array}\right.$$
\end{corl}

It must be noted that for $h = 7,\, 9,\, 11,$ we only need to compute the kernels of the Kameko squaring operations $(\widetilde {Sq^0_*})_{11}: QP^{\otimes h}_{11}\longrightarrow QP^{\otimes h}_{(11-h)/2}.$ Furthermore, we can readily deduce the following:
$$ {\rm Ker}((\widetilde {Sq^0_*})_{11})\cong  \left\{ \begin{array}{ll}
(QP^{\otimes h}_{11})^{0}(\widehat{\omega}_{(1)})\bigoplus \bigg(\bigoplus_{2\leq j\leq 4}QP^{\otimes h}_{11}(\widehat{\omega}_{(j)})\bigg) &\mbox{if $h = 7$},\\[1mm]
\bigg(\bigoplus_{1\leq j\leq 4}(QP^{\otimes h}_{11})^{0}(\widehat{\omega}_{(j)})\bigg)\bigoplus QP^{\otimes h}_{11}(\widehat{\omega}_{(5)})\bigoplus
(QP^{\otimes h}_{11})^{>0}(\widehat{\omega}_{(6)})&\mbox{if $h = 9$},\\[1mm]
\bigoplus_{1\leq j\leq 6}(QP^{\otimes h}_{11})^{0}(\widehat{\omega}_{(j)}) &\mbox{if $h = 11$}.
\end{array}\right.$$
\medskip

The following relevant observation is indispensable in order to establish the theorem: For each positive integer $n,$ the \textit{up Kameko map} $\psi: P_{n}^{\otimes 6}\longrightarrow P_{2n+6}^{\otimes 6}$ is an injective linear map defined on monomials by $\psi(t) = t_1t_2t_3t_4t_5t_6t^2.$ Then, from the calculations in \cite{Mothebe}, we deduce that for each $1\leq d\leq 4,\, d\in\mathbb Z,$ if $t\in (\mathscr {C}^{\otimes 5}_{n+1-2^{d}})^{>0},$ then $t_l^{2^{d}-1}\mathsf{q}_{l}(t)\in (\mathscr {C}^{\otimes 6}_{n})^{>0}$ for any $1\leq l\leq 6.$ In other words, if $t$ is an admissible monomial of degree $n+1-2^{d}$ in the $\mathcal A$-module $P^{\otimes 5},$ then $t_l^{2^{d}-1}\mathsf{q}_{l}(t)$ is an admissible monomial of degree $n$ in $\mathcal A$-module $P^{\otimes 6}.$ 

We set $ (\mathscr C(d, n))^{>0}:= \big\{t_l^{2^{d}-1}\mathsf{q}_{l}(t)|\, t\in (\mathscr {C}^{\otimes 5}_{n+1-2^{d}})^{>0},\, 1\leq l\leq 6\big\},\ 1\leq d\leq 4.$ Notice that when $n = n_1$ and $h = 6$, Kameko's maps can be rewritten as $ 
 (\widetilde {Sq^0_*})_{n_1}: QP^{\otimes h}_{n_1}  \longrightarrow QP^{\otimes h}_{\frac{n_1-h}{2}},\ \  \psi: P_{\frac{n_1-h}{2}}^{\otimes h}\longrightarrow P_{n_1}^{\otimes h}.$ So, $\psi\big(\mathscr {C}^{\otimes 6}_{\frac{n_1-h}{2}}\big)\subset  (\mathscr C(d, n_1))^{>0}$ and $\widetilde {Sq^0_*}([u]) = [0]$ for all $u\in (\mathscr C(d, n_1))^{>0}\setminus \psi\big(\mathscr {C}^{\otimes 6}_{\frac{n_1-h}{2}}\big).$ These lead to $[u]\in {\rm Ker}((\widetilde {Sq^0_*})_{n_1}).$ According to the works in \cite{Sum3, Tin}, we have
$$ 
 \big|(\mathscr {C}^{\otimes 5}_{n_1+1-2^{d}})^{>0}\big| =\left \{\begin{array}{ll}
720 &\mbox{if $d = 1$},\\[1mm]
610 &\mbox{if $d = 2$},\\[1mm]
642 &\mbox{if $d = 3$},\\[1mm]
75&\mbox{if $d = 4$}.
\end{array}\right.$$
Thanks to the results, a straightforward computation yields 
$$ \bigcup_{1\leq d\leq 4}(\mathscr C(d, n_1))^{>0} = D_1\bigcup D_2 \bigcup E,$$
where
$$ \begin{array}{ll}
&D_1 = \big\{c_j:\ 1\leq j\leq 336\big\},\ \ D_1\subset (\mathscr {C}^{\otimes 6}_{n_1})^{>0}(4,5,1,1),\\[1mm]
&D_2 = \big\{d_j:\ 1\leq j\leq 210 \big\},\ \ D_2\subset (\mathscr {C}^{\otimes 6}_{n_1})^{>0}(4,5,3),\\[1mm]
&E\subset \big((\mathscr {C}^{\otimes 6}_{n_1})^{>0}(4,3,2,1)\bigcup (\mathscr {C}^{\otimes 6}_{n_1})^{>0}(4,3,4)\bigcup \psi(\mathscr {C}^{\otimes 6}_{n_0})\big),\\[1mm]
\end{array}$$
and the admissible monomials $c_j$ and $d_j$ are respectively described by Appendices \ref{s51} and \ref{s52}. The results obtained are based on filtering and removing the same monomials. For each monomial $u\in D_1,$ we notice that for any $[t]_{(4,5,1,1)}\in (QP^{\otimes 6}_{n_1})^{>0}(4,5,1,1),$ if $t$ is not equal to $u,$ then $t$ must take one of the following forms:
\begin{enumerate}
\item[$*$] $t_i^{3}t_j^{15}t_k^{2}t_{\ell}t_m^{2}t_p^{3}$,\  $t_i^{3}t_j^{15}t_k^{2}t_{\ell}^{2}t_mt_p^{3}$,\ $t_i^{3}t_j^{2}t_kt_{\ell}^{2}t_m^{7}t_p^{11}$,\ $t_i^{3}t_j^{2}t_k^{2}t_{\ell}t_m^{7}t_p^{11}$ where $j < k < \ell,$ and $(i, j, k, \ell, m, p)$ is a premutation of $(1, 2, 3, 4, 5, 6);$
\medskip

\item [$*$] $t_i^{3}t_j^{2}t_kt_{\ell}^{14}t_m^{3}t_p^{3}$,\ $t_i^{3}t_j^{2}t_k^{14}t_{\ell}t_m^{3}t_p^{3}$,\ $t_i^{3}t_j^{14}t_kt_{\ell}^{2}t_m^{3}t_p^{3}$,\ $t_i^{3}t_j^{14}t_k^{2}t_{\ell}t_m^{3}t_p^{3}$,\ $t_1t_2^{2}t_3^{14}t_{4}^{3}t_5^{3}t_6^{3},$  where $(i, j, k, \ell, m, p)$ is a premutation of $(1, 2, 3, 4, 5, 6)$ and $j < k < \ell;$ 

\medskip

\item [$*$] $t_i^{3}t_j^{2}t_kt_{\ell}^{7}t_m^{3}t_p^{10},$\ $t_i^{3}t_j^{2}t_k^{7}t_{\ell}t_m^{3}t_p^{10},$\ $t_i^{3}t_j^{7}t_kt_{\ell}^{2}t_m^{3}t_p^{10},$\ $t_i^{3}t_j^{7}t_k^{2}t_{\ell}t_m^{3}t_p^{10},$\ $t_i^{3}t_j^{2}t_k^{5}t_{\ell}^{10}t_m^{3}t_p^{3},$\ $t_i^{3}t_j^{2}t_k^{6}t_{\ell}^{9}t_m^{3}t_p^{3},$\ $t_i^{3}t_j^{6}t_k^{2}t_{\ell}^{9}t_m^{3}t_p^{3},$\ $t_i^{3}t_j^{6}t_k^{9}t_{\ell}^{2}t_m^{3}t_p^{3},$ where $j < k < \ell,$ and $(i, j, k, \ell, m, p)$ is a premutation of $(1, 2, 3, 4, 5, 6);$

\medskip

\item [$*$] $t_i^{3}t_j^{3}t_k^{6}t_{\ell}t_m^{3}t_p^{10},$\ $t_1t_2^{6}t_3^{3}t_{4}^{3}t_5^{3}t_6^{10},$\ $t_1t_2^{6}t_3^{3}t_{4}^{3}t_5^{10}t_6^{3},$\  $t_1t_2^{6}t_3^{3}t_{4}^{10}t_5^{3}t_6^{3},$\  $t_1t_2^{6}t_3^{10}t_{4}^{3}t_5^{3}t_6^{3},$\ $t_i^{3}t_j^{2}t_k^{2}t_{\ell}^{5}t_m^{3}t_p^{11},$\ $t_i^{3}t_j^{2}t_k^{5}t_{\ell}^{2}t_m^{3}t_p^{11},$\\ $t_i^{3}t_j^{2}t_k^{13}t_{\ell}^{2}t_m^{3}t_p^{3},$\ $t_i^{3}t_j^{2}t_k^{2}t_{\ell}^{13}t_m^{3}t_p^{3},$\ $t_i^{3}t_j^{2}t_k^{2}t_{\ell}^{9}t_m^{7}t_p^{3},$\ $t_i^{3}t_j^{2}t_k^{9}t_{\ell}^{2}t_m^{7}t_p^{3},$ where $j < k < \ell,$ and $(i, j, k, \ell, m, p)$ is a premutation of $(1, 2, 3, 4, 5, 6);$

\medskip

\item [$*$] $t_1t_j^{2}t_k^{6}t_{\ell}^{3}t_m^{3}t_p^{11},$\  $t_1t_j^{6}t_k^{2}t_{\ell}^{3}t_m^{3}t_p^{11},$ $t_1t_j^{6}t_k^{3}t_{\ell}^{2}t_m^{3}t_p^{11},$ $t_1^{3}t_j^{2}t_kt_{\ell}^{6}t_m^{3}t_p^{11},$\ $t_1^{3}t_j^{2}t_k^{6}t_{\ell}t_m^{3}t_p^{11},$ $t_1^{3}t_j^{6}t_kt_{\ell}^{2}t_m^{3}t_p^{11},$\ $t_1^{3}t_j^{6}t_k^{2}t_{\ell}t_m^{3}t_p^{11}$ with $j < k < \ell,$ and $(j, k, \ell, m, p)$ is a premutation of $(2, 3, 4, 5, 6);$

\medskip

\item [$*$] $t_i^{3}t_j^{3}t_k^{3}t_{\ell}^{3}t_m^{2}t_p^{12},$\ $t_i^{3}t_j^{3}t_k^{3}t_{\ell}^{2}t_m^{4}t_p^{11},$\ $t_i^{3}t_j^{3}t_k^{3}t_{\ell}^{7}t_m^{2}t_p^{8},$\ $t_i^{3}t_j^{3}t_k^{3}t_{\ell}^{3}t_m^{4}t_p^{10},$\ $t_i^{3}t_j^{3}t_k^{3}t_{\ell}^{3}t_m^{6}t_p^{8},$ where $(i, j, k, \ell, m, p)$ is a premutation of $(1, 2, 3, 4, 5, 6).$
\end{enumerate}

It is straightforward to check that these monomials are strictly inadmissible, and so, they are inadmissible. To exemplify, let us consider the monomials $t_it_j^{6}t_k^{3}t_{\ell}^{3}t_m^{3}t_p^{10}$ and $t_i^{3}t_j^{3}t_k^{3}t_{\ell}^{3}t_m^{2}t_p^{12}.$ It is clear that $\omega(t_it_j^{6}t_k^{3}t_{\ell}^{3}t_m^{3}t_p^{10}) = \omega(t_i^{3}t_j^{3}t_k^{3}t_{\ell}^{3}t_m^{2}t_p^{12}) = (4,5,1,1).$ As well known, the action of the Steenrod algebra $\mathcal A$ on the polynomial algebra $P^{\otimes 6}$ is given by the rule
$$ Sq^k(t_j) = \left\{\begin{array}{ll}
t_j &\mbox{if $k = 0,$}\\[1mm]
t_j^2& \mbox{if $k = 1,$}\ \ \mbox{(\textit{the instability condition})},\\[1mm]
0 &\mbox{otherwise,} 
\end{array}\right.$$
and the Cartan formula $Sq^k(fg) = \sum_{a+b = k}Sq^{a}(f)Sq^{b}(g),$ for all $f,\, g\in P^{\otimes 6}.$ Note that for each $t\in P^{\otimes 1}$ and each positive integer $n,$ $Sq^{a}(t^{n}) = \binom{n}{a}t^{n+a},$ where the binomial coefficients $\binom{n}{a}$ are to be interpreted modulo 2 with the usual convention $\binom{n}{a} = 0$ if $n < a.$ Therefore, through a direct calculation, we obtain:
$$ \begin{array}{ll}
\medskip
t_it_j^{6}t_k^{3}t_{\ell}^{3}t_m^{3}t_p^{10} &=  Sq^{1}(t_i^{2}t_j^{3}t_k^{5}t_{\ell}^{3}t_m^{3}t_p^{9} + t_i^{2}t_j^{5}t_k^{3}t_{\ell}^{3}t_m^{3}t_p^{9} + t_i^{2}t_j^{3}t_k^{3}t_{\ell}^{3}t_m^{5}t_p^{9} +  t_i^{2}t_j^{3}t_k^{3}t_{\ell}^{5}t_m^{3}t_p^{9} )\\
\medskip
&\quad + Sq^{2}(t_it_j^{3}t_k^{5}t_{\ell}^{3}t_m^{3}t_p^{9} + t_it_j^{5}t_k^{3}t_{\ell}^{3}t_m^{3}t_p^{9} + t_it_j^{3}t_k^{3}t_{\ell}^{3}t_m^{5}t_p^{9} + t_it_j^{3}t_k^{3}t_{\ell}^{5}t_m^{3}t_p^{9})\\
\medskip
&\quad + t_it_j^{3}t_k^{6}t_{\ell}^{3}t_m^{3}t_p^{10} + t_it_j^{3}t_k^{3}t_{\ell}^{3}t_m^{6}t_p^{10} +  t_it_j^{3}t_k^{3}t_{\ell}^{6}t_m^{3}t_p^{10}  \mod P_{26}^{\otimes 6}(< (4,5,1,1)),\\
t_i^{3}t_j^{3}t_k^{3}t_{\ell}^{3}t_m^{2}t_p^{12}&=Sq^{1}(t_i^{3}t_j^{3}t_k^{3}t_{\ell}^{3}t_mt_p^{12}) \mod P_{n_1}^{\otimes 6}(< (4,5,1,1)),
\end{array}$$
which imply $t_it_j^{6}t_k^{3}t_{\ell}^{3}t_m^{3}t_p^{10}\equiv_{(4,5,1,1)} (t_it_j^{3}t_k^{6}t_{\ell}^{3}t_m^{3}t_p^{10} + t_it_j^{3}t_k^{3}t_{\ell}^{3}t_m^{6}t_p^{10} +  t_it_j^{3}t_k^{3}t_{\ell}^{6}t_m^{3}t_p^{10}),$ and  $t_i^{3}t_j^{3}t_k^{3}t_{\ell}^{3}t_m^{2}t_p^{12}\equiv_{(4,5,1,1)} 0.$ Hence, the monomials $ t_it_j^{6}t_k^{3}t_{\ell}^{3}t_m^{3}t_p^{10}$ and $t_i^{3}t_j^{3}t_k^{3}t_{\ell}^{3}t_m^{2}t_p^{12}$ are strictly inadmissible and $(4,5,1,1)$-hit, respectively. Since the monomials in $D_1$ are admissible, $\mathscr (C_{n_1}^{\otimes 6})^{>0}(4,5,1,1) = D_1.$ Thus, it may be claimed that $\dim (QP^{\otimes 6}_{n_1})^{>0}(4,5,1,1) = |D_1| = 336.$ Next,  we observe that for each monomial $\widehat{u}\in D_2,$ if $[\widehat{t}]_{(4,5,3)}\in (QP^{\otimes 6}_{n_1})^{>0}(4,5,3),$ and $\widehat{t}\neq \widehat{u},$ then $\widehat{t}$ must have one of the following forms:
\begin{enumerate}
\item[$-$] $t_i^{3}t_j^{2}t_k^{2}t_{\ell}^{5}t_m^{7}t_p^{7}$,\  $t_i^{3}t_j^{2}t_k^{5}t_{\ell}^{2}t_m^{7}t_p^{7}$,\ $t_i^{7}t_j^{2}t_kt_{\ell}^{2}t_m^{7}t_p^{7}$,\ $t_i^{7}t_j^{2}t_k^{2}t_{\ell}t_m^{7}t_p^{7},$\ $t_i^{3}t_j^{6}t_k^{5}t_{\ell}^{6}t_m^{3}t_p^{3},$\ $t_i^{3}t_j^{6}t_k^{6}t_{\ell}^{5}t_m^{3}t_p^{3},$ where $j < k < \ell,$ and $(i, j, k, \ell, m, p)$ is a premutation of $(1, 2, 3, 4, 5, 6);$

\medskip

\item[$-$] $t_it_j^{2}t_k^{6}t_{\ell}^{3}t_m^{7}t_p^{7},$\ $t_it_j^{6}t_k^{2}t_{\ell}^{3}t_m^{7}t_p^{7},$\ $t_it_j^{6}t_k^{3}t_{\ell}^{2}t_m^{7}t_p^{7},$\ $t_i^{3}t_j^{2}t_k^{6}t_{\ell}t_m^{7}t_p^{7},$\ $t_i^{3}t_j^{6}t_k^{2}t_{\ell}t_m^{7}t_p^{7},$ where $j < k < \ell,$ and $(i, j, k, \ell, m, p)$ is a premutation of $(1, 2, 3, 4, 5, 6);$

\medskip

\item[$-$] $t_i^{3}t_j^{6}t_kt_{\ell}^{6}t_m^{3}t_p^{7}$,\ $t_i^{3}t_j^{6}t_k^{6}t_{\ell}t_m^{3}t_p^{7},$\ $t_it_j^{6}t_k^{3}t_{\ell}^{6}t_m^{3}t_p^{7},$\ $t_it_j^{6}t_k^{6}t_{\ell}^{3}t_m^{3}t_p^{7},$\ $t_i^{3}t_j^{2}t_k^{5}t_{\ell}^{6}t_m^{3}t_p^{7},$\ $t_i^{3}t_j^{2}t_k^{6}t_{\ell}^{5}t_m^{3}t_p^{7},$\ $t_i^{3}t_j^{6}t_k^{2}t_{\ell}^{5}t_m^{3}t_p^{7},$\\ $t_i^{3}t_j^{6}t_k^{5}t_{\ell}^{2}t_m^{3}t_p^{7},$ where $j < k < \ell,$ and $(i, j, k, \ell, m, p)$ is a premutation of $(1, 2, 3, 4, 5, 6);$

\medskip

\item[$-$] $t_i^{3}t_j^{3}t_k^{7}t_{\ell}^{3}t_m^{4}t_p^{6}$,\ $t_i^{3}t_j^{3}t_k^{7}t_{\ell}^{7}t_m^{2}t_p^{4},$ where $(i, j, k, \ell, m, p)$ is a premutation of $(1, 2, 3, 4, 5, 6).$

\end{enumerate}

It is indeed facile to simply assert that these monomials are inadmissible. As an illustration, let us consider the monomials $t_it_j^{2}t_k^{6}t_{\ell}^{3}t_m^{7}t_p^{7}$ and $t_i^{3}t_j^{3}t_k^{7}t_{\ell}^{3}t_m^{4}t_p^{6}.$ Then, $\omega(t_it_j^{2}t_k^{6}t_{\ell}^{3}t_m^{7}t_p^{7}) = \omega(t_i^{3}t_j^{3}t_k^{7}t_{\ell}^{3}t_m^{4}t_p^{6}) = (4,5,3)$ and by a simple computation, we get
$$ \begin{array}{ll}
\medskip
 t_it_j^{2}t_k^{6}t_{\ell}^{3}t_m^{7}t_p^{7} &= Sq^{2}(t_it_jt_k^{3}t_{\ell}^{5}t_m^{7}t_p^{7} + t_it_jt_k^{5}t_{\ell}^{3}t_m^{7}t_p^{7} + t_it_jt_k^{3}t_{\ell}^{3}t_m^{7}t_p^{9} + t_it_jt_k^{3}t_{\ell}^{3}t_m^{9}t_p^{7})\\
\medskip
&\quad + Sq^{1}( t_i^{2}t_jt_k^{3}t_{\ell}^{5}t_m^{7}t_p^{7} + t_i^{2}t_jt_k^{5}t_{\ell}^{3}t_m^{7}t_p^{7})  +  t_it_j^{2}t_k^{3}t_{\ell}^{6}t_m^{7}t_p^{7} \mod P_{n_1}^{\otimes 6}(< (4,5,3)),\\
t_i^{3}t_j^{3}t_k^{7}t_{\ell}^{3}t_m^{4}t_p^{6}&= Sq^{1}(t_i^{3}t_j^{3}t_k^{7}t_{\ell}^{3}t_m^{4}t_p^{5}) \mod P_{n_1}^{\otimes 6}(< (4,5,3)).
\end{array} $$
These equalities show that $ t_it_j^{2}t_k^{6}t_{\ell}^{3}t_m^{7}t_p^{7}$ is strictly inadmissible (since $t_it_j^{2}t_k^{3}t_{\ell}^{6}t_m^{7}t_p^{7}< t_it_j^{2}t_k^{6}t_{\ell}^{3}t_m^{7}t_p^{7}$) and that $t_i^{3}t_j^{3}t_k^{7}t_{\ell}^{3}t_m^{4}t_p^{6}$ is $(4,5,3)$-hit.

\medskip

Since the monomials in $D_2$ are admissible, $\mathscr (C_{n_1}^{\otimes 6})^{>0}(4,5,3) = D_2.$ So, $\dim (QP^{\otimes 6}_{n_1})^{>0}(4,5,3) = |D_2| = 210.$ Incorporating the above-mentioned computation, we obtain
$$ \dim U_1 =  |D_1| + |D_2| = 336 + 210 = 546.$$

\medskip

The next step is to determine the dimension of $U_2$.  Directly computing the dimension of this cohit module is a task of considerable complexity. Our calculations show that
$$ |\mathscr (C_{n_1}^{\otimes 6})^{>0}(4,3,2,1)| = 2880, \ \ \mbox{and}\ \   |\mathscr (C_{n_1}^{\otimes 6})^{>0}(4,3,4)| = 210.$$
Consequently, $$ \dim U_2 = 2880 + 210 = 3090.$$ 
Nevertheless, the methods employed to compute it are akin to those utilized in our earlier publications \cite{P.S1, Phuc4, Phuc6, Phuc10}. 

\medskip

\medskip

Thus, from the above calculations, ${\rm Ker}((\widetilde {Sq^0_*})_{n_1})\cap (QP^{\otimes 6}_{n_1})^{>0}$ is an $\mathbb F_2$-vector of dimension $3636.$ Finally, as $QP^{\otimes 6}_{n_1} \cong  (QP^{\otimes 6}_{n_1})^{0}\bigoplus \big({\rm Ker}((\widetilde {Sq^0_*})_{n_1})\cap (QP^{\otimes 6}_{n_1})^{>0}\big)\bigoplus QP^{\otimes 6}_{n_0},$ we conclude that 
$$ \begin{array}{ll}
\medskip
\dim QP_{n_1}^{\otimes 6} &= \dim (QP_{n_1}^{\otimes 6})^{0} + \dim {\rm Ker}((\widetilde {Sq^0_*})_{n_1})\cap (QP^{\otimes 6}_{n_1})^{>0} + \dim QP_{n_0}^{\otimes 6}\\
&=\dim (QP_{n_1}^{\otimes 6})^{0} + \dim U_1 + \dim U_2 + \dim QP_{n_0}^{\otimes 6}\\
& =  5184 + 546 + 3090 + 945 = 9765.
\end{array}$$ 
The proof of the theorem is complete.

To close this subsection, we would like to provide some insightful observations and remarks regarding the indecomposables $Q^{\otimes h}_{11}.$

\begin{nx}
\begin{itemize}

\item[(i)] With Corollary \ref{hqT} and a result from \cite{Tin} concerning $\dim Q^{\otimes 5}_{11}$ as our basis, we assert that the localized version of Kameko's conjecture in Note \ref{cyP}(i) holds true for all $h$ and degree $11$.

\item[(ii)] Clearly, Corollary \ref{hqT} implies that the coinvariant $(\mathbb F_2\otimes_{GL_h}{\rm Ann}_{\overline{\mathcal A}}[P^{\otimes h}]^{*})_{11}$ is trivial for all $h\geq 12.$ Hence, Singer's Conjecture \ref{gtSinger} is true for bidegrees $(h, h+11)$ with $h>11.$ Moreover, in \cite{Phuc11}, we have demonstrated that the conjecture is also true for $6\leq h\leq 8.$ This is achieved by explicitly computing the dimensions of the invariants $[{\rm Ker}((\widetilde {Sq^0_*})_{11})]^{GL_7},\, [QP^{\otimes 7}_2]^{GL_7},$ and $[QP^{\otimes h}_{11}(\widehat{\omega}_{(j)})]^{GL_h}$ for $h = 6,\, 8,\, 1\leq j\leq 5.$ We then prove that the transfer homomorphism $Tr_h^{\mathcal A}$ is a monomorphism if $6\leq h\leq 7$ and is a trivial isomorphism if $h = 8.$ Here the weight vectors $\widehat{\omega}_{(j)}$ are given as in Corollary \ref{hqT}. It is noteworthy to mention that, the works of Ch\ohorn n and H\`a \cite{CHa, CHa0} demonstrated the non-surjectivity of the transfer in the bidegrees $(6, 6+11)$ and $(7, 7+11).$

According to the research conducted by \cite{Tangora, Bruner, Bruner2, Lin2}, it can be concluded that for every $h\geq 6,$
$$
{\rm Ext}_{\mathcal A}^{h, h+11}(\mathbb F_2, \mathbb F_2) =  \left\{\begin{array}{ll}
\mathbb F_2h_0Ph_2&\mbox{if $h = 6$},\\[1mm]
\mathbb F_2h_0^{2}Ph_2 = \mathbb F_2 h_1^{2}Ph_1&\mbox{if $h = 7$},\\[1mm]
0&\mbox{if $h\geq 8$}.
\end{array}\right.$$

For $h = 9,$ according to Corollary \ref{hqT}, we have $QP^{\otimes 9}_{1} = (QP^{\otimes 9}_{1})^{0} = \langle \{[t_i]\}_{1\leq i\leq 9}\rangle,$  and $QP^{\otimes 9}_{11}(\widehat{\omega}_{(j)}) = (QP^{\otimes 9}_{11})^{0}(\widehat{\omega}_{(j)})$ for $1\leq j\leq 4.$ So the invariants $[QP^{\otimes 9}_{1}]^{GL_9}$ and $QP^{\otimes 9}_{11}(\widehat{\omega}_{(j)}),\, 1\leq j\leq 4$ are trivial. By combining these data with Corollary \ref{hqT} and taking into account that the mapping $(\widetilde {Sq^0_*})_{11}: QP^{\otimes 9}_{11}\longrightarrow QP^{\otimes 9}_2$ is a surjective, we can derive an estimate $$\dim (\mathbb F_2\otimes_{GL_{9}}{\rm Ann}_{\overline{\mathcal A}}[P^{\otimes 9}]^{*})_{11}=\dim [QP^{\otimes 9}_{11}]^{GL_9}\leq \dim [{\rm Ker}((\widetilde {Sq^0_*})_{11})]^{GL_9}\leq \dim [QP^{\otimes 9}_{11}(\widehat{\omega}_{(5)})]^{GL_9}.$$
Owing to Corollary \ref{hqT}, one has an isomorphism $QP^{\otimes 9}_{11}(\widehat{\omega}_{(5)})\cong (QP^{\otimes 9}_{11})^{0}(\widehat{\omega}_{(5)})\bigoplus (QP^{\otimes 9}_{11})^{>0}(\widehat{\omega}_{(5)})$ where $\dim (QP^{\otimes 9}_{11})^{0}(\widehat{\omega}_{(5)}) = 21\binom{9}{7} + 48\binom{9}{8} = 1188$ and $\dim (QP^{\otimes 9}_{11})^{>0}(\widehat{\omega}_{(5)}) = 27.$ 

When $h=10$, it follows from Corollary \ref{hqT} that the invariants $QP^{\otimes 10}_{11}(\widehat{\omega}_{(j)})$ are trivial for $1\leq j\leq 5.$ As a result, we can establish an inequality $$\dim (\mathbb F_2\otimes_{GL_{10}}{\rm Ann}_{\overline{\mathcal A}}[P^{\otimes 10}]^{*})_{11}=\dim [QP^{\otimes 10}_{11}]^{GL_{10}}\leq \dim [QP^{\otimes 10}_{11}(\widehat{\omega}_{(6)})]^{GL_{10}}.$$
According to Corollary \ref{hqT}, $QP^{\otimes 10}_{11}(\widehat{\omega}_{(6)})\cong (QP^{\otimes 10}_{11})^{0}(\widehat{\omega}_{(6)})\bigoplus (QP^{\otimes 10}_{11})^{>0}(\widehat{\omega}_{(6)})$ where $$\dim (QP^{\otimes 10}_{11})^{0}(\widehat{\omega}_{(6)}) = 9\binom{10}{9} =90\ \mbox{and}\ \dim (QP^{\otimes 10}_{11})^{>0}(\widehat{\omega}_{(6)}) = 9.$$ 

In the case where $h=11$, it is possible to derive from Corollary \ref{hqT} that $$QP^{\otimes 11}_{11}\cong \bigoplus_{1\leq j\leq 6}(QP^{\otimes 11}_{11})^{0}(\widehat{\omega}_{(j)})\bigoplus (QP^{\otimes 11}_{11})^{>0}(\widehat{\omega}_{(7)}),$$ where $(QP^{\otimes 11}_{11})^{>0}(\widehat{\omega}_{(7)}) \cong \mathbb F_2[t_1t_2\ldots t_{11}]_{\widehat{\omega}_{(7)}}.$ Using the $\mathcal A$-homomorphisms $\sigma_d: P^{\otimes 11}\longrightarrow P^{\otimes 11},$ for $1\leq d\leq 11,$ one gets $\dim (\mathbb F_2\otimes_{GL_{11}}{\rm Ann}_{\overline{\mathcal A}}[P^{\otimes 11}]^{*})_{11} =\dim [QP^{\otimes 11}_{11}]^{GL_{11}} = 0 = \dim {\rm Ext}_{\mathcal A}^{11, 22}(\mathbb F_2, \mathbb F_2).$ Hence, Singer's transfer is a trivial isomorphism in bidegree $(11, 22).$ Thus if the invariants $[QP^{\otimes 9}_{11}(\widehat{\omega}_{(5)})]^{GL_9}$ and $[QP^{\otimes 10}_{11}(\widehat{\omega}_{(6)})]^{GL_{10}}$ are trivial, then Singer's Conjecture \ref{gtSinger} also holds for bidegrees $(h, h+11)$ with $9\leq h\leq 11.$ This matter will be thoroughly investigated and discussed in another context.
\end{itemize}
\end{nx}

 \subsection{Proof of Theorem \ref{dlc4}}

In what follows, suppose that $\omega$ is a weight vector of degree $n_1.$ We denote by $\mathscr C^{\otimes 6}_{n_1}(\omega)$ the set of all admissible monomials in $P^{\otimes 6}_{n_1}(\omega)$ and by $[\mathscr C^{\otimes 6}_{n_1}(\omega)]_{\omega} = \{[t]_{\omega}:\ t\in \mathscr C^{\otimes 6}_{n_1}(\omega)\}.$ For $z_1, z_2,\ldots, z_m\in P^{\otimes 6}_{n_1}(\omega)$ with $m\geq 1,$ we put
$$ \begin{array}{ll}
\medskip
\Sigma_6(z_1, \ldots, z_m) &= \big\{\theta(z_j):\ \theta\in \Sigma_6,\, 1\leq j\leq m\big\},\\
\medskip
[\mathscr C(z_1, \ldots, z_m)]_{\omega}&= [\mathscr C^{\otimes 6}_{n_1}(\omega)]_{\omega}\cap \langle [\Sigma_6(z_1, \ldots, z_m)]_{\omega} \rangle,\\
\widehat{p(z)}&= \sum_{y\in \mathscr C^{\otimes 6}_{n_1}(\omega)\cap \Sigma_6(z)}y.
\end{array}$$
$\langle [\Sigma_6(z_1, \ldots, z_m)]_{\omega} \rangle$ is manifestly a $\Sigma_6$-submodule of $QP^{\otimes 6}_{n_1}(\omega).$ As we have pointed out before,
$$ 
 {\rm Ker}((\widetilde {Sq^0_*})_{n_1})\cong (QP_{n_1}^{\otimes 6})(4,3,2,1)\bigoplus QP_{n_1}^{\otimes 6}(4,3,4)\bigoplus QP_{n_1}^{\otimes 6}(4,5,1,1)\bigoplus QP_{n_1}^{\otimes 6}(4,5,3),
$$
where 
$$ \begin{array}{ll}
\medskip
&(QP_{n_1}^{\otimes 6})(4,3,4) = (QP_{n_1}^{\otimes 6})^{>0}(4,3,4),\ \ QP_{n_1}^{\otimes 6}(4,5,1,1) = (QP_{n_1}^{\otimes 6})^{>0}(4,5,1,1),\\
&QP_{n_1}^{\otimes 6}(4,5,3) = (QP_{n_1}^{\otimes 6})^{>0}(4,5,3).
\end{array}$$ 
So, one gets an estimate
$$\begin{array}{ll}
\dim {\rm [}{\rm Ker}((\widetilde {Sq^0_*})_{n_1}){\rm]}^{GL_6}&\leq \dim {\rm [}QP_{n_1}^{\otimes 6}(4,3,2,1){\rm ]}^{GL_6} +\dim {\rm [}QP_{n_1}^{\otimes 6}(4,3,4){\rm ]}^{GL_6} \\
&\quad + \dim {\rm [}QP_{n_1}^{\otimes 6}(4,5,1,1){\rm ]}^{GL_6} + \dim {\rm [}QP_{n_1}^{\otimes 6}(4,5,3){\rm ]}^{GL_6}.
\end{array}$$
By using the monomial basis of the space $QP_{n_1}^{\otimes 6}(\omega)$ with $\omega\in \{(4,3,2,1), (4,3,4), (4,5,1,1), (4,5,3)\}$ (see Theorem \ref{dlc1}) and the homomorphisms $\sigma_d$ for $1\leq d\leq 6,$ we find that the invariants $[QP_{n_1}^{\otimes 6}(\omega)]^{GL_6}$ are zero. Indeed, we will prove this claim for the invariants $[QP_{n_1}^{\otimes 6}(4,5,1,1)]^{GL_6}$ and $[QP_{n_1}^{\otimes 6}(4,5,3)]^{GL_6}$ in detail. Similarly, we also obtain the results for the other spaces.
\medskip

We set $\widetilde{\omega}:= (4,5,1,1)$ and $\omega:= (4,5,3).$ We first describe the space $[QP_{n_1}^{\otimes 6}(\widetilde{\omega})]^{GL_6}.$ Following the proof of Theorem \ref{dlc1}, the space $QP_{n_1}^{\otimes 6}(\widetilde{\omega}) = (QP_{n_1}^{\otimes 6})^{>0}(\widetilde{\omega})$ is $336$-dimensional with the monomial basis $\{[c_j]_{\widetilde{\omega}}:\ 1\leq j\leq 336\}.$ Note that the admissible monomials $c_j$ are explicitly described as in Subsect. \ref{s51}. Let us consider the following admissible monomials:
$$ \begin{array}{ll} 
\medskip
c_1 &= t_1^{3}t_2t_3^{15}t_4^{2}t_5^{2}t_6^{3}, \ \ c_{61} = t_1^{3}t_2^{7}t_3^{11}t_4t_5^{2}t_6^{2},\ \ c_{121} = t_1t_2^{3}t_3^{14}t_4^{2}t_5^{2}t_6^{3},\\
\medskip
c_{156} &= t_1^{3}t_2^{3}t_3^{13}t_4^{2}t_5^{2}t_6^{3}, \ \ c_{166} = t_1^{3}t_2t_3^{2}t_4^{3}t_5^{6}t_6^{11},\ \ c_{196} = t_1^{3}t_2^{5}t_3^{11}t_4^{2}t_5^{2}t_6^{3},\\
c_{211} &= t_1t_2^{7}t_3^{10}t_4^{2}t_5^{3}t_6^{3}, \ \ c_{286} = t_1^{3}t_2^{7}t_3^{9}t_4^{2}t_5^{2}t_6^{3},\ \ c_{301} = t_1^{3}t_2t_3^{3}t_4^{3}t_5^{6}t_6^{10}.
\end{array}$$
The following spaces are $\Sigma_6$-submodules of $QP_{n_1}^{\otimes 6}(\widetilde{\omega})$:
$$ \begin{array}{ll}
\medskip
\langle[\Sigma_6(c_{1})]_{\widetilde{\omega}}\rangle &= \langle \{[c_j]_{\widetilde{\omega}}:\ 1\leq j\leq 60\} \rangle,\ \ \langle[\Sigma_6(c_{61})]_{\widetilde{\omega}}\rangle = \langle \{[c_j]_{\widetilde{\omega}}:\ 61\leq j\leq 120\} \rangle,\\
\medskip
\langle[\Sigma_6(c_{121},\, c_{211})]_{\widetilde{\omega}}\rangle &= \langle \{[c_j]_{\widetilde{\omega}}:\ 121\leq j\leq 165,\, 211\leq j\leq 300\} \rangle,\\
 \langle[\Sigma_6(c_{166})]_{\widetilde{\omega}}\rangle &= \langle \{[c_j]_{\widetilde{\omega}}:\ 166\leq j\leq 210\} \rangle,\ \ \langle[\Sigma_6(c_{301})]_{\widetilde{\omega}}\rangle = \langle \{[c_j]_{\widetilde{\omega}}:\ 301\leq j\leq 336\} \rangle.
\end{array}$$
So, we have an isomorphism
$$ \begin{array}{ll}
\medskip
 QP_{n_1}^{\otimes 6}(\widetilde{\omega}) &\cong \langle[\Sigma_6(c_{1})]_{\widetilde{\omega}}\rangle\bigoplus \langle[\Sigma_6(c_{61})]_{\widetilde{\omega}}\rangle \bigoplus \langle[\Sigma_6(c_{121},\, c_{211})]_{\omega}\rangle\\
&\quad\bigoplus \langle[\Sigma_6(c_{166})]_{\widetilde{\omega}}\rangle \bigoplus \langle[\Sigma_6(c_{301})]_{\widetilde{\omega}}\rangle.
\end{array}$$
By direct calculations using the homomorphisms $\sigma_i: P^{\otimes 6}\longrightarrow P^{\otimes 6}$ for $1\leq i\leq 5,$ we obtain the following results:
$$ \begin{array}{ll}
\medskip
\langle[\Sigma_6(c_{j})]_{\widetilde{\omega}}\rangle^{\Sigma_6} &= 0,\ \ \mbox{for $j = 166,$}\\
\langle[\Sigma_6(c_{j})]_{\widetilde{\omega}}\rangle^{\Sigma_6} &= \langle [\widehat{p(c_j)}]_{\widetilde{\omega}}\rangle,\ \ \mbox{for $j = 1,\, 61,\, 301$},\\
\langle[\Sigma_6(c_{121},\, c_{211})]_{\widetilde{\omega}}\rangle^{\Sigma_6} &\langle [q]_{\widetilde{\omega}}\rangle, \ \mbox{where}\\ 
&q = \sum_{121\leq j\leq 165}{\rm c}_j + \sum_{211\leq j\leq 216}{\rm c}_j + \sum_{254\leq j\leq 261}{\rm c}_j +  \sum_{265\leq j\leq 285}{\rm c}_j + \sum_{289\leq j\leq 291}{\rm c}_j + \sum_{294\leq j\leq 300}{\rm c}_j.
\end{array}$$
where $\widehat{p(c_j)} = \sum_{c_j\in \mathscr C(c_j)} c_j$ with $\mathscr C(c_1) = \{c_j:\ 1\leq j\leq 60\},$\ $\mathscr C(c_{61}) = \{c_j:\ 61\leq j\leq 120\},$ and $\mathscr C(c_{301}) = \{c_j:\ 301\leq j\leq 336\}.$ Note that the sets $[\mathscr C(c_j)]_{\widetilde{\omega}}$ are the bases of the spaces $\langle[\Sigma_6(c_{j})]_{\widetilde{\omega}}\rangle$ for $j  = 1,\, 61,\, 211,\, 301.$ Thus, one gets $$[QP_{n_1}^{\otimes 6}(\widetilde{\omega})]^{\Sigma_6} = \langle \{[\widehat{p(c_j)}]_{\widetilde{\omega}},\, [q]_{\widetilde{\omega}}:\ j  =1,\, 61,\, 301\}\rangle.$$ Now, assume that $[g]_{\widetilde{\omega}}\in [QP_{n_1}^{\otimes 6}(\widetilde{\omega})]^{GL_6},$ then because $\Sigma_6\subset GL_6,$ we must have that $$g\equiv_{\widetilde{\omega}}\gamma_1\widehat{p(c_1)} + \gamma_2\widehat{p(c_{61})} + \gamma_3\widehat{p(c_{301})} + \gamma_4q, \,\, \gamma_i\in \mathbb F_2,\,\, 1\leq i\leq 4.$$ By a simple computation using the homomorphism $\sigma_6$ and the relation $\sigma_6(g) +g\equiv_{\widetilde{\omega}}0,$ we get $\gamma_i = 0,\, \forall i.$  Hence, $[QP_{n_1}^{\otimes 6}(\widetilde{\omega})]^{GL_6} = 0.$
\medskip

Next, we compute the space $[QP_{n_1}^{\otimes 6}(\omega)]^{GL_6}.$ Following the proof of Theorem \ref{dlc1}, $\dim QP_{n_1}^{\otimes 6}(\omega) = 210,$ and $QP_{n_1}^{\otimes 6}(\omega)  = \langle \{[d_j]_{\omega}:\ 1\leq j\leq 210\}\rangle,$ where the admissible monomials $d_j$ are explicitly described as in Subsect. \ref{s52}. By a simple computation, we have a direct summand decomposition of the $\Sigma_6$-submodules:
$$ QP_{n_1}^{\otimes 6}(\omega) = \langle[\Sigma_6(d_1)]_{\omega}\rangle\bigoplus \langle[\Sigma_6(d_{21})]_{\omega}\rangle \bigoplus \langle[\Sigma_6(d_{111})]_{\omega}\rangle\bigoplus \langle[\Sigma_6(d_{201})]_{\omega}\rangle,$$
where $$ \begin{array}{ll}
\medskip
 \langle[\Sigma_6(d_1)]_{\omega}\rangle\ &=\langle \{[d_j]_{\omega}:\ 1\leq j\leq 20\} \rangle, \ \ \ \langle[\Sigma_6(d_{21})]_{\omega}\rangle =\langle \{[d_j]_{\omega}:\ 21\leq j\leq 110\} \rangle,\\
\langle[\Sigma_6(d_{111})]_{\omega}\rangle &=\langle \{[d_j]_{\omega}:\ 111\leq j\leq 200\} \rangle, \ \ \ \langle[\Sigma_6(d_{201})]_{\omega}\rangle =\langle \{[d_j]_{\omega}:\ 201\leq j\leq 210\} \rangle.
\end{array}$$
We first compute the action of the symmetric group $\Sigma_6$ on $ QP_{n_1}^{\otimes 6}(\omega).$ We find that $$[QP_{n_1}^{\otimes 6}(\omega)]^{\Sigma_6} = \big\langle \{[\sum_{1\leq j\leq 20}d_j]_{\omega},\ [\sum_{21\leq j\leq 110}d_j]_{\omega},\ [\sum_{111\leq j\leq 200}d_j]_{\omega}\}\big \rangle.$$ This is immediate from the following assertions: 

\begin{enumerate}
\item[i)] $ \langle[\Sigma_6(d_1)]_{\omega}\rangle^{\Sigma_6} = \langle [\widehat{p(d_1)}]_{\omega}$ with $\widehat{p(d_1)}:= \sum_{1\leq j\leq 20}d_j;$ 
\medskip

\item[ii)] $\langle[\Sigma_6(d_{21})]_{\omega}\rangle^{\Sigma_6} = \langle [\widehat{p(d_{21})}]_{\omega} \rangle$ with $\widehat{p(d_{21})}:= \sum_{21\leq j\leq 110}d_j;$
\medskip
	
\item[iii)]  $\langle[\Sigma_6(d_{111})]_{\omega}\rangle^{\Sigma_6} = \langle [\widehat{p(d_{111})}]_{\omega} \rangle$ with $\widehat{p(d_{111})}:= \sum_{111\leq j\leq 200}d_j;$
\medskip

\item[iv)]  $\langle[\Sigma_6(d_{201})]_{\omega}\rangle^{\Sigma_6} = 0.$
\end{enumerate}

We compute the cases i) and iv) and leave the rest to the reader. It is straightforward to see that the sets $[\mathscr C(d_1)]_{\omega} = \{[d_j]_{\omega}:\ 1\leq j\leq 20\}$ and $[\mathscr C(d_{201})]_{\omega} =\{[d_j]_{\omega}:\ 201\leq j\leq 210\}$ are the bases of the spaces $\langle[\Sigma_6(d_1)]_{\omega}\rangle$ and $\langle[\Sigma_6(d_{201})]_{\omega}\rangle,$ respectively. Suppose that $[f]_{\omega}\in \langle[\Sigma_6(d_1)]_{\omega}\rangle^{\Sigma_6}$ and $[g]_{\omega}\in \langle[\Sigma_6(d_{201})]_{\omega}\rangle^{\Sigma_6}.$ Then, one has that $f\equiv_{\omega}\sum_{1\leq j\leq 20}\gamma_jd_j,$ and $g\equiv_{\omega}\sum_{201\leq j\leq 210}\beta_jd_j,$ in which the coefficients $\gamma_j$ and $\beta_j$ belong to $\mathbb F_2$ for all $j.$ By direct calculations using the homomorphisms $\sigma_i: P^{\otimes 6}\longrightarrow P^{\otimes 6}$ for $1\leq i\leq 5,$ and the relations $\sigma_i(f) + f\equiv_{\omega}0,$ and $\sigma_i(g) + g\equiv_{\omega}0,$ we obtain the following equalities:

\begin{align*}
\sigma_1(f) +f& \equiv_{\omega} (\gamma_5 + \gamma_{11})(d_5 + d_{11}) +  (\gamma_6 + \gamma_{12})(d_6 + d_{12}) +  (\gamma_7 + \gamma_{13})(d_7 + d_{13}) \\
\medskip
&\quad +  (\gamma_8 + \gamma_{14})(d_8 + d_{14})+  (\gamma_9 + \gamma_{15})(d_9 + d_{15}) +  (\gamma_{10} + \gamma_{16})(d_{10} + d_{16})  \equiv_{\omega} 0,\\
\sigma_2(f) + f&\equiv_{\omega} (\gamma_2 + \gamma_{5})(d_2 + d_{5}) + (\gamma_3 + \gamma_{6})(d_3 + d_{6}) + (\gamma_4 + \gamma_{7})(d_4 + d_{7})\\
\medskip
 &\quad +  (\gamma_{14} + \gamma_{17})(d_{14} + d_{17}) +  (\gamma_{15} + \gamma_{18})(d_{15} + d_{18}) +  (\gamma_{16} + \gamma_{19})(d_{16} + d_{19})\equiv_{\omega} 0,\\
\sigma_3(f) + f&\equiv_{\omega} (\gamma_1 + \gamma_{2})(d_1 + d_{2}) +  (\gamma_6 + \gamma_{8})(d_6 + d_{8}) +  (\gamma_7 + \gamma_{9})(d_7 + d_{9}) \\
\medskip
&\quad +  (\gamma_{12} + \gamma_{14})(d_{12} + d_{14}) +  (\gamma_{13} + \gamma_{15})(d_{13} + d_{15}) +  (\gamma_{19} + \gamma_{20})(d_{19} + d_{20})\equiv_{\omega} 0,\\
\sigma_4(f) + f&\equiv_{\omega} (\gamma_2 + \gamma_{3})(d_2 + d_{3}) +  (\gamma_5 + \gamma_{6})(d_5 + d_{6}) +  (\gamma_9 + \gamma_{10})(d_9 + d_{10})\\
\medskip
&\quad +  (\gamma_{11} + \gamma_{12})(d_{11} + d_{12}) +  (\gamma_{15} + \gamma_{16})(d_{15} + d_{16}) +  (\gamma_{18} + \gamma_{19})(d_{18} + d_{19})\equiv_{\omega} 0,\\
\sigma_5(f) + f&\equiv_{\omega} (\gamma_3 + \gamma_{4})(d_3 + d_{4}) +  (\gamma_6 + \gamma_{7})(d_6 + d_{7}) +  (\gamma_8 + \gamma_{9})(d_8 + d_{9})\\
\medskip
&\quad +  (\gamma_{12} + \gamma_{13})(d_{12} + d_{13}) +  (\gamma_{14} + \gamma_{15})(d_{14} + d_{15}) +  (\gamma_{17} + \gamma_{18})(d_{17} + d_{18}) \equiv_{\omega} 0,\\
\medskip
\sigma_1(g) + g &\equiv_{\omega} \beta_{205}d_{205} +  (\beta_{205} + \beta_{208} + \beta_{209})(d_{205} + d_{208} + d_{209})\\
&\quad + (\beta_{206} + \beta_{208} + \beta_{210})(d_{206} + d_{208} + d_{210})\\
\medskip
&\quad + (\beta_{207} + \beta_{209} + \beta_{210})(d_{207} + d_{209} + d_{210})  \equiv_{\omega} 0,\\
\sigma_2(g) + g &\equiv_{\omega} (\beta_{202} + \beta_{205})(d_{202} + d_{205}) + (\beta_{203} + \beta_{206})(d_{203} + d_{206})\\
\medskip
&\quad + (\beta_{204} + \beta_{207})(d_{204} + d_{207})  \equiv_{\omega} 0,\\
\sigma_3(g) + g &\equiv_{\omega} (\beta_{201} + \beta_{202})(d_{201} + d_{202}) + (\beta_{206} + \beta_{208})(d_{206} + d_{208})\\
\medskip
&\quad + (\beta_{207} + \beta_{209})(d_{207} + d_{209})\equiv_{\omega} 0,\\
\sigma_4(g) + g &\equiv_{\omega} (\beta_{202} + \beta_{203})(d_{202} + d_{203}) +  (\beta_{205} + \beta_{206})(d_{205} + d_{206})\\
\medskip
&\quad + (\beta_{209} + \beta_{210})(d_{209} + d_{210})\equiv_{\omega} 0,\\
\sigma_5(g) + g &\equiv_{\omega} (\beta_{203} + \beta_{204})(d_{203} + d_{204}) + (\beta_{206} + \beta_{207})(d_{206} + d_{207})\\
&\quad + (\beta_{208} + \beta_{209})(d_{208} + d_{209}) \equiv_{\omega} 0.
\end{align*}

The above equalities imply that $\gamma_{j} = \gamma_{1}$ for all $j,\, 2\leq j\leq 20,$ and $\beta_{201} = \beta_{202} = \cdots = \beta_{210} = 0.$ Next, we compute the action of the general linear group $GL_6$ on $QP_{n_1}^{\otimes 6}(\omega).$ Since $\Sigma_6\subset GL_6,$ if the equivalence class  $[h]_{\omega}$ belongs to the invariant $[QP_{n_1}^{\otimes 6}(\omega)]^{GL_6},$ then $$h\equiv_{\omega} \xi_1\widehat{p(d_1)} + \xi_2\widehat{p(d_{21})} + \xi_3\widehat{p(d_{111})},\ \ \xi_j\in \mathbb F_2,\ j  = 1,\, 2,\, 3.$$ 
Using the homomorphism $\sigma_6: P^{\otimes 6}\longrightarrow P^{\otimes 6}$ and the relation $\sigma_6(h) + h\equiv_{\omega} 0,$ we get 
$$ \sigma_6(h) + h\equiv_{\omega} \big(\xi_1\big(\sum_{5\leq j\leq 10}d_j\big) +   \xi_2d_{27} + \xi_3d_{111} + \mbox{other terms }\big)\equiv_{\omega} 0.$$
The above equality indicates that $\xi_1 = \xi_2 = \xi_3 = 0,$ and therefore  $[QP_{n_1}^{\otimes 6}(\omega)]^{GL_6}$ is zero.
\medskip

Thus, from the above calculations, one gets $[{\rm Ker}((\widetilde {Sq^0_*})_{n_1})]^{GL_6} = 0$. On the other hand, because $$\dim [QP_{n_1}^{\otimes 6}]^{GL_6}\leq \dim [{\rm Ker}((\widetilde {Sq^0_*})_{n_1})]^{GL_6} + \dim [QP_{n_0}^{\otimes 6}]^{GL_6}$$
and $[QP_{n_0}^{\otimes 6}]^{GL_6} \cong (\mathbb F_2 \otimes_{GL_6} {\rm Ann}_{\overline{\mathcal A}}[P^{\otimes 6}]^{*})_{n_0} =  0$ (see \cite{Phuc11}), one gets $[QP_{n_1}^{\otimes 6}]^{GL_6} = 0.$ By this and Corollary \ref{hqs22}), the coinvariant $(\mathbb F_2 \otimes_{GL_6} {\rm Ann}_{\overline{\mathcal A}}[P^{\otimes 6}]^{*})_{n_s}$ vanishes for every positive integer $s.$ Now, it is well-known (see, for instance, Tangora \cite{Tangora}) that the only elements $h_2^{2}g_1 = h_4Ph_2$, and $D_2$ are non-zero in ${\rm Ext}_{\mathcal A}^{6, 6+n_1}(\mathbb F_2, \mathbb F_2)$, and ${\rm Ext}_{\mathcal A}^{6, 6+n_2}(\mathbb F_2, \mathbb F_2),$ respectively. These data and the above calculations indicate that the sixth algebraic transfer, $Tr_6^{\mathcal A}$ is a monomorphism, but not an epimorphism in degrees $n_s,$ for $0 < s\leq 2.$ Therefore, Singer's transfer $Tr_6^{\mathcal A}$ does not detect the non-zero elements $h_4Ph_2$ and $D_2.$ One can observe that the result for $s=3$ is a consequence of the fact (as referenced in \cite{Bruner, Bruner2, Lin2}) that the sixth cohomology group ${\rm Ext}_{\mathcal A}^{6,6+n_3}(\mathbb F_2, \mathbb F_2)$ is trivial. The proof of the theorem is complete.

\section{Conclusion}\label{s5}

The central emphasis of our work is to further investigate the hit problem for the polynomial algebra $P^{\otimes 6} = \mathbb F_2[t_1, \ldots, t_6]$ in degree $n_s = 16\cdot 2^{s}-6$ for general $s$, building on the results from a previous work in \cite{MKR} which addressed the case $s = 0$.  We have shown that the cohit $\mathbb F_2$-module $(P^{\otimes h}/\overline{\mathcal A}P^{\otimes h}){n_{h, s}}$ is equal to the order of the factor group $GL_{h-1}/B_{h-1}$ in general degrees $n_{h, s}=2^{s+4}-h$ with $h\geq 6$ and $s\geq h-5$.
Also, based on a previously established result in \cite{Hai}, we have indicated that the cohit $\mathbb F_q$-module $(\mathbb F_q[t_1,\ldots, t_h]/\overline{\pmb A}_q\mathbb F_q[t_1,\ldots, t_h])_{q^{h-1}-h}$ is equal to the order of the factor group $GL_{h-1}(\mathbb F_q)/B^{*}_{h-1}(\mathbb F_q).$ 
As applications, we have established the dimension result for the cohit module $(\mathbb F_2\otimes_{\mathcal A}P^{\otimes 7})_{n_{s+5}}$ and confirmed Singer's Conjecture \ref{gtSinger} for bidegrees $(h, h+n)$ with $h\geq 1$ and $1\leq n\leq n_0$, as well as for bidegrees $(6, 6+n_s)$ with $s\geq 1$. One of the important corollaries then states that the sixth algebraic transfer does not detect the non-zero elements $h_2^{2}g_1 = h_4Ph_2\in {\rm Ext}_{\mathcal A}^{6, 6+n_1}(\mathbb F_2, \mathbb F_2)$ and $D_2\in {\rm Ext}_{\mathcal A}^{6, 6+n_2}(\mathbb F_2, \mathbb F_2).$
The research conducted in this article advances the understanding of the Peterson hit problem and Singer's algebraic transfer, representing a significant contribution to the literature on this topic. In particular, the application of the hit problem technique in this study has showcased its effectiveness as a powerful tool for investigating the algebraic transfer. We therefore can expect that this approach will lead to more significant breakthroughs in the future.

\section{Appendix}\label{s6}


In this section, we enumerate all admissible monomials in the spaces $P^{\otimes 6}_{n_1}(\omega),$ where
$$ \omega\in \{(4,3,2,1),\ \ (4,3,4), \ \ (4,5,1,1), \ \ (4,5,3)\}.$$ 
We also provide all $\Sigma_6$-invariants of $ (QP_{n_1}^{\otimes 6})(4,3,2,1)$ and $ (QP_{n_1}^{\otimes 6})(4,3,4).$ It is worth noting that all results were verified using the \textbf{OSCAR} computer algebra system \cite{Oscar}.

\section*{A. All admissible monomials in the spaces \mbox{\boldmath $(P^{\otimes 6}_{n_1})(\omega)$}}

\subsection{Admissible monomials in \mbox{\boldmath $(P^{\otimes 6}_{n_1})(4,3,2,1)$}}\label{s50}

We have $(QP_{n_1}^{\otimes 6})(4,3,2,1) = (QP_{n_1}^{\otimes 6})^{0}(4,3,2,1)\oplus (QP_{n_1}^{\otimes 6})^{>0}(4,3,2,1),$ with $\dim (QP_{n_1}^{\otimes 6})^{0}(4,3,2,1) = 5184,$ and $\dim (QP_{n_1}^{\otimes 6})^{>0}(4,3,2,1) = 2880.$   

$\bullet$ The set $\mathscr (C_{n_1}^{\otimes 6})^{0}(4,3,2,1)$ consists of the following 5184 admissible monomials ${\rm a}_j$:

\begin{longtable}{@{}L L L L@{}}
\({\rm a}_{57}=t_{1}\,\allowbreak t_{2}\,\allowbreak t_{3}^{15}\,\allowbreak t_{4}^{2}\,\allowbreak t_{5}^{7}\) & \({\rm a}_{58}=t_{1}\,\allowbreak t_{2}\,\allowbreak t_{3}^{15}\,\allowbreak t_{4}^{2}\,\allowbreak t_{6}^{7}\) & \({\rm a}_{67}=t_{1}\,\allowbreak t_{2}\,\allowbreak t_{3}^{15}\,\allowbreak t_{4}^{7}\,\allowbreak t_{5}^{2}\) & \({\rm a}_{68}=t_{1}\,\allowbreak t_{2}\,\allowbreak t_{3}^{15}\,\allowbreak t_{4}^{7}\,\allowbreak t_{6}^{2}\) \\
\({\rm a}_{69}=t_{1}\,\allowbreak t_{2}\,\allowbreak t_{3}^{15}\,\allowbreak t_{5}^{2}\,\allowbreak t_{6}^{7}\) & \({\rm a}_{72}=t_{1}\,\allowbreak t_{2}\,\allowbreak t_{3}^{15}\,\allowbreak t_{5}^{7}\,\allowbreak t_{6}^{2}\) & \({\rm a}_{93}=t_{1}\,\allowbreak t_{2}\,\allowbreak t_{3}^{2}\,\allowbreak t_{4}^{15}\,\allowbreak t_{5}^{7}\) & \({\rm a}_{94}=t_{1}\,\allowbreak t_{2}\,\allowbreak t_{3}^{2}\,\allowbreak t_{4}^{15}\,\allowbreak t_{6}^{7}\) \\
\({\rm a}_{133}=t_{1}\,\allowbreak t_{2}\,\allowbreak t_{3}^{2}\,\allowbreak t_{4}^{7}\,\allowbreak t_{5}^{15}\) & \({\rm a}_{142}=t_{1}\,\allowbreak t_{2}\,\allowbreak t_{3}^{2}\,\allowbreak t_{4}^{7}\,\allowbreak t_{6}^{15}\) & \({\rm a}_{143}=t_{1}\,\allowbreak t_{2}\,\allowbreak t_{3}^{2}\,\allowbreak t_{5}^{15}\,\allowbreak t_{6}^{7}\) & \({\rm a}_{144}=t_{1}\,\allowbreak t_{2}\,\allowbreak t_{3}^{2}\,\allowbreak t_{5}^{7}\,\allowbreak t_{6}^{15}\) \\
\({\rm a}_{309}=t_{1}\,\allowbreak t_{2}\,\allowbreak t_{3}^{7}\,\allowbreak t_{4}^{15}\,\allowbreak t_{5}^{2}\) & \({\rm a}_{310}=t_{1}\,\allowbreak t_{2}\,\allowbreak t_{3}^{7}\,\allowbreak t_{4}^{15}\,\allowbreak t_{6}^{2}\) & \({\rm a}_{317}=t_{1}\,\allowbreak t_{2}\,\allowbreak t_{3}^{7}\,\allowbreak t_{4}^{2}\,\allowbreak t_{5}^{15}\) & \({\rm a}_{326}=t_{1}\,\allowbreak t_{2}\,\allowbreak t_{3}^{7}\,\allowbreak t_{4}^{2}\,\allowbreak t_{6}^{15}\) \\
\end{longtable}

\medskip

$\bullet$ The set $\mathscr (C_{n_1}^{\otimes 6})^{>0}(4,3,2,1)$ consists of the following 2880 admissible monomials ${\rm a}_j$:

\begin{longtable}{@{}L L L L@{}}
\({\rm a}_{5}=t_{1}\,\allowbreak t_{2}\,\allowbreak t_{3}\,\allowbreak t_{4}^{15}\,\allowbreak t_{5}^{2}\,\allowbreak t_{6}^{6}\) & \({\rm a}_{6}=t_{1}\,\allowbreak t_{2}\,\allowbreak t_{3}\,\allowbreak t_{4}^{15}\,\allowbreak t_{5}^{6}\,\allowbreak t_{6}^{2}\) & \({\rm a}_{8}=t_{1}\,\allowbreak t_{2}\,\allowbreak t_{3}\,\allowbreak t_{4}^{2}\,\allowbreak t_{5}^{15}\,\allowbreak t_{6}^{6}\) & \({\rm a}_{9}=t_{1}\,\allowbreak t_{2}\,\allowbreak t_{3}\,\allowbreak t_{4}^{2}\,\allowbreak t_{5}^{6}\,\allowbreak t_{6}^{15}\) \\
\({\rm a}_{16}=t_{1}\,\allowbreak t_{2}\,\allowbreak t_{3}\,\allowbreak t_{4}^{6}\,\allowbreak t_{5}^{15}\,\allowbreak t_{6}^{2}\) & \({\rm a}_{17}=t_{1}\,\allowbreak t_{2}\,\allowbreak t_{3}\,\allowbreak t_{4}^{6}\,\allowbreak t_{5}^{2}\,\allowbreak t_{6}^{15}\) & \({\rm a}_{49}=t_{1}\,\allowbreak t_{2}\,\allowbreak t_{3}^{15}\,\allowbreak t_{4}\,\allowbreak t_{5}^{2}\,\allowbreak t_{6}^{6}\) & \({\rm a}_{50}=t_{1}\,\allowbreak t_{2}\,\allowbreak t_{3}^{15}\,\allowbreak t_{4}\,\allowbreak t_{5}^{6}\,\allowbreak t_{6}^{2}\) \\
\({\rm a}_{51}=t_{1}\,\allowbreak t_{2}\,\allowbreak t_{3}^{15}\,\allowbreak t_{4}^{2}\,\allowbreak t_{5}\,\allowbreak t_{6}^{6}\) & \({\rm a}_{56}=t_{1}\,\allowbreak t_{2}\,\allowbreak t_{3}^{15}\,\allowbreak t_{4}^{2}\,\allowbreak t_{5}^{6}\,\allowbreak t_{6}\) & \({\rm a}_{63}=t_{1}\,\allowbreak t_{2}\,\allowbreak t_{3}^{15}\,\allowbreak t_{4}^{6}\,\allowbreak t_{5}\,\allowbreak t_{6}^{2}\) & \({\rm a}_{64}=t_{1}\,\allowbreak t_{2}\,\allowbreak t_{3}^{15}\,\allowbreak t_{4}^{6}\,\allowbreak t_{5}^{2}\,\allowbreak t_{6}\) \\
\end{longtable}

\subsection{Admissible monomials in \mbox{\boldmath $(P^{\otimes 6}_{n_1})(4,3,4) = (P^{\otimes 6}_{n_1})^{> 0}(4,3,4)$}}\label{s51-0}

We have $(QP_{n_1}^{\otimes 6})(4,3,4) = (QP_{n_1}^{\otimes 6})^{>0}(4,3,4),$ with $\dim (QP_{n_1}^{\otimes 6})^{>0}(4,3,4) = 210.$ Consequently, the set $\mathscr (C_{n_1}^{\otimes 6})^{>0}(4,3,4)$ consists of the following 210 admissible monomials ${\rm b}_j$:

\begin{longtable}{@{}L L L L@{}}
\({\rm b}_{1}=t_{1}\,\allowbreak t_{2}^{2}\,\allowbreak t_{3}^{4}\,\allowbreak t_{4}^{5}\,\allowbreak t_{5}^{7}\,\allowbreak t_{6}^{7}\) & \({\rm b}_{2}=t_{1}\,\allowbreak t_{2}^{2}\,\allowbreak t_{3}^{4}\,\allowbreak t_{4}^{7}\,\allowbreak t_{5}^{5}\,\allowbreak t_{6}^{7}\) & \({\rm b}_{3}=t_{1}\,\allowbreak t_{2}^{2}\,\allowbreak t_{3}^{4}\,\allowbreak t_{4}^{7}\,\allowbreak t_{5}^{7}\,\allowbreak t_{6}^{5}\) & \({\rm b}_{4}=t_{1}\,\allowbreak t_{2}^{2}\,\allowbreak t_{3}^{5}\,\allowbreak t_{4}^{4}\,\allowbreak t_{5}^{7}\,\allowbreak t_{6}^{7}\) \\
\({\rm b}_{9}=t_{1}\,\allowbreak t_{2}^{2}\,\allowbreak t_{3}^{5}\,\allowbreak t_{4}^{7}\,\allowbreak t_{5}^{4}\,\allowbreak t_{6}^{7}\) & \({\rm b}_{12}=t_{1}\,\allowbreak t_{2}^{2}\,\allowbreak t_{3}^{5}\,\allowbreak t_{4}^{7}\,\allowbreak t_{5}^{7}\,\allowbreak t_{6}^{4}\) & \({\rm b}_{13}=t_{1}\,\allowbreak t_{2}^{2}\,\allowbreak t_{3}^{7}\,\allowbreak t_{4}^{4}\,\allowbreak t_{5}^{5}\,\allowbreak t_{6}^{7}\) & \({\rm b}_{14}=t_{1}\,\allowbreak t_{2}^{2}\,\allowbreak t_{3}^{7}\,\allowbreak t_{4}^{4}\,\allowbreak t_{5}^{7}\,\allowbreak t_{6}^{5}\) \\
\({\rm b}_{15}=t_{1}\,\allowbreak t_{2}^{2}\,\allowbreak t_{3}^{7}\,\allowbreak t_{4}^{5}\,\allowbreak t_{5}^{4}\,\allowbreak t_{6}^{7}\) & \({\rm b}_{18}=t_{1}\,\allowbreak t_{2}^{2}\,\allowbreak t_{3}^{7}\,\allowbreak t_{4}^{5}\,\allowbreak t_{5}^{7}\,\allowbreak t_{6}^{4}\) & \({\rm b}_{19}=t_{1}\,\allowbreak t_{2}^{2}\,\allowbreak t_{3}^{7}\,\allowbreak t_{4}^{7}\,\allowbreak t_{5}^{4}\,\allowbreak t_{6}^{5}\) & \({\rm b}_{20}=t_{1}\,\allowbreak t_{2}^{2}\,\allowbreak t_{3}^{7}\,\allowbreak t_{4}^{7}\,\allowbreak t_{5}^{5}\,\allowbreak t_{6}^{4}\) \\
\end{longtable}

\subsection{Admissible monomials in \mbox{\boldmath $(P^{\otimes 6}_{n_1})(4,5,1,1) = (P^{\otimes 6}_{n_1})^{> 0}(4,5,1,1)$}}\label{s51}

According to the proof of Theorem \ref{dlc1},  the dimension of $(QP^{\otimes 6}_{n_1})^{>0}(4,5,1,1)$ is equal to the cardinality of $D_1,$ where $D_1 = \big\{c_j:\ 1\leq j\leq 336\big\}$ and the admissible monomials $c_j$ are given as follows:

\begin{center}
\begin{tabular}{llll}
${\rm c}_{1}=t_1^{3}t_2t_3^{15}t_4^{2}t_5^{2}t_6^{3}$, & ${\rm c}_{2}=t_1^{3}t_2t_3^{15}t_4^{2}t_5^{3}t_6^{2}$, & ${\rm c}_{3}=t_1^{3}t_2t_3^{15}t_4^{3}t_5^{2}t_6^{2}$, & ${\rm c}_{4}=t_1^{3}t_2^{3}t_3^{15}t_4t_5^{2}t_6^{2}$, \\
${\rm c}_{5}=t_1t_2^{3}t_3^{15}t_4^{2}t_5^{2}t_6^{3}$, & ${\rm c}_{6}=t_1t_2^{3}t_3^{15}t_4^{2}t_5^{3}t_6^{2}$, & ${\rm c}_{7}=t_1t_2^{3}t_3^{15}t_4^{3}t_5^{2}t_6^{2}$, & ${\rm c}_{8}=t_1t_2^{2}t_3^{15}t_4^{2}t_5^{3}t_6^{3}$, \\
${\rm c}_{9}=t_1t_2^{2}t_3^{15}t_4^{3}t_5^{2}t_6^{3}$, & ${\rm c}_{10}=t_1t_2^{2}t_3^{15}t_4^{3}t_5^{3}t_6^{2}$, & ${\rm c}_{11}=t_1^{3}t_2t_3^{2}t_4^{2}t_5^{3}t_6^{15}$, & ${\rm c}_{12}=t_1^{3}t_2t_3^{2}t_4^{2}t_5^{15}t_6^{3}$, \\
${\rm c}_{13}=t_1^{3}t_2t_3^{2}t_4^{3}t_5^{2}t_6^{15}$, & ${\rm c}_{14}=t_1^{3}t_2t_3^{2}t_4^{3}t_5^{15}t_6^{2}$, & ${\rm c}_{15}=t_1^{3}t_2t_3^{2}t_4^{15}t_5^{2}t_6^{3}$, & ${\rm c}_{16}=t_1^{3}t_2t_3^{2}t_4^{15}t_5^{3}t_6^{2}$, \\
${\rm c}_{17}=t_1^{3}t_2t_3^{3}t_4^{2}t_5^{2}t_6^{15}$, & ${\rm c}_{18}=t_1^{3}t_2t_3^{3}t_4^{2}t_5^{15}t_6^{2}$, & ${\rm c}_{19}=t_1^{3}t_2t_3^{3}t_4^{15}t_5^{2}t_6^{2}$, & ${\rm c}_{20}=t_1^{3}t_2^{3}t_3t_4^{2}t_5^{2}t_6^{15}$, \\
${\rm c}_{21}=t_1^{3}t_2^{3}t_3t_4^{2}t_5^{15}t_6^{2}$, & ${\rm c}_{22}=t_1^{3}t_2^{3}t_3t_4^{15}t_5^{2}t_6^{2}$, & ${\rm c}_{23}=t_1^{3}t_2^{15}t_3t_4^{2}t_5^{2}t_6^{3}$, & ${\rm c}_{24}=t_1^{3}t_2^{15}t_3t_4^{2}t_5^{3}t_6^{2}$, \\
${\rm c}_{25}=t_1^{3}t_2^{15}t_3t_4^{3}t_5^{2}t_6^{2}$, & ${\rm c}_{26}=t_1^{3}t_2^{15}t_3^{3}t_4t_5^{2}t_6^{2}$, & ${\rm c}_{27}=t_1^{15}t_2t_3^{2}t_4^{2}t_5^{3}t_6^{3}$, & ${\rm c}_{28}=t_1^{15}t_2t_3^{2}t_4^{3}t_5^{2}t_6^{3}$, \\
${\rm c}_{29}=t_1^{15}t_2t_3^{2}t_4^{3}t_5^{3}t_6^{2}$, & ${\rm c}_{30}=t_1^{15}t_2t_3^{3}t_4^{2}t_5^{2}t_6^{3}$, & ${\rm c}_{31}=t_1^{15}t_2t_3^{3}t_4^{2}t_5^{3}t_6^{2}$, & ${\rm c}_{32}=t_1^{15}t_2t_3^{3}t_4^{3}t_5^{2}t_6^{2}$, \\
\end{tabular}
\end{center}

\newpage
\begin{center}
\begin{tabular}{llll}
${\rm c}_{157}=t_1^{3}t_2^{3}t_3^{13}t_4^{2}t_5^{3}t_6^{2}$, & ${\rm c}_{158}=t_1^{3}t_2^{3}t_3^{13}t_4^{3}t_5^{2}t_6^{2}$, & ${\rm c}_{159}=t_1^{3}t_2^{3}t_3^{3}t_4^{13}t_5^{2}t_6^{2}$, & ${\rm c}_{160}=t_1^{3}t_2^{13}t_3^{3}t_4^{2}t_5^{2}t_6^{3}$, \\
${\rm c}_{161}=t_1^{3}t_2^{13}t_3^{3}t_4^{2}t_5^{3}t_6^{2}$, & ${\rm c}_{162}=t_1^{3}t_2^{13}t_3^{3}t_4^{3}t_5^{2}t_6^{2}$, & ${\rm c}_{163}=t_1^{3}t_2^{13}t_3^{2}t_4^{2}t_5^{3}t_6^{3}$, & ${\rm c}_{164}=t_1^{3}t_2^{13}t_3^{2}t_4^{3}t_5^{2}t_6^{3}$, \\
${\rm c}_{165}=t_1^{3}t_2^{13}t_3^{2}t_4^{3}t_5^{3}t_6^{2}$, & ${\rm c}_{166}=t_1^{3}t_2t_3^{2}t_4^{3}t_5^{6}t_6^{11}$, & ${\rm c}_{167}=t_1^{3}t_2t_3^{3}t_4^{2}t_5^{6}t_6^{11}$, & ${\rm c}_{168}=t_1^{3}t_2t_3^{3}t_4^{6}t_5^{2}t_6^{11}$, \\
${\rm c}_{169}=t_1^{3}t_2t_3^{3}t_4^{6}t_5^{11}t_6^{2}$, & ${\rm c}_{170}=t_1^{3}t_2^{3}t_3t_4^{2}t_5^{6}t_6^{11}$, & ${\rm c}_{171}=t_1^{3}t_2^{3}t_3t_4^{6}t_5^{2}t_6^{11}$, & ${\rm c}_{172}=t_1^{3}t_2^{3}t_3t_4^{6}t_5^{11}t_6^{2}$, \\
${\rm c}_{173}=t_1t_2^{3}t_3^{2}t_4^{3}t_5^{6}t_6^{11}$, & ${\rm c}_{174}=t_1t_2^{3}t_3^{3}t_4^{2}t_5^{6}t_6^{11}$, & ${\rm c}_{175}=t_1t_2^{3}t_3^{3}t_4^{6}t_5^{2}t_6^{11}$, & ${\rm c}_{176}=t_1t_2^{3}t_3^{3}t_4^{6}t_5^{11}t_6^{2}$, \\
${\rm c}_{177}=t_1t_2^{2}t_3^{3}t_4^{3}t_5^{6}t_6^{11}$, & ${\rm c}_{178}=t_1t_2^{2}t_3^{3}t_4^{6}t_5^{11}t_6^{3}$, & ${\rm c}_{179}=t_1t_2^{3}t_3^{2}t_4^{6}t_5^{11}t_6^{3}$, & ${\rm c}_{180}=t_1t_2^{3}t_3^{6}t_4^{2}t_5^{11}t_6^{3}$, \\
${\rm c}_{181}=t_1t_2^{3}t_3^{6}t_4^{11}t_5^{2}t_6^{3}$, & ${\rm c}_{182}=t_1^{3}t_2t_3^{2}t_4^{6}t_5^{11}t_6^{3}$, & ${\rm c}_{183}=t_1^{3}t_2t_3^{6}t_4^{2}t_5^{11}t_6^{3}$, & ${\rm c}_{184}=t_1^{3}t_2t_3^{6}t_4^{11}t_5^{2}t_6^{3}$, \\
${\rm c}_{185}=t_1t_2^{2}t_3^{3}t_4^{6}t_5^{3}t_6^{11}$, & ${\rm c}_{186}=t_1t_2^{3}t_3^{2}t_4^{6}t_5^{3}t_6^{11}$, & ${\rm c}_{187}=t_1t_2^{3}t_3^{6}t_4^{2}t_5^{3}t_6^{11}$, & ${\rm c}_{188}=t_1t_2^{3}t_3^{6}t_4^{11}t_5^{3}t_6^{2}$, \\
\end{tabular}%
\end{center}

\newpage
\subsection{Admissible monomials in \mbox{\boldmath $(P^{\otimes 6}_{n_1})(4,5,3) = (P^{\otimes 6}_{n_1})^{> 0}(4,5,3)$}}\label{s52}

By the proof of Theorem \ref{dlc1}, $\dim (QP^{\otimes 6}_{n_1})^{>0}(4,5,3) = |D_2| =  210,$ where $D_2 = \big\{d_j:\ 1\leq j\leq 210\big\}$ and the admissible monomials $d_j$ are determined as follows:

\begin{center}
\begin{tabular}{llrr}
${\rm d}_{1}=t_1t_2^{2}t_3^{2}t_4^{7}t_5^{7}t_6^{7}$, & ${\rm d}_{2}=t_1t_2^{2}t_3^{7}t_4^{2}t_5^{7}t_6^{7}$, & \multicolumn{1}{l}{${\rm d}_{3}=t_1t_2^{2}t_3^{7}t_4^{7}t_5^{2}t_6^{7}$,} & \multicolumn{1}{l}{${\rm d}_{4}=t_1t_2^{2}t_3^{7}t_4^{7}t_5^{7}t_6^{2}$,} \\
${\rm d}_{5}=t_1t_2^{7}t_3^{2}t_4^{2}t_5^{7}t_6^{7}$, & ${\rm d}_{6}=t_1t_2^{7}t_3^{2}t_4^{7}t_5^{2}t_6^{7}$, & \multicolumn{1}{l}{${\rm d}_{7}=t_1t_2^{7}t_3^{2}t_4^{7}t_5^{7}t_6^{2}$,} & \multicolumn{1}{l}{${\rm d}_{8}=t_1t_2^{7}t_3^{7}t_4^{2}t_5^{2}t_6^{7}$,} \\
${\rm d}_{9}=t_1t_2^{7}t_3^{7}t_4^{2}t_5^{7}t_6^{2}$, & ${\rm d}_{10}=t_1t_2^{7}t_3^{7}t_4^{7}t_5^{2}t_6^{2}$, & \multicolumn{1}{l}{${\rm d}_{11}=t_1^{7}t_2t_3^{2}t_4^{2}t_5^{7}t_6^{7}$,} & \multicolumn{1}{l}{${\rm d}_{12}=t_1^{7}t_2t_3^{2}t_4^{7}t_5^{2}t_6^{7}$,} \\
${\rm d}_{13}=t_1^{7}t_2t_3^{2}t_4^{7}t_5^{7}t_6^{2}$, & ${\rm d}_{14}=t_1^{7}t_2t_3^{7}t_4^{2}t_5^{2}t_6^{7}$, & \multicolumn{1}{l}{${\rm d}_{15}=t_1^{7}t_2t_3^{7}t_4^{2}t_5^{7}t_6^{2}$,} & \multicolumn{1}{l}{${\rm d}_{16}=t_1^{7}t_2t_3^{7}t_4^{7}t_5^{2}t_6^{2}$,} \\
${\rm d}_{17}=t_1^{7}t_2^{7}t_3t_4^{2}t_5^{2}t_6^{7}$, & ${\rm d}_{18}=t_1^{7}t_2^{7}t_3t_4^{2}t_5^{7}t_6^{2}$, & \multicolumn{1}{l}{${\rm d}_{19}=t_1^{7}t_2^{7}t_3t_4^{7}t_5^{2}t_6^{2}$,} & \multicolumn{1}{l}{${\rm d}_{20}=t_1^{7}t_2^{7}t_3^{7}t_4t_5^{2}t_6^{2}$,} \\
${\rm d}_{21}=t_1t_2^{2}t_3^{3}t_4^{6}t_5^{7}t_6^{7}$, & ${\rm d}_{22}=t_1t_2^{2}t_3^{3}t_4^{7}t_5^{6}t_6^{7}$, & \multicolumn{1}{l}{${\rm d}_{23}=t_1t_2^{2}t_3^{3}t_4^{7}t_5^{7}t_6^{6}$,} & \multicolumn{1}{l}{${\rm d}_{24}=t_1t_2^{2}t_3^{7}t_4^{3}t_5^{6}t_6^{7}$,} \\
${\rm d}_{25}=t_1t_2^{2}t_3^{7}t_4^{3}t_5^{7}t_6^{6}$, & ${\rm d}_{26}=t_1t_2^{2}t_3^{7}t_4^{7}t_5^{3}t_6^{6}$, & \multicolumn{1}{l}{${\rm d}_{27}=t_1t_2^{3}t_3^{2}t_4^{6}t_5^{7}t_6^{7}$,} & \multicolumn{1}{l}{${\rm d}_{28}=t_1t_2^{3}t_3^{2}t_4^{7}t_5^{6}t_6^{7}$,} \\
\end{tabular}
\end{center}

\newpage
\begin{center}
\begin{tabular}{llll}
${\rm d}_{161}=t_1^{3}t_2^{3}t_3^{5}t_4^{2}t_5^{6}t_6^{7}$, & ${\rm d}_{162}=t_1^{3}t_2^{3}t_3^{5}t_4^{2}t_5^{7}t_6^{6}$, & \multicolumn{1}{l}{${\rm d}_{163}=t_1^{3}t_2^{3}t_3^{5}t_4^{6}t_5^{2}t_6^{7}$,} & \multicolumn{1}{l}{${\rm d}_{164}=t_1^{3}t_2^{3}t_3^{5}t_4^{6}t_5^{7}t_6^{2}$,} \\
${\rm d}_{165}=t_1^{3}t_2^{3}t_3^{5}t_4^{7}t_5^{2}t_6^{6}$, & ${\rm d}_{166}=t_1^{3}t_2^{3}t_3^{5}t_4^{7}t_5^{6}t_6^{2}$, & \multicolumn{1}{l}{${\rm d}_{167}=t_1^{3}t_2^{5}t_3^{2}t_4^{3}t_5^{6}t_6^{7}$,} & \multicolumn{1}{l}{${\rm d}_{168}=t_1^{3}t_2^{5}t_3^{2}t_4^{3}t_5^{7}t_6^{6}$,} \\
${\rm d}_{169}=t_1^{3}t_2^{5}t_3^{2}t_4^{6}t_5^{3}t_6^{7}$, & ${\rm d}_{170}=t_1^{3}t_2^{5}t_3^{2}t_4^{6}t_5^{7}t_6^{3}$, & \multicolumn{1}{l}{${\rm d}_{171}=t_1^{3}t_2^{5}t_3^{2}t_4^{7}t_5^{3}t_6^{6}$,} & \multicolumn{1}{l}{${\rm d}_{172}=t_1^{3}t_2^{5}t_3^{2}t_4^{7}t_5^{6}t_6^{3}$,} \\
${\rm d}_{173}=t_1^{3}t_2^{5}t_3^{3}t_4^{2}t_5^{6}t_6^{7}$, & ${\rm d}_{174}=t_1^{3}t_2^{5}t_3^{3}t_4^{2}t_5^{7}t_6^{6}$, & \multicolumn{1}{l}{${\rm d}_{175}=t_1^{3}t_2^{5}t_3^{3}t_4^{6}t_5^{2}t_6^{7}$,} & \multicolumn{1}{l}{${\rm d}_{176}=t_1^{3}t_2^{5}t_3^{3}t_4^{6}t_5^{7}t_6^{2}$,} \\
${\rm d}_{177}=t_1^{3}t_2^{5}t_3^{3}t_4^{7}t_5^{2}t_6^{6}$, & ${\rm d}_{178}=t_1^{3}t_2^{5}t_3^{3}t_4^{7}t_5^{6}t_6^{2}$, & \multicolumn{1}{l}{${\rm d}_{179}=t_1^{3}t_2^{5}t_3^{6}t_4^{2}t_5^{3}t_6^{7}$,} & \multicolumn{1}{l}{${\rm d}_{180}=t_1^{3}t_2^{5}t_3^{6}t_4^{2}t_5^{7}t_6^{3}$,} \\
${\rm d}_{181}=t_1^{3}t_2^{5}t_3^{6}t_4^{3}t_5^{2}t_6^{7}$, & ${\rm d}_{182}=t_1^{3}t_2^{5}t_3^{6}t_4^{3}t_5^{7}t_6^{2}$, & \multicolumn{1}{l}{${\rm d}_{183}=t_1^{3}t_2^{5}t_3^{6}t_4^{7}t_5^{2}t_6^{3}$,} & \multicolumn{1}{l}{${\rm d}_{184}=t_1^{3}t_2^{5}t_3^{6}t_4^{7}t_5^{3}t_6^{2}$,} \\
\end{tabular}%
\end{center}

\medskip

\section*{B.  All $\Sigma_6$-invariants of $ (QP_{n_1}^{\otimes 6})(4,3,2,1)$ and $ (QP_{n_1}^{\otimes 6})(4,3,4).$}

\subsection{All $\Sigma_6$-invariants of $ (QP_{n_1}^{\otimes 6})(4,3,2,1)$}

As shown above, $\mathscr (C_{n_1}^{\otimes 6})(4,3,2,1) = \mathscr (C_{n_1}^{\otimes 6})^{0}(4,3,2,1)\cup \mathscr (C_{n_1}^{\otimes 6})^{>0}(4,3,2,1),$ with $|\mathscr (C_{n_1}^{\otimes 6})^{0}(4,3,2,1)|  = 5184$ and $\mathscr (C_{n_1}^{\otimes 6})^{>0}(4,3,2,1) = 2880.$ By direct calculations using these results and the $\mathcal A$-homomorphisms $\sigma_d: P_6\longrightarrow P_6,\, 1\leq d\leq 5,$ we obtain
$$  [(QP_{n_1}^{\otimes 6})(4,3,2,1)]^{\Sigma_6} = \mathbb F_2\cdot ([\widehat{q}_j]:\ 1\leq j\leq 10),$$
where the invariant polynomials $\widehat{q}_j$ are determined as follows:



\subsection{All $\Sigma_6$-invariants of $ (QP_{n_1}^{\otimes 6})(4,3,4)$}

Set $\widetilde{\omega}^{*}:= (4,3,4).$ As shown above, $\mathscr (C_{n_1}^{\otimes 6})(\widetilde{\omega}^{*}) =  \mathscr (C_{n_1}^{\otimes 6})^{>0}(\widetilde{\omega}^{*}),$ with $|\mathscr (C_{n_1}^{\otimes 6})^{>0}(\widetilde{\omega}^{*})|  = 210.$ By direct calculations using this result and the $\mathcal A$-homomorphisms $\sigma_d: P_6\longrightarrow P_6,\, 1\leq d\leq 5,$ we obtain
$$  [(QP_{n_1}^{\otimes 6})(\widetilde{\omega}^{*})]^{\Sigma_6} = \mathbb F_2\cdot [\widehat{q}^{*}]_{\widetilde{\omega}^{*}},$$
where 
\begin{align*}
\widehat{q}^{*} &= t_{1}^{7} t_{2}^{7} t_{3}^{3} t_{4} t_{5}^{4} t_{6}^{4} + t_{1}^{7} t_{2}^{3} t_{3}^{7} t_{4} t_{5}^{4} t_{6}^{4} + t_{1}^{3} t_{2}^{7} t_{3}^{7} t_{4} t_{5}^{4} t_{6}^{4} + t_{1}^{7} t_{2}^{7} t_{3} t_{4}^{3} t_{5}^{4} t_{6}^{4} \\
    &\quad + t_{1}^{7} t_{2} t_{3}^{7} t_{4}^{3} t_{5}^{4} t_{6}^{4} + t_{1} t_{2}^{7} t_{3}^{7} t_{4}^{3} t_{5}^{4} t_{6}^{4} + t_{1}^{7} t_{2}^{3} t_{3} t_{4}^{7} t_{5}^{4} t_{6}^{4} + t_{1}^{3} t_{2}^{7} t_{3} t_{4}^{7} t_{5}^{4} t_{6}^{4} \\
    &\quad + t_{1}^{7} t_{2} t_{3}^{3} t_{4}^{7} t_{5}^{4} t_{6}^{4} + t_{1} t_{2}^{7} t_{3}^{3} t_{4}^{7} t_{5}^{4} t_{6}^{4} + t_{1}^{3} t_{2} t_{3}^{7} t_{4}^{7} t_{5}^{4} t_{6}^{4} + t_{1} t_{2}^{3} t_{3}^{7} t_{4}^{7} t_{5}^{4} t_{6}^{4} \\
    &\quad + t_{1}^{7} t_{2}^{7} t_{3} t_{4}^{2} t_{5}^{5} t_{6}^{4} + t_{1}^{7} t_{2} t_{3}^{7} t_{4}^{2} t_{5}^{5} t_{6}^{4} + t_{1} t_{2}^{7} t_{3}^{7} t_{4}^{2} t_{5}^{5} t_{6}^{4} + t_{1}^{7} t_{2} t_{3}^{2} t_{4}^{7} t_{5}^{5} t_{6}^{4} \\
    &\quad + t_{1} t_{2}^{7} t_{3}^{2} t_{4}^{7} t_{5}^{5} t_{6}^{4} + t_{1} t_{2}^{2} t_{3}^{7} t_{4}^{7} t_{5}^{5} t_{6}^{4} + t_{1}^{7} t_{2}^{3} t_{3} t_{4}^{4} t_{5}^{7} t_{6}^{4} + t_{1}^{3} t_{2}^{7} t_{3} t_{4}^{4} t_{5}^{7} t_{6}^{4} \\
    &\quad + t_{1}^{7} t_{2} t_{3}^{3} t_{4}^{4} t_{5}^{7} t_{6}^{4} + t_{1} t_{2}^{7} t_{3}^{3} t_{4}^{4} t_{5}^{7} t_{6}^{4} + t_{1}^{3} t_{2} t_{3}^{7} t_{4}^{4} t_{5}^{7} t_{6}^{4} + t_{1} t_{2}^{3} t_{3}^{7} t_{4}^{4} t_{5}^{7} t_{6}^{4} \\
    &\quad + t_{1}^{7} t_{2} t_{3}^{2} t_{4}^{5} t_{5}^{7} t_{6}^{4} + t_{1} t_{2}^{7} t_{3}^{2} t_{4}^{5} t_{5}^{7} t_{6}^{4} + t_{1} t_{2}^{2} t_{3}^{7} t_{4}^{5} t_{5}^{7} t_{6}^{4} + t_{1}^{3} t_{2} t_{3}^{4} t_{4}^{7} t_{5}^{7} t_{6}^{4} \\
    &\quad + t_{1} t_{2}^{3} t_{3}^{4} t_{4}^{7} t_{5}^{7} t_{6}^{4} + t_{1} t_{2}^{2} t_{3}^{5} t_{4}^{7} t_{5}^{7} t_{6}^{4} + t_{1}^{7} t_{2}^{7} t_{3} t_{4}^{2} t_{5}^{4} t_{6}^{5} + t_{1}^{7} t_{2} t_{3}^{7} t_{4}^{2} t_{5}^{4} t_{6}^{5} \\
    &\quad + t_{1} t_{2}^{7} t_{3}^{7} t_{4}^{2} t_{5}^{4} t_{6}^{5} + t_{1}^{7} t_{2} t_{3}^{2} t_{4}^{7} t_{5}^{4} t_{6}^{5} + t_{1} t_{2}^{7} t_{3}^{2} t_{4}^{7} t_{5}^{4} t_{6}^{5} + t_{1} t_{2}^{2} t_{3}^{7} t_{4}^{7} t_{5}^{4} t_{6}^{5} \\
    &\quad + t_{1}^{7} t_{2} t_{3}^{2} t_{4}^{4} t_{5}^{7} t_{6}^{5} + t_{1} t_{2}^{7} t_{3}^{2} t_{4}^{4} t_{5}^{7} t_{6}^{5} + t_{1} t_{2}^{2} t_{3}^{7} t_{4}^{4} t_{5}^{7} t_{6}^{5} + t_{1} t_{2}^{2} t_{3}^{4} t_{4}^{7} t_{5}^{7} t_{6}^{5} \\
    &\quad + t_{1}^{7} t_{2}^{3} t_{3} t_{4}^{4} t_{5}^{4} t_{6}^{7} + t_{1}^{3} t_{2}^{7} t_{3} t_{4}^{4} t_{5}^{4} t_{6}^{7} + t_{1}^{7} t_{2} t_{3}^{3} t_{4}^{4} t_{5}^{4} t_{6}^{7} + t_{1} t_{2}^{7} t_{3}^{3} t_{4}^{4} t_{5}^{4} t_{6}^{7} \\
    &\quad + t_{1}^{3} t_{2} t_{3}^{7} t_{4}^{4} t_{5}^{4} t_{6}^{7} + t_{1} t_{2}^{3} t_{3}^{7} t_{4}^{4} t_{5}^{4} t_{6}^{7} + t_{1}^{7} t_{2} t_{3}^{2} t_{4}^{5} t_{5}^{4} t_{6}^{7} + t_{1} t_{2}^{7} t_{3}^{2} t_{4}^{5} t_{5}^{4} t_{6}^{7} \\
    &\quad + t_{1} t_{2}^{2} t_{3}^{7} t_{4}^{5} t_{5}^{4} t_{6}^{7} + t_{1}^{3} t_{2} t_{3}^{4} t_{4}^{7} t_{5}^{4} t_{6}^{7} + t_{1} t_{2}^{3} t_{3}^{4} t_{4}^{7} t_{5}^{4} t_{6}^{7} + t_{1} t_{2}^{2} t_{3}^{5} t_{4}^{7} t_{5}^{4} t_{6}^{7} \\
    &\quad + t_{1}^{7} t_{2} t_{3}^{2} t_{4}^{4} t_{5}^{5} t_{6}^{7} + t_{1} t_{2}^{7} t_{3}^{2} t_{4}^{4} t_{5}^{5} t_{6}^{7} + t_{1} t_{2}^{2} t_{3}^{7} t_{4}^{4} t_{5}^{5} t_{6}^{7} + t_{1} t_{2}^{2} t_{3}^{4} t_{4}^{7} t_{5}^{5} t_{6}^{7} \\
    &\quad + t_{1}^{3} t_{2} t_{3}^{4} t_{4}^{4} t_{5}^{7} t_{6}^{7} + t_{1} t_{2}^{3} t_{3}^{4} t_{4}^{4} t_{5}^{7} t_{6}^{7} + t_{1} t_{2}^{2} t_{3}^{5} t_{4}^{4} t_{5}^{7} t_{6}^{7} + t_{1} t_{2}^{2} t_{3}^{4} t_{4}^{5} t_{5}^{7} t_{6}^{7}.
\end{align*}

\end{document}